\documentclass[10pt]{amsart}
\usepackage{amssymb}
\usepackage{amsbsy}
\usepackage{amscd}
\usepackage{amsmath}
\usepackage{amsthm}

\theoremstyle{plain}
 \newtheorem{thm}{Theorem}[section]
 \newtheorem{lem}[thm]{Lemma}
 \newtheorem{prop}[thm]{Proposition}
 \newtheorem{cor}[thm]{Corollary}
 \newtheorem{guess}[thm]{Conjecture}
\theoremstyle{definition}
 \newtheorem{defn}{Definition}[section]
\theoremstyle{remark}
 \newtheorem{rem}{Remark}[section]
 \newtheorem{ex}{Example}[section]

\def\Bbb{\mathbb}
\def\frak{\mathfrak}
\def\cal{\mathcal}

\newcommand{ \Supp}{\operatorname{Supp}}
\newcommand{\Ext}{\operatorname{Ext}}
\newcommand{\Hom}{\operatorname{Hom}}
\newcommand{\codim}{\operatorname{codim}}
\newcommand{\im}{\operatorname{im}}

\newcommand{\rk}{\operatorname{rk}}

\newcommand{\NS}{\operatorname{NS}}
\newcommand{\coker}{\operatorname{coker}}
\newcommand{\Pic}{\operatorname{Pic}}
\newcommand{\ch}{\operatorname{ch}}
\newcommand{\td}{\operatorname{td}}
\newcommand{\Alb}{\operatorname{Alb}}
\newcommand{\Hilb}{\operatorname{Hilb}}
\newcommand{\Quot}{\operatorname{Quot}}

\newcommand{\Spec}{\operatorname{Spec}}
\newcommand{\Amp}{\operatorname{Amp}}
\newcommand{\Km}{\operatorname{Km}}
\newcommand{\IT}{\operatorname{IT}}
\newcommand{\WIT}{\operatorname{WIT}}
\newcommand{\Spl}{\operatorname{Spl}}
\newcommand{\Proj}{\operatorname{Proj}}

\font\b=cmr10 scaled \magstep5
\def\bigzerou{\smash{\lower1.7ex\hbox{\b 0}}}

\numberwithin{equation}{section}
\begin{document}

\title{
Moduli spaces of stable sheaves on abelian surfaces
}
\author{K\={o}ta Yoshioka}
 
\address{
Department of mathematics, Faculty of Science, Kobe University,
Kobe, 657, Japan}

\email{yoshioka@math.kobe-u.ac.jp}
 \subjclass{14D20}
 \maketitle

\section{Introduction}

Let $X$ be a smooth projective surface defined over ${\Bbb C}$
and $H$ an ample line bundle on $X$.
If $K_X$ is trivial, that is, $X$ is an abelian or a K3 surface,
Mukai \cite{Mu:4} introduced a quite useful notion now called 
Mukai lattice $(H^{ev}(X,{\Bbb Z}),\langle \;\;,\;\;\rangle)$,
where $H^{ev}(X,{\Bbb Z})=\oplus_i H^{2i}(X,{\Bbb Z})$
and $\langle x,y \rangle=-\int_X (x^{\vee} y)$ (see Defn. 1.1). 
$\langle \quad,\quad \rangle$ is an even unimodular bilinear form. 
For a coherent sheaf $E$ on $X$,
we can attach an element of $H^{ev}(X,{\Bbb Z})$
called Mukai vector $v(E):=\ch(E)\sqrt{\td_X}$,
where $\ch(E)$ is the Chern character of $E$ and 
$\td_X$ is the Todd class of $X$.
If $X$ is an abelian surface, then $\td_X=1$, and hence
$v(E)$ is nothing but the Chern character of $E$.
For an element $v \in H^*(X,{\Bbb Z})$,
we denote the $0$-th component $v_0 \in H^0(X,{\Bbb Z})$
by $\rk v$ and the second component $v_1 \in H^2(X,{\Bbb Z})$
by $c_1(v)$.
We set $\ell(v):=\gcd(\rk v,c_1(v)) \in {\Bbb Z}_{\geq 0}$.
Then $v$ is written as 
$v=\ell(v)(r+\xi)+a \omega$,
where $\omega$ is the fundamental class of $X$,
$r \in {\Bbb Z}$, $\xi \in H^2(X,{\Bbb Z})$ and $r+\xi$ is primitive. 
We denote the moduli space of (Gieseker) stable sheaves $E$
of $v(E)=v$ by $M_H(v)$ and its Gieseker compactification by
$\overline{M}_H(v)$.
In \cite{Mu:3}, Mukai proved that $M_H(v)$ is smooth
of dimension $\langle v^2 \rangle+2$.
Moreover, he
constructed a natural symplectic structure on $M_H(v)$.
If $H$ is a general element of the ample cone $\Amp(X)$
(i.e. there are hyperplanes called walls in $\NS(X)$
 and $H$ does not lie on these walls \cite{Y:2})
and $v$ is primitive, then $\overline{M}_H(v)=M_H(v)$,
in particular
$M_H(v)$ is a projective scheme.
Hence under this condition,
$M_H(v)$ is a projective symplectic manifold.
We remark that the number of walls is locally finite.
Thus this condition is not strong.
If $X$ is a K3 surface and $v$ is primitive, then
$M_H(v)$ is extensively studied by many authors.
In particular,   
$M_H(v)$ is an irreducible symplectic
manifold and the period of $M_H(v)$ is written down in
terms of Mukai lattice (\cite{Mu:4}, \cite{Mu:6}, \cite{O:1}, 
\cite{Y:5},\cite{Y:8}).

In this paper, we shall treat mainly the case where $X$ is an abelian surface. 
Let $v=r+c_1+a \omega \in H^{ev}(X,{\Bbb Z})$, $c_1 \in \NS(X)$,
be a Mukai vector,
where we identify $H^0(X,{\Bbb Z})$ with ${\Bbb Z}$.
\begin{defn}
A Mukai vector $v$ is {\it positive} ($v>0$), if
(1) $r>0$, or (2) $r=0$, $c_1$ is effective and $a \ne 0$ or
(3) $r=c_1=0$ and $a < 0$.
\end{defn}
Let $E$ be a stable sheaf of $v(E)=v$ on $X$.
Then obviously $T^*_x(E) \otimes L \not \cong E$ for a general 
$(x, L) \in X \times \Pic^0(X)$.
Thus $\dim M_H(v)=\langle v^2 \rangle+2 >0$. 
Since $M_H(v)$ is a symplectic manifold,
$\dim M_H(v)$ is even, and hence $\langle v^2 \rangle \geq 0$. 
If $\langle v^2 \rangle =0$, then Mukai showed that
$M_H(v)$ is an abelian surface (see \cite[(5.13)]{Mu:6}).

In \cite{Y:3}, we studied $H^i(M_H(v),{\Bbb Z})$
for $i=1,2$ under the assumption $\ell(v)=1$.
We also constructed a morphism ${\frak a}_v:M_H(v) \to X \times \widehat{X}$
and proved that ${\frak a}_v$ is an albanese map
for $\langle v^2 \rangle \geq 2$,
where $\widehat{X}$ is the dual of $X$. 
In this paper, we remove the technical assumptions in \cite{Y:3}
(Theorem \ref{thm:H2}) and compute the deformation type of 
$M_H(v)$. 
\begin{thm}\label{thm:deform equiv}
Let $X$ be an abelian surface.
Let $v$ be a primitive Mukai vector such that
$v>0$, $c_1(v) \in \NS(X)$ and $\langle v^2 \rangle \geq 2$.
Then for an ample divisor $H$ such that $\overline{M}_H(v)=M_H(v)$,
\begin{enumerate}
\item[(1)]
${\frak a}_v:M_H(v) \to X \times \widehat{X}$ is the albanese map.
\item[(2)]
$M_H(v)$ is deformation equivalent to
$\widehat{X} \times \Hilb_X^{\langle v^2 \rangle/2}$.
\end{enumerate}
\end{thm}

Then our next interest is the fiber of the albanese map.
If $\langle v^2 \rangle=0$, then Mukai showed that
 ${\frak a}_v$ is an immersion.
If $\langle v^2 \rangle=2$, then Mukai \cite{Mu:2} and the author
\cite[Prop. 4.2]{Y:3} showed that ${\frak a}_v:M_H(v) \to X \times \widehat{X}$
is an isomorphism.
Hence we assume that $\langle v^2 \rangle \geq 4$. 
Let $K_H(v)$ be a fiber of ${\frak a}_v$.
Then $\dim K_H(v)=\langle v^2 \rangle-2$.
Hence if $\langle v^2 \rangle \geq 6$, then 
$\dim K_H(v) \geq 4$.
In this case, we get the following,
which is an analogous result to that for a K3 surface.

\begin{thm}\label{thm:period}
Let $X$ be an abelian surface.
Let $v\in
H^{ev}(X,{\Bbb Z})$
be a primitive Mukai vector such that $v>0$, $c_1(v) \in \NS(X)$
and $\langle v^2 \rangle \geq 6$.

\begin{enumerate}
\item[(1)]
For a general ample line bundle $H$,
$K_H(v)$ is 
deformation equivalent to a generalized Kummer variety
$K_{\langle v^2 \rangle/2-1}$ constructed by Beauville.
In particular, $K_H(v)$ is 
an irreducible symplectic manifold.
\item[(2)]
Let $B_{K_H(v)}$ be Beauville's bilinear form on
$H^2(K_H(v),{\Bbb Z})$.
Then
\begin{equation}
\theta_v:(v^{\perp}, \langle \;\;,\;\; \rangle) \to 
(H^2(K_H(v),{\Bbb Z}),B_{K_H(v)})
\end{equation}
is an isometry of Hodge structures,
where $\theta_v$ is the composition of
Mukai homomorphism $v^{\perp} \to H^2(M_H(v),{\Bbb Z})$ 
and the restriction map $H^2(M_H(v),{\Bbb Z}) \to H^2(K_H(v),{\Bbb Z})$. 
\end{enumerate}
\end{thm}

Our theorem shows that Mukai lattice for an abelian surface
is as important as that for a K3 surface.
This is our main motivation of this paper.  
As an application of this theorem, we shall show that
for some $v$, $M_H(v)$ is not birationally equivalent
to $\widehat{Y} \times \Hilb_Y^n$ for any $Y$
( Example \ref{ex:ex1}).  

In section 1, we collect some known facts which will be used in this paper.
Also we define chamber structure of polarizations for moduli spaces of
stable sheaves of pure dimension 1.
In section 2, we collect elementary facts on Fourier-Mukai transforms. 
In K3 surfaces cases, we know that isometries of Mukai lattices 
are quite useful to compute the period of moduli spaces
\cite{Mu:4}, \cite{Y:5}, \cite{Y:8}. 
Hence it is also important to study isometries in our cases.
Fourier-Mukai transforms are good examples of isometries
\cite{Mu:5}. So we also consider the relation between Fourier-Mukai functors
and the homomorphism $\theta_v$ (Proposition \ref{prop:comm}
and \ref{prop:comm2}).
Indeed, if we can find enough examples  to prove Theorem \ref{thm:period} (1)
of birational maps $M_H(v) \cdots \to M_{H'}(w)$
induced by Fourier-Mukai functor,
then (2) follows immediately
from these propositions.

Unfortunately, even for original Fourier-Mukai functor ${\cal F}_{\cal P}$
investigated by Mukai \cite{Mu:2}, \cite{Mu:5}
(also see \cite{A:1},\cite{F-L:1}),
only a few example of birational maps induced by ${\cal F}_{\cal P}$
are known.
Hence it is very interesting to construct lots of examples of
birational maps.
In section 3.1, we shall construct enough examples 
(Proposition \ref{lem:fourier4}, \ref{lem:fourier3})
to prove 
Theorem \ref{thm:deform equiv} and \ref{thm:period} 
for the case where $\ell(v)=1$.
Let us briefly explain the main idea of our construction.

Let $(X,H)$ be a polarized abelian surface of 
$\NS(X)={\Bbb Z}H$.
Let $\widehat{X}$ be the dual abelian surface of $X$ and
${\cal P}$ the Poincar\'{e} line bundle on $\widehat{X} \times X$.
We assume that $v(E)=r+H+a \omega$ and
$\chi(E)=r+a<0$.
Then for any $E \in M_H(v)$, we can show that 
\begin{enumerate}
\item[(1)]
${\cal F}_{\cal P}^i(E):=
R^i p_{\widehat{X}*}(p_X^* E \otimes {\cal P})=0$ for $i \ne 1$ and
\item[(2)]
${\cal F}_{\cal P}^1(E)$ is a stable sheaf,
\end{enumerate}
where $p_X:\widehat{X} \times X \to X$ 
(resp. $p_{\widehat{X}}: \widehat{X} \times X \to \widehat{X}$)
is the projection.
We shall explain how to prove (1).
Since $E$ is stable and $(c_1(E),H)>0$,
${\cal F}_{\cal P}^2(E )=0$.
Since $p_X^* E \otimes {\cal P}$ is torsion free,
it is sufficient to prove 
\begin{equation}\label{eq:wit1}
\#\{x \in \widehat{X}|H^0(X,E \otimes {\cal P}_x) \ne 0 \} < \infty.
\end{equation}
Assume that $H^0(X,E \otimes {\cal P}_{x_i}) \ne 0$
for distinct points $x_1,x_2,\dots,x_n \in \widehat{X}$.
Our main idea is to consider the evaluation map:
\begin{equation}
 \phi:\oplus_{i=1}^n {\cal P}_{x_i}^{\vee} \otimes 
 H^0(X,E \otimes {\cal P}_{x_i}) \to E.
\end{equation}
If $n >r$,
then $\phi$ is surjective in codimension 1 and
$\ker \phi$ is stable (see \cite[Lem. 2.1]{Y:5}).
Then we see that
$\langle v(\ker \phi)^2 \rangle<-2$, if $n \gg r$.
On the other hand, $\langle v(F)^2 \rangle \geq 0$
for any stable sheaf $F$ on $\widehat{X}$.
Thus $n$ is bounded above. 
Therefore \eqref{eq:wit1} holds.
Motivated by a recent work of Markman \cite{Mr:1} and
a work of Mukai \cite{Mu:5},
we shall also treat the composition of
${\cal F}_{\cal P}$ and the ``taking dual'' functor.

In order to treat the case where $\ell(v)>1$, we need to consider
Fourier-Mukai transform on elliptic abelian surfaces.
Hence we also consider (relative) Fourier-Mukai transform on elliptic surfaces
(Theorem \ref{thm:FM}).

In section 4, we shall resume moduli spaces of stable sheaves 
on abelian surfaces.
Since the canonical bundle of $M_H(v)$ is trivial,
$M_H(v)$ has a Bogomolov decomposition.
We shall first construct a decomposition which will become
a Bogomolov decomposition for $M_H(v)$.
In section 4.2, we discuss Fourier-Mukai functor on 
abelian surfaces again.
In particular, we consider the relation of Fourier-Mukai functor
and the albanese map of moduli spaces.
In section 4.3, we prove Theorem  
\ref{thm:deform equiv} and 
\ref{thm:period}. 
We shall first treat rank 1 case.
In this case, $K_H(v)$ is the generalized Kummer variety
$K_{n-1}$
constructed by Beauville \cite{B:1},
where $n=\langle v^2 \rangle/2$.
Hence Theorem \ref{thm:period} follows from 
Beauville's description of $H^2(K_{n-1},{\Bbb Q})$ and
some computations.  
Higher rank cases follow from
deformation arguments as in \cite{G-H:1}, \cite{O:1}, \cite{Y:3}
and isomorphisms constructed in section 3.

As another application of Theorem \ref{thm:period},
we shall construct a non-K\"{a}hler compact symplectic manifold
which is an elementary transform of $K_H(v)$ (cf. \cite[2.5]{H:2}).
In section 4.5, we shall treat the remaining case, 
i.e. $\langle v^2 \rangle=4$.
In this case, we shall prove that $K_H(v)$ is isomorphic to
a moduli space of stable sheaves on the Kummer surface associated to $X$.

Sections 5--8 are appendices.
In section 5, we consider family of moduli spaces
of stable sheaves induced by a family of abelian (or K3)
surfaces,
which is necessary for deformation arguments in section 4.3.
Section 6 is devoted to getting more examples of birational maps
of moduli spaces.
In section 7, we generalize results in section 3.1
to more general Fourier-Mukai functor
(\cite{BBH:1}, \cite{Br:2}, \cite{Mu:7}, \cite{Mu:8}, \cite{Or:1}).
In particular, we can treat Fourier-Mukai functor 
on K3 surfaces. 
Under similar minimality conditions 
for $c_1(E)$ (cf. \cite{Y:5}),
our method in section 3.1 also works for these cases
(Theorem \ref{thm:fourier}). 
As an interesting corollary of Theorem \ref{thm:fourier},
we shall construct two different K3 surfaces $X,Y$
such that $\Hilb_X^2 \cong \Hilb_Y^2$ (see Example \ref{ex:birat}).
By a modification of the proof of Theorem \ref{thm:deform equiv},
we can also consider deformation types of 
moduli spaces of sheaves on K3 surfaces.
This is treated in section 8.

This paper is an extended version of my papers \cite{Y:6}.

\section{Preliminaries}

{\it Notation.}
 
Let $M$ be a complex manifold.
For a cohomology class $x \in H^*(M,{\Bbb Z})$,
$[x]_i \in H^{2i}(X,{\Bbb Z})$ denotes the $2i$-th
component of $x$.
For a projective manifold $M$, $\Amp(M) \subset H^2(M,{\Bbb Q})$
is the ample cone of $M$.

Let $p:X \to \Spec({\Bbb C})$ be an abelian surface  or a K3 surface
over ${\Bbb C}$.
We denote the projection $S \times X \to S$ by $p_S$.
In this paper, we identify a divisor class $D$ with the associated 
line bundle ${\cal O}_X(D)$ in many cases.

If $\NS(X)\cong {\Bbb Z}$, then
for a coherent sheaf $E$ on $X$,
we set
\begin{equation}
\deg(E):=\frac{(c_1(E),c_1(H))}{(c_1(H)^2)} \in {\Bbb Z},
\end{equation}
where $H$ is the ample generator of $\NS(X)$.

\subsection{Mukai lattice}

We shall recall the Mukai lattice [Mu4].
\begin{defn}
We define a symmetric bilinear form on 
$H^{ev}(X,{\Bbb Z}):=\oplus _i H^{2i}(X,{\Bbb Z})$:
\begin{align*}
\langle x, y \rangle:=&-\int_X(x^{\vee} y)\\
=&\; \int_X(x_1y_1-x_0y_2-x_2y_0)
\end{align*}
where $x=x_0+x_1+x_2, x_i \in  H^{2i}(X,{\Bbb Z})$
(resp. $y=y_0+y_1+y_2, y_i \in  H^{2i}(X,{\Bbb Z})$) and
$x^{\vee}=x_0-x_1+x_2$.
We define weight 2 Hodge structure by
\begin{equation}
 \begin{cases}
  H^{0,2}(H^{ev}(X,{\Bbb C}))=H^{0,2}(X)\\
  H^{1,1}(H^{ev}(X,{\Bbb C}))=H^{0,0}(X) \oplus H^{1,1}(X) \oplus H^{2,2}(X)\\
  H^{2,0}(H^{ev}(X,{\Bbb C}))=H^{2,0}(X).
 \end{cases}
\end{equation}
We call this lattice Mukai lattice.
\end{defn}
For a coherent sheaf $E$ on $X$,
\begin{equation}
v(E):=\ch(E)\sqrt{\td_X}=\ch(E)(1+\varepsilon \omega)
\end{equation}
is the Mukai vector of $E$, 
where $\omega$ is the fundamental class of $X$ and 
$\varepsilon=0,1$ according as $X$ is of type abelian or K3.
Then Riemann-Roch theorem is written as follows:
\begin{equation}\label{eq:RR}
\chi(E,F)=-\langle v(E),v(F) \rangle,
\end{equation}
where $E$ and $F$ are coherent sheaves on $X$.

For an abelian surface $X$,
Mukai lattice has a decomposition 
\begin{equation}\label{eq:decomp}
H^{ev}(X,{\Bbb Z}) =H^2(X,{\Bbb Z}) \oplus H^0(X,{\Bbb Z}) \oplus
H^4(X,{\Bbb Z})=U^{\oplus 4},
\end{equation}
where $U$ is the hypabolic lattice.

\subsection{Moduli space of stable sheaves}
Let $X$ be an abelian or a K3 surface, and    
$H$ an ample line bundle on $X$.
For $v \in H^{ev}(X,{\Bbb Z})$, 
let $M_H(v)$ be the moduli space of (Gieseker) stable sheaves $E$
of Mukai vector $v(E)=v$ and $\overline{M}_H(v)$ the Gieseker
compactification of $M_H(v)$ obtained by adding semi-stable sheaves.
By Mukai \cite{Mu:3}, $M_H(v)$ is smooth of dimension $
\langle v^2 \rangle+2$ and has a symplectic structure.
We choose an ample divisor $H$ on $X$
which does not lie on walls with respect to $v$ (\cite{Y:2}).
Then we have
\begin{itemize}
\item[$(\natural)$]
for every $\mu$-semi-stable sheaf $E$ of $v(E)=v$,
if $F \subset E$ satisfies
$(c_1(F),H)/\rk F=(c_1(E),H)/\rk E$, then
$c_1(F)/\rk F=c_1(E)/\rk E$,
\end{itemize}
by the definition of walls.
In particular, if $v$ is primitive, then 
$\overline{M}_H(v)=M_H(v)$, and hence $M_H(v)$ is compact.

We set
$$
v^{\perp}:=\{x \in H^{ev}(X,{\Bbb Z})| \langle v, x \rangle=0 \}.
$$
Let $\theta_v: v^{\perp} \to H^2(M_H(v),{\Bbb Z})$
be Mukai homomorphism defined by
\begin{equation}
\theta_v(x):=-\frac{1}{\rho}\left[p_{M_H(v)*}((\ch {\cal E})\sqrt{\td_X}
x^{\vee})\right]_1
\end{equation}
where ${\cal E}$ is a quasi-universal family of similitude 
$\rho$.
For a line bundle $L$ on $X$,
let $T_L:H^{ev}(X,{\Bbb Z}) \to H^{ev}(X,{\Bbb Z})$
be the homomorphism sending $x$ to $x\ch(L)$.
Then $T_L$ is an isometry of Mukai lattice.
\begin{lem}\label{lem:isomT}
If $\rk v>0$ and $H$ is general, then 
$T_L$ defines an isomorphism 
$M_H(v) \to M_H(T_L(v))$ sending $E \in M_H(v)$ to
$E \otimes L \in M_H(T_L(v))$.
Under this identification,
\begin{equation}\label{eq:isomT}
\theta_{T_L(v)}(T_L(x))=\theta_v(x), 
\end{equation}
for $x \in v^{\perp}$.
\end{lem}

\begin{proof}
Let $E$ be a stable sheaf of $v(E)=v$.
Assume that $E \otimes L$ is not stable.
Since the operation $E \mapsto E \otimes L$ preserves 
$\mu$-semi-stability,
there is a subsheaf $F \subset E$ such that
\begin{equation}\label{eq:T}
\begin{split}
 \frac{(c_1(F),H)}{\rk F}&=\frac{(c_1(E),H)}{\rk E},\\
 \frac{\chi(F \otimes L)}{\rk F}& \geq \frac{\chi(E \otimes L)}{\rk E}.
\end{split}
\end{equation}
Since $H$ is general, \eqref{eq:T} implies that
$c_1(F)/\rk F=c_1(E)/\rk E$.
A simple calculation shows that
\begin{equation}
 \begin{split}
  \frac{\chi(E \otimes L)}{\rk E}-\frac{\chi(F \otimes L)}{\rk F}
  &=\frac{\chi(E)}{\rk E}-\frac{\chi(F)}{\rk F}+
  \frac{(c_1(E),L)}{\rk E}-\frac{(c_1(F),L)}{\rk F}\\
  &=\frac{\chi(E)}{\rk E}-\frac{\chi(F)}{\rk F}.
 \end{split}
\end{equation}
This means that $E$ is not stable.
Hence $E \otimes L$ must be a stable sheaf.
\eqref{eq:isomT} follows from direct computations.
\end{proof}
 
\begin{rem}
If $\rk v=0$, then twisting by a line bundle $L$ does not
preserve the stability in general.
\end{rem}

\subsection{Beauville's bilinear form}
Let $M$ be an irreducible symplectic manifold
of dimension $2n$.
Beauville \cite{B:1} constructed a primitive symmetric bilinear
form 
\begin{equation}
B_M:H^2(M,{\Bbb Z}) \times H^2(M,{\Bbb Z}) \to {\Bbb Z}.
\end{equation}
Up to multiplication by positive constants,
$q_M(x):=B_M(x,x)$ satisfies that
\begin{equation}
q_M(x)=\frac{n}{2}\int_M \phi^{n-1} \overline{\phi}^{n-1} x^2
+(1-n)\int_M \phi^{n} \overline{\phi}^{n-1} x
\int_M \phi^{n-1}  \overline{\phi}^{n} x,
\end{equation}
where $\phi$ is a holomorphic 2 form with 
$\int_M \phi^{n}\overline{\phi}^{n}=1$.
For $\lambda, x \in H^2(M,{\Bbb C})$,
the following relation holds \cite[Thm. 5]{B:1}:
\begin{equation}\label{eq:Bform}
v(\lambda)^2q_M(x)=q_M(\lambda)\left[
(2n-1)v(\lambda)\int_M \lambda^{2n-2} x^2
-(2n-2)\left(\int_M \lambda^{2n-1} x\right)^2 \right],
\end{equation}
where $v(\lambda)=\int_M \lambda^{2n}$.

\subsection{Stable sheaves of pure dimension 1}
Let $X$ be an arbitrary surface and 
$H$ an effective divisor of $(H^2)>0$ on $X$.
Let $E$ be a purely 1-dimensional sheaf with $c_1(E)=c_1(H)$ and
$\chi(E)=\chi$.
For an ample divisor $L$,
$\chi(E \otimes L^{\otimes n})=(H,L)n+\chi$.
Hence $E$ is semi-stable with respect to $L$, if
$$
\frac{\chi(F)}{(c_1(F),L)} \leq \frac{\chi(E)}{(c_1(E),L)} 
$$
for any proper subsheaf $F \ne 0$ of $E$,
and $E$ is stable if the inequality is strict.
We assume that $\chi \ne 0$.
Under this assumption,
we shall generalize the concept of the chamber structure of polarizations.

For a subsheaf $F$ of $E$,
we set $\xi:=\chi(F)c_1(E)-\chi(E)c_1(F)$. 
Since $\chi(E) \ne 0$, we get $\xi \ne 0$ if $c_1(F) \not \in
{\Bbb Q}\; c_1(E)$.
For such a $\xi \ne 0$, we set
$W_{\xi}:=\{x \in \Amp(X)|(x,\xi)=0\}$.
Then 
\begin{equation}
(\xi^2)=\chi(F)^2(c_1(E)^2)-2\chi(E)\chi(F)(c_1(E),c_1(F))
+\chi(E)^2(c_1(F)^2).
\end{equation}
If $W_{\xi}$ is not empty,
then the Hodge index theorem implies that $(\xi^2) \leq 0$.
We claim that the choice of $c_1(F)$ is finite and only depends on the pair
$(c_1(H),\chi)$.
Then the choice of $\chi(F)$ is finite,
which shows that the number of
non-empty walls $W_{\xi}$ is finite.

Proof of the claim:
We fix an ample line bundle $L_0$ on $X$.
Let $D$ be an effective divisor on $X$ such that
$\det(F)={\cal O}_X(D)$. 
We note that $D$ satisfies the inequality
$0<(D,L_0) \leq (c_1(E),L_0)$.
For an integer $d$ with
$0<d \leq (c_1(E),L_0)$,
we set
\begin{equation}
 {\cal C}_d:=\{D|\text{$D$ is an effective divisor of $(D,L_0)=d$ } \}.
\end{equation}
It is well known that ${\cal C}_d$ is bounded.
Hence the choice of $c_1(F)$ is finite.
Thus our claim holds.

We shall call a connected component of $\Amp(X) \setminus
\cup_{\xi} W_{\xi}$ {\it a chamber}.
\begin{lem}
Assume that $(c_1(H),\chi)$ is primitive in $\NS(X) \times {\Bbb Z}$
and $\chi \ne 0$. 
Then the moduli space of stable sheaves $E$ 
of $(\rk(E),c_1(E),\chi(E))=(0,c_1(H),\chi)$ 
with respect to $L$ is compact for a general ample divisor $L$.
\end{lem}

\section{Fourier-Mukai functor}

\subsection{Notation}

Let $((X_1,H_1), (X_2,H_2),{\cal P})$ be a triple of polarized
K3 or abelian surfaces $(X_1,H_1), (X_2,H_2)$ and 
a coherent sheaf ${\cal P}$ on $X_1 \times X_2$ such that
\begin{enumerate}
\item
${\cal P}$ is flat over $X_1$ and $X_2$,
\item
\begin{equation}\label{eq:str-simple1}
\begin{split}
\Hom({\cal P}_{|\{x_1\} \times X_2},{\cal P}_{|\{x_1\} \times X_2})
&={\Bbb C}_{x_1},\\
\Ext^i({\cal P}_{|\{x_1\} \times X_2},{\cal P}_{|\{y_1\} \times X_2})
&=0,\;x_1 \ne y_1\; \text{ for $0 \leq i \leq 2$}.
\end{split}
\end{equation}
\item
\begin{equation}\label{eq:str-simple2}
\begin{split}
\Hom({\cal P}_{|X_1 \times \{x_2\} },{\cal P}_{|X_1 \times \{x_2\}})
&={\Bbb C}_{x_2},\\
\Ext^i({\cal P}_{|X_1 \times \{x_2\} },{\cal P}_{|X_1 \times \{y_2\}})
&=0,\;x_2 \ne y_2\; \text{ for $0 \leq i \leq 2$}.
\end{split}
\end{equation}
\end{enumerate}
We denote the projections $X_1 \times X_2 \to X_i$, $i=1,2$ 
 by $p_{X_i}$.
Let ${\mathbf D}(X_1)$ and ${\mathbf D}(X_2)$ be 
the bounded derived category of
coherent sheaves on $X_1$ and $X_2$ respectively.
Then the functor ${\cal F}_{\cal P}:{\mathbf D}(X_1) \to {\mathbf D}(X_2)$ 
defined by
\begin{equation}
{\cal F}_{\cal P}(x)={\mathbf R}p_{X_2*}({\cal P} \otimes p_{X_1}^*(x)), 
x \in {\mathbf D}(X_1)
\end{equation}
gives an equivalence of categories.
${\cal F}_{\cal P}$ is called Fourier-Mukai functor defined by ${\cal P}$.
We define $\widehat{{\cal F}}_{\cal P}:{\mathbf D}(X_2) \to {\mathbf D}(X_1)$ 
by
\begin{equation}
 \widehat{{\cal F}}_{\cal P}(y)=
{\mathbf R}\Hom_{p_{X_1}}({\cal P},p_{X_2}^*(y)),
 y \in {\mathbf D}(X_2).
\end{equation}
Then $\widehat{{\cal F}}_{\cal P}[2]$ is the inverse of ${\cal F}_{\cal P}$.
For a coherent sheaf $E$ on $X_1$ (resp. a
coherent sheaf $F$ on $X_2$),
we set
\begin{equation}
 \begin{split}
  {\cal F}^i_{\cal P}(E):&=H^i({\cal F}_{\cal P}(E))=
  R^ip_{X_2*}({\cal P} \otimes p_{X_1}^*(E)),\\
  \widehat{{\cal F}}^i_{\cal P}(F):&=H^i(\widehat{{\cal F}}_{\cal P}(F))
  =\Ext^i_{p_{X_1}}({\cal P},p_{X_2}^*(F)).
 \end{split}
\end{equation}
\begin{defn}
 $E$ (resp. $F$) satisfies $\WIT_i$ with respect to ${\cal F}_{\cal P}$
 (resp. $\widehat{{\cal F}}_{\cal P}$), if
 ${\cal F}_{\cal P}^j(E)=0$ (resp. $\widehat{{\cal F}}_{\cal P}^j(F)=0$)
 for $j \ne i$.
\end{defn}
Let ${\cal F}_{\cal P}^H:H^{*}(X_1,{\Bbb Q}) \to H^{*}(X_2,{\Bbb Q})$ and
$\widehat{{\cal F}}_{\cal P}^H:H^{*}(X_2,{\Bbb Q}) \to H^{*}(X_1,{\Bbb Q})$ 
be homomorphisms such that
\begin{align}
 {\cal F}_{\cal P}^H(x)&=
 p_{X_2*}\left((\ch{\cal P})p_{X_1}^*\sqrt{\td_{X_1}}p_{X_2}^*\sqrt{\td_{X_2}}
 p_{X_1}^*(x)\right), x \in H^{*}(X_1,{\Bbb Q}),\\
 \widehat{{\cal F}}_{\cal P}^H(y)&=
 p_{X_1*}\left((\ch{\cal P})^{\vee}p_{X_1}^*\sqrt{\td_{X_1}}p_{X_2}^*\sqrt{\td_{X_2}}
 p_{X_2}^*(y)\right), y \in H^{*}(X_2,{\Bbb Q}).
\end{align}
By Grothendieck-Riemann-Roch theorem,
the following diagram is commutative.
\begin{equation}
 \begin{CD}
  {\mathbf D}(X_1) @>{{\cal F}_{\cal P}}>> {\mathbf D}(X_2)\\
  @V{\sqrt{\td_{X_1}}\ch}VV @VV{\sqrt{\td_{X_2}}\ch}V\\
  H^{ev}(X_1,{\Bbb Q}) @>{{\cal F}_{\cal P}^H}>> H^{ev}(X_2,{\Bbb Q})
 \end{CD}
\end{equation}

\begin{lem}\label{lem:comm}
For $x \in H^{ev}(X_1,{\Bbb Z})$, $y \in H^{ev}(X_2,{\Bbb Z})$,
we get $\langle {\cal F}_{\cal P}^H(x),y\rangle=
\langle x,\widehat{{\cal F}}_{\cal P}^H(y)\rangle$.
\end{lem}

\begin{proof}
\begin{equation}
 \begin{split}
  \langle {\cal F}_{\cal P}^H(x),y\rangle &=
  -\int_{X_2}p_{X_2*}((\ch{\cal P})p_{X_1}^*\sqrt{\td_{X_1}}p_{X_2}^*\sqrt{\td_{X_2}}
  p_{X_1}^*(x))y^{\vee}\\
  &=-\int_{X_1 \times X_2}
  ((\ch{\cal P})p_{X_1}^*\sqrt{\td_{X_1}}p_{X_2}^*\sqrt{\td_{X_2}}
  p_{X_1}^*(x)p_{X_2}^*(y)^{\vee})\\
  &=-\int_{X_1 \times X_2}p_{X_1}^*(x)
  ((\ch{\cal P})^{\vee}p_{X_1}^*\sqrt{\td_{X_1}}p_{X_2}^*\sqrt{\td_{X_2}}
  p_{X_2}^*(y))^{\vee}\\
  &=-\int_{X_1}x\{
  p_{X_1*}((\ch{\cal P})^{\vee}p_{X_1}^*\sqrt{\td_{X_1}}p_{X_2}^*\sqrt{\td_{X_2}}
  p_{X_2}^*(y))\}^{\vee}\\
  &=\langle x,\widehat{{\cal F}}_{\cal P}^H(y)\rangle.
 \end{split}
\end{equation}
\end{proof}
\begin{lem}\label{lem:isometry of F_P}
For $x \in H^{ev}(X_1,{\Bbb Z})$,
${\cal F}_{\cal P}^H(x)$ belongs to $H^{ev}(X_2,{\Bbb Z})$.
In particular ${\cal F}_{\cal P}^H$ is an isometry of Mukai lattice.
\end{lem}
\begin{proof}
By \cite[sect. 2]{Y:4}, 
$[{\cal F}_{\cal P}^H(x)]_1$ belongs to $H^2(X_2,{\Bbb Z})$.
By Lemma \ref{lem:comm}, we get
\begin{equation}
 \begin{split}
  \langle {\cal F}_{\cal P}^H(x), 1 \rangle &=
  \langle x,\widehat{\cal F}_{\cal P}^H(1) \rangle \in {\Bbb Z},\\
  \langle {\cal F}_{\cal P}^H(x), \omega_2 \rangle &=
  \langle x,\widehat{\cal F}_{\cal P}^H(\omega_2) \rangle \in {\Bbb Z}.
 \end{split}
\end{equation}
Hence ${\cal F}_{\cal P}^H(x) \in H^{ev}(X_2,{\Bbb Z})$.
\end{proof}

We are also interested in the composition of ${\cal F}_{\cal P}$ and the
``taking-dual'' functor ${\cal D}_{X_2}:{\mathbf D}(X_2) \to
{\mathbf D}(X_2)_{op}$
sending $x \in {\mathbf D}(X_2)$ to 
${\mathbf R}{\cal H}om(x,{\cal O}_{X_2})$,
where ${\mathbf D}(X_2)_{op}$ is the opposite category of 
${\mathbf D}(X_2)$.
By Grothendieck-Serre duality, 
${\cal G}_{\cal P}:=({\cal D}_{X_2} \circ {\cal F}_{\cal P})[2]$
is defined by
\begin{equation}
 {\cal G}_{\cal P}(x):={\mathbf R}\Hom_{p_{X_2}}
 ({\cal P} \otimes p_{X_1}^*(x),{\cal O}_{X_1 \times X_2}), 
 x \in {\mathbf D}(X_1). 
\end{equation}
Let 
$\widehat{{\cal G}}_{\cal P}:{\mathbf D}(X_2)_{op} \to 
{\mathbf D}({X}_1)$ be the
inverse of ${\cal G}_{\cal P}$:
\begin{equation}
 \widehat{{\cal G}}_{\cal P}(y):={\mathbf R}\Hom_{p_{{X}_1}}
 ({\cal P} \otimes p_{X_2}^*(y),{\cal O}_{X_1 \times X_2}), 
 y \in {\mathbf D}(X_2). 
\end{equation}
For a coherent sheaf $E$ on $X_1$ (resp. a
coherent sheaf $F$ on $X_2$),
we set
\begin{equation}
 \begin{split}
  {\cal G}^i_{\cal P}(E):&=H^i({\cal G}_{\cal P}(E)),\\
  \widehat{{\cal G}}^i_{\cal P}(F):&=H^i(\widehat{{\cal G}}_{\cal P}(F)).
 \end{split}
\end{equation}
Then there are spectral sequences
\begin{equation}\label{eq:spectral3}
E_2^{p,q}=\widehat{{\cal G}}_{\cal P}^p({\cal G}_{\cal P}^{-q}(E)) 
\Rightarrow
\begin{cases}
 E, \;p+q=0\\
0, \text{ otherwise},
\end{cases}
\end{equation}
\begin{equation}\label{eq:spectral4}
E_2^{p,q}={\cal G}_{\cal P}^p(\widehat{{\cal G}}_{\cal P}^{-q}(F)) 
\Rightarrow
\begin{cases}
 F,\; p+q=0\\
0, \text{ otherwise}.
\end{cases}
\end{equation}
In particular
\begin{equation}\label{eq:vanish2}
\begin{cases}
{\cal G}_{\cal P}^p(\widehat{\cal G}_{\cal P}^0(F))=0,\;p=1,2,\\
{\cal G}_{\cal P}^p(\widehat{\cal G}_{\cal P}^2(F))=0,\;p=0,1.
\end{cases}
\end{equation}
%
Let ${\cal G}_{\cal P}^H:H^{ev}(X_1,{\Bbb Z}) \to H^{ev}(X_2,{\Bbb Z})$ be the
isomorphism of lattices defined by
\begin{equation}
\begin{split}
 {\cal G}_{\cal P}^H(x):=p_{X_2*}\left(\ch ({\cal P})^{\vee}
 p_{X_1}^* \sqrt{\td_{X_1}} p_{X_2}^* \sqrt{\td_{X_2}} p_{X_1}^*(x^{\vee})\right),
 x \in H^{ev}(X_1,{\Bbb Z}).
\end{split}
\end{equation}
Since ${\cal G}_{\cal P}^H$ is the composition of
${\cal F}_{\cal P}^H$ and the isometry
$D_{X_2}:H^{ev}(X_2,{\Bbb Z}) \to H^{ev}(X_2,{\Bbb Z})$ sending
$x$ to $x^{\vee}$,
Lemma \ref{lem:isometry of F_P} implies that
${\cal G}_{\cal P}^H$ is well-defined.

\subsection{Relation to Mukai homomorphism}
Keep notations in 2.1.
Let $v \in H^{ev}(X_1,{\Bbb Z})$ be a primitive Mukai vector.
For a fixed $i$, we set $w=(-1)^i{\cal F}_{\cal P}^H(v)$.
We define an open subscheme of $M_{H_1}(v)$ by
\begin{equation}
\begin{split}
 U:=\left\{E \in M_{H_1}(v) \left|
\begin{split}
  &\text{$\text{WIT}_i$ holds for $E$ with respect to ${\cal F}_{\cal P}$}\\
  &\text{and 
  ${\cal F}_{\cal P}^i(E)$ belongs to $M_{H_2}(w)$}
\end{split} 
 \right. 
 \right\}.
\end{split}
\end{equation} 
The following lemma which follows from 
the proof of base change theorem
was proved by Mukai \cite[Thm. 1.6]{Mu:5}. 
\begin{lem}\label{lem:flat}
Let $\{{\cal E}_s\}_{s \in S}$ be a flat family of stable sheaves
on $X_1$ such that ${\cal E}_s \in U$.
Then $\{{\cal F}_{\cal P}^i({\cal E}_s)\}_{s \in S}$ is also
a flat family of stable sheaves on $X_2$.
\end{lem}
Then ${\cal F}_{\cal P}$ induces a morphism $f:U \to M_{H_2}(w)$.
By properties of ${\cal F}_{\cal P}$, $f$ is an open immersion.
We assume that $\codim_{M_{H_1}(v)}(M_{H_1}(v) \setminus U) \geq 2$
and $U \to M_{H_2}(w)$ is birational.
We denote the image of $U$ by $V$.
We set $U=X_0$ and we denote projections
$U \times X_1 \times X_2 \to X_j$ and
$U \times X_1 \times X_2 \to X_j \times X_k$ by
$q_j$ and $q_{jk}$ respectively.
We also denote the projection
$U \times X_j \to U$ by $r_j$
and the projection $U \times X_j \to X_j$ by $s_j$.
Let ${\cal E}$ be a quasi-universal family of similitude 
$\rho$ on $U \times X_1$.
By the identification $U \to V$,
$R^iq_{02*}(q_{12}^*{\cal P} \otimes q_{01}^*{\cal E})$
becomes a quasi-universal family of similitude $\rho$ on $V \times X_2$.
The following proposition shows the importance of Mukai lattice.
Indeed, it will play an important role to compute
the period of $K_H(v)$.
\begin{prop}\label{prop:comm}
${\cal F}_{\cal P}^H$ induces an isometry $v^{\perp} \to w^{\perp}$
and the following diagram is commutative.
\begin{equation}\label{eq:diagram}
 \begin{CD}
  v^{\perp} @>{(-1)^i{\cal F}_{\cal P}^H}>> w^{\perp}\\
  @V{\theta_v}VV @VV{\theta_w}V\\
  H^2(M_{H_1}(v),{\Bbb Z}) @= H^2(M_{H_2}(w),{\Bbb Z})
 \end{CD}
\end{equation}
\end{prop}
\begin{proof}
The first assertion follows from Lemma \ref{lem:comm}.
For $y \in w^{\perp}$, we see that
\begin{equation}
\begin{split}
\rho \theta_w(y) &=(-1)^{i+1}\left[r_{2*}
(\ch((-1)^i R^iq_{02*}(q_{12}^*({\cal P})
 \otimes q_{01}^*({\cal E})))
s_2^*(\sqrt{\td_{X_2}}y^{\vee}))\right]_1\\ 
&=(-1)^{i+1}\left[r_{2*}(q_{02*}(q_{12}^*(\ch{\cal P})q_{01}^*(\ch{\cal E})
q_1^*( \td_{X_1}))
s_2^*(\sqrt{\td_{X_2}}y^{\vee}))\right]_1\\
&=(-1)^{i+1}\left[q_{0*}(q_{12}^*(\ch{\cal P})q_{01}^*(\ch{\cal E})
q_1^*( \td_{X_1})
q_2^*(\sqrt{\td_{X_2}}y^{\vee}))\right]_1\\
&=(-1)^{i+1} \left[r_{1*}(q_{01*}(q_{01}^*(\ch {\cal E})
q_1^*(\sqrt{ \td_{X_1}})
q_{12}^*((\ch {\cal P}) p_{X_1}^*(\sqrt{ \td_{X_1}})
p_{X_2}^*(\sqrt{\td_{X_2}}y^{\vee}))))\right]_1\\
&=(-1)^{i+1}\left[r_{1*}((\ch{\cal E})s_1^*(\sqrt{\td_{X_1}})
s_1^*( p_{X_1*}((\ch{\cal P})^{\vee}
p_{X_1}^*(\sqrt{ \td_{X_1}})p_{X_2}^*(\sqrt{\td_{X_2}}y))^{\vee}))\right]_1\\
&=(-1)^{i+1} \left[r_{1*}((\ch{\cal E})s_1^*(\sqrt{\td_{X_1}})
s_1^*(\widehat{{\cal F}}_{\cal P}^H (y))^{\vee})\right]_1\\
&=\rho \theta_v((-1)^i \widehat{{\cal F}}_{\cal P}^H (y)).
\end{split}
\end{equation}
Since $\widehat{{\cal F}}_{\cal P}^H \circ {\cal F}_{\cal P}^H=
1_{H^{ev}(X_1,{\Bbb Z})}$,
we get \eqref{eq:diagram}.
\end{proof}
For ${\cal G}_{\cal P}$, we set
\begin{equation}
 U':=\left\{E \in M_H(v) \left|
\begin{split}
  &\text{$\text{WIT}_i$ holds for $E$ with respect to ${\cal G}_{\cal P}$}\\
  &\text{and 
  ${\cal G}_{\cal P}^i(E)$ belongs to $M_{H_2}(w^{\vee})$}
\end{split} 
\right. 
 \right\}.
\end{equation}
Then under similar assumptions on $U'$,
we get the following.
\begin{prop}\label{prop:comm2}
${\cal G}_{\cal P}^H$ induces an isometry $v^{\perp} \to (w^{\vee})^{\perp}$
and the following diagram is commutative.
\begin{equation}\label{eq:diagram2}
 \begin{CD}
  v^{\perp} @>{(-1)^{i+1}{\cal G}_{\cal P}^H}>> (w^{\vee})^{\perp}\\
  @V{\theta_v}VV @VV{\theta_{w^{\vee}}}V\\
  H^2(M_{H_1}(v),{\Bbb Z}) @= H^2(M_{H_2}(w^{\vee}),{\Bbb Z})
 \end{CD}
\end{equation}
\end{prop}

\section{Isomorphisms induced by ${\cal F}_{\cal P}$}

\subsection{Original Fourier-Mukai functor}

We start with original Fourier-Mukai functor.
So we assume that $X$ is an abelian surface.
Let $\widehat{X}$ be the dual of $X$ and ${\cal P}$ the Poincar\'{e}
line bundle on $X \times \widehat{X}$.
We shall consider functors
${\cal F}_{\cal P}$ and ${\cal G}_{\cal P}$.
%
%
In \cite[Prop. 1.17]{Mu:5}, Mukai proved the following.
\begin{lem}\label{lem:homology}
By the canonical identification
$H^i(X,{\Bbb Z})=H_{i}(\widehat{X},{\Bbb Z})$,
\begin{equation}
 {\cal F}_{\cal P}(x_i)=(-1)^{i(i+1)/2}{\mathrm {PD}}(x_i),\;x_i \in 
 H^{i}(X,{\Bbb Z}),
\end{equation}
where ${\mathrm {PD}}(x_i)$ is the Poincar\'{e} dual of 
$x_i$.
\end{lem}
We assume that $\NS(X)={\Bbb Z} H$,
where $H$ is an ample generator.
Then the dual of $X$ also satisfies the same condition.
We set $\widehat{H}:=\det(-{\cal F}_{\cal P}(H))$.
Then $\widehat{H}$ is the ample generator of $\NS(\widehat{X})$.
By Lemma \ref{lem:homology},
\begin{equation}
 {\cal F}_{\cal P}^H(r+d c_1(H)+a \omega)
 =a-d c_1(\widehat{H})+r \widehat{\omega},
\end{equation}
where $\widehat{\omega}$ is the fundamental class of $\widehat{X}$.

\subsubsection{The case where $\langle v,1 \rangle<0$.}

In this subsection, we treat the case where $\langle v,1 \rangle<0$.
We first treat the functor ${\cal G}_{\cal P}$.
\begin{prop}\label{lem:fourier4}
Let $E$ be a $\mu$-stable sheaf of 
Mukai vector $v(E)=r+c_1(H)+a \omega:=v$.
If $a>0$, then $E$ satisfies $\WIT_2$ with respect to ${\cal G}_{\cal P}$
and ${\cal G}_{\cal P}^2(E)$ is 
a $\mu$-stable sheaf of $v({\cal G}_{\cal P}^2(E))=a+c_1(\widehat{H})+r \widehat{\omega}$.
In particular, ${\cal G}_{\cal P}$ induces an isomorphism
$M_H(v) \to M_{\widehat{H}}({\cal F}_{\cal P}^H(v)^{\vee})$, if $r, a >0$.
\end{prop}
\begin{proof}
(1) $E$ satisfies $\WIT_2$:
Clearly $E \boxtimes {\cal O}_{\widehat{X}}$ is flat over $\widehat{X}$.
Hence we can use base change theorem to show 
${\cal G}_{\cal P}^0(E)={\cal G}_{\cal P}^1(E)=0$.
We first prove the following two claims:
\begin{enumerate}
\item
$ \Hom(E \otimes{\cal P}_x,{\cal O}_X)=0$
for all $x \in \widehat{X}$. 
\item
$\Ext^1(E \otimes{\cal P}_x,{\cal O}_X)=0$
except for finitely many points $x \in \widehat{X}$.
\end{enumerate}
By the stability of $E$, we get claim (i).
Suppose that $\Ext^1(E, {\cal P}_x^{\vee})=
\Ext^1(E \otimes{\cal P}_x,{\cal O}_X) \ne  0$
for distinct points $x=x_1,x_2,\dots,x_n$.
By \cite[Lem. 2.1]{Y:5}, we get a $\mu$-stable extension sheaf $G$:
\begin{equation}
0 \to \oplus_{i=1}^n {\cal P}_{x_i}^{\vee} \to G
\to E \to 0.
\end{equation}
Since $v(G)=v(E)+n$,
we see that
$\langle v(G)^2 \rangle=\langle v(E)^2 \rangle-2na$.
Since $\dim M_H(v(G)) \geq 2$, $v(G)$ satisfies 
$\langle v(G)^2 \rangle \geq 0$.
Hence $n$ must satisfy the inequality
 $n \leq \langle v(E)^2 \rangle/2a$.
In particular claim (ii) holds.

Applying base change theorem, we see that ${{\cal G}}_{\cal P}^0(E)=0$ and
${{\cal G}}_{\cal P}^1(E)$ is of dimension $0$.
This means that ${\cal G}_{\cal P}^1(E)$ satisfies $\IT_2$.
In order to prove ${\cal G}^1_{\cal P}(E)=0$, it is sufficient to prove
$\widehat{{\cal G}}_{\cal P}^2({{\cal G}}_{\cal P}^1(E))=0$.
Since ${{\cal G}}_{\cal P}^0(E)=0$,
$\widehat{{\cal G}}_{\cal P}^0({{\cal G}}_{\cal P}^0(E))=0$.
By using the spectral sequence \eqref{eq:spectral3}, 
we conclude that $\widehat{{\cal G}}_{\cal P}^2({{\cal G}}_{\cal P}^1(E))=0$.

(2) 
${\cal G}_{\cal P}^2(E)$ is torsion free\footnote{
This claim also follows from the proof of
base change theorem.}: 
Indeed, let $T$ be the torsion submodule of ${\cal G}_{\cal P}^2(E)$.
Since ${\cal G}_{\cal P}^2(E)$ is locally free in codimension 1,
$T$ is of dimension 0.
Hence $T$ satisfies $\IT_2$ and
$\widehat{{\cal G}}_{\cal P}^2(T)$ is a locally free sheaf of
$\deg(\widehat{{\cal G}}_{\cal P}^2(T))=0$.
Since $\widehat{{\cal G}}_{\cal P}^2(T)$ is a quotient of $E$,
$\widehat{{\cal G}}_{\cal P}^2(T)$ must be $0$.
Hence $T={\cal G}_{\cal P}^2(\widehat{{\cal G}}_{\cal P}^2(T))=0$.
Thus ${\cal G}_{\cal P}^2(E)$ is torsion free.

(3)
${\cal G}_{\cal P}^2(E)$ is $\mu$-stable:
If ${\cal G}_{\cal P}^2(E)$ is not $\mu$-stable, then
there is an exact sequence
\begin{equation}
0 \to A \to {\cal G}_{\cal P}^2(E) \to B \to 0,
\end{equation} 
where $B (\ne 0)$ is a $\mu$-stable sheaf of 
$\deg(B) \leq 0$.
Then we get 
\begin{equation}
\begin{split}
\widehat{{\cal G}}_{\cal P}^0(B) &=0,\\
\widehat{{\cal G}}_{\cal P}^1(B) &=\widehat{{\cal G}}_{\cal P}^0(A),\\
\end{split}
\end{equation}
and an exact sequence
\begin{equation}
 0 \to \widehat{{\cal G}}_{\cal P}^1(A) \to 
 \widehat{{\cal G}}_{\cal P}^2(B) \to
 E \to \widehat{{\cal G}}_{\cal P}^2(A) \to 0.
\end{equation}
If $B \ne {\cal P}_x^{\vee}$ for any $x \in X$,
then by using base change theorem again,
we get $\widehat{{\cal G}}_{\cal P}^2(B)=0$.
This implies that  $B$ satisfies $\WIT_1$.
By \eqref{eq:vanish2},
${\cal G}_{\cal P}^1(\widehat{{\cal G}}_{\cal P}^1(B))=
{\cal G}_{\cal P}^1(\widehat{{\cal G}}_{\cal P}^0(A))=0$.
Hence $B=0$, which is a contradiction.
If $B ={\cal P}_x^{\vee}$ for some $x \in X$,
then $\widehat{{\cal G}}_{\cal P}^1(B)=0$ and
$\widehat{{\cal G}}_{\cal P}^2(B) \cong {\Bbb C}_x$.
Hence $\widehat{{\cal G}}_{\cal P}^0(A)=0$ and
$\widehat{{\cal G}}_{\cal P}^1(A)=
\widehat{{\cal G}}_{\cal P}^2(B)= {\Bbb C}_x$.
So we get $ {\cal G}_{\cal P}^2(\widehat{{\cal G}}_{\cal P}^1(A)) \ne 0$,
which contradicts \eqref{eq:spectral4}.
Therefore ${\cal G}_{\cal P}^2(E)$ is $\mu$-stable.

If $r>0$, then ${\cal G}_{\cal P}^H(v)$ satisfies the same condition
as $v$.
Taking into account of Lemma \ref{lem:flat},
we get an isomorphism $M_H(v) \to M_{\widehat{H}}({\cal G}_{\cal P}^H(v))$.
\end{proof}
\begin{rem}\label{rem:r=0}
If $r=0$ and $a>0$, then
${\cal G}_{\cal P}$ induces an open immersion 
$M_H(v) \to M_{\widehat{H}}({\cal G}_{\cal P}^H(v))$.
In \cite{Y:3}, we proved that
$M_{\widehat{H}}({\cal G}_{\cal P}^H(v))$ is irreducible
(which will also follow from the proof of Theorem \ref{thm:period}).
Hence $M_H(v) \to M_{\widehat{H}}({\cal G}_{\cal P}^H(v))$ is an isomorphism
in this case too.
\end{rem}
By the proof of this proposition, $E$ satisfies $\IT_2$ for ${\cal G}_{\cal P}$
if $a>\langle v^2 \rangle/2$.
By Serre duality,
 $E$ satisfies $\IT_0$ for ${\cal F}_{\cal P}$.
Since ${\cal F}^0_{\cal P}(E) \cong {\cal G}_{\cal P}^2(E)^{\vee}$,
${\cal F}^0_{\cal P}(E)$ is also $\mu$-stable.
Thus we get the following.
\begin{cor}\label{lem:fourier1}
Let $E$ be a $\mu$-stable sheaf of 
Mukai vector $v(E)=r+c_1(H)+a \omega$.
We assume that $a>\langle v^2 \rangle/2$.
Then $E$ satisfies $\IT_0$ with respect to ${\cal F}_{\cal P}$
and ${\cal F}_{\cal P}^0(E)$ is 
a $\mu$-stable vector bundle of 
$v({\cal F}_{\cal P}^0(E))=a-c_1(\widehat{H})+r \widehat{\omega}$.
\end{cor}
Since $\widehat{{\cal F}}_{\cal P}=\widehat{{\cal G}}_{\cal P}
 \circ {\cal D}_{\widehat{X}}$,
we also obtain the following.
\begin{cor}\label{lem:fourier2}
Let $E$ be a $\mu$-stable vector bundle of 
Mukai vector $v(E)=r-c_1(\widehat{H})+a \widehat{\omega}$ on $\widehat{X}$.
We assume that $a>0$.
Then $E$ satisfies $\WIT_2$ with respect to $\widehat{{\cal F}}_{\cal P}$
and $\widehat{\cal F}_{\cal P}^2(E)$ is 
$\mu$-stable.
\end{cor}
\begin{rem}
We set $v:=r-c_1(\widehat{H})+a\widehat{\omega}$.
Then,
$$
\codim_{M_{\widehat{H}}(v)}\{E \in M_{\widehat{H}}(v)|\text{$E$ is not locally free}\}=r-1.
$$
Indeed, let $F$ be a locally free sheaf of rank $r$ on $X$ 
and $\Quot^n_{F/X}$ the quot-scheme parametrizing quotients
$F \to A$ such that the Hilbert polynomial of $A$ is constant $n$.
By \cite[Thm. 0.4]{Y:1}, we see that $\dim \Quot^n_{F/X}=(r+1)n$.
For $E \in M_{\widehat{H}}(v)$, $E^{\vee \vee}$ is also $\mu$-stable.
We set ${\frak Q}^n:=\cup_{F \in M_{\widehat{H}}(v+n \omega)}\Quot^n_{F/X}$.
Then ${\frak Q}^n$ has a scheme structure and $\dim {\frak Q}^n=
\dim M_{\widehat{H}}(v+n \omega)+(r+1)n$.
Since we have a natural injective morphism
$\sqcup_{n>0}{\frak Q}^n \to M_{\widehat{H}}(v)$ whose image consists of
non-locally free sheaves, we get our claim.  
\end{rem}

\subsubsection{The case where $\langle v,1 \rangle>0$.}

As in the previous subsection, we assume that $\NS(X)={\Bbb Z}H$.

\begin{prop}\label{lem:fourier3}
Let $E$ be a $\mu$-stable sheaf of 
Mukai vector $v(E)=r+c_1(H)+a \omega$.
We assume that $a<0$.
Then $E$ satisfies $\WIT_1$ with respect to ${\cal F}_{\cal P}$
and ${\cal F}_{\cal P}^1(E)$ is 
a $\mu$-stable sheaf of $v({\cal F}_{\cal P}^1(E))=
-a+c_1(\widehat{H})-r \widehat{\omega}$.
In particular,
Fourier-Mukai functor induces an isomorphism
$M_H(v) \to M_{\widehat{H}}(-{\cal F}_{\cal P}^H(v))$.
\end{prop}
\begin{proof}
(1)
$E$ satisfies $\WIT_1$:
We first show that $H^0(X,E \otimes {\cal P}_x) =0$
except for finitely many points $x \in \widehat{X}$.
Suppose that $k_i:=h^0(X,E \otimes {\cal P}_{x_i}) \ne 0$
for distinct points $x_1,x_2,\dots,x_n$.
We shall consider the evaluation map
\begin{equation}
\phi:\oplus _{i=1}^n {\cal P}_{x_i}^{\vee} \otimes 
H^0(X,E \otimes {\cal P}_{x_i}) \to E.
\end{equation}
We assume that $\sum_i k_i >r$,
that is, $\rk(\oplus _{i=1}^n {\cal P}_{x_i}^{\vee} \otimes 
H^0(X,E \otimes {\cal P}_{x_i}))>\rk(E)$.
By the proof of [Y5, Lem. 2.1],
$\phi$ is surjective in codimension 1
and $\ker \phi$ is $\mu$-stable.
We set $b:=\dim(\coker \phi)$.
Then $v(\ker \phi)=\sum_{i=1}^n k_i-(v(E)-b \omega)$.
Since $\sum_i k_i>r$, we get
\begin{equation}
\begin{split}
\langle v(\ker \phi)^2 \rangle &=
\langle v(E)^2 \rangle +2a \sum_i k_i-2b \sum_i k_i+2br\\
& \leq \langle v(E)^2 \rangle +2a \sum_i k_i.
\end{split}
\end{equation}
Since $\langle v(\ker \phi)^2 \rangle \geq 0$, we get
$\sum_i k_i \leq \langle v(E)^2 \rangle/(-2a)$.
Therefore $H^0(X,E \otimes {\cal P}_x) =0$
except for finitely many points $x \in \widehat{X}$.
Base change theorem implies that 
${\cal F}^0_{\cal P}(E)$ is a torsion sheaf of dimension 0.
Since $E$ is a torsion free sheaf on the integral scheme
$\Supp(E)$, ${\cal F}^0_{\cal P}(E)$ is torsion free.
Hence we get ${\cal F}^0_{\cal P}(E)=0$.
By the stability of $E$ and Serre duality, 
$H^2(X,E \otimes {\cal P}_x)=0$
for all $x \in \widehat{X}$.
Hence ${\cal F}^2_{\cal P}(E)=0$.
Therefore $E$ satisfies $\WIT_1$.

(2)
${\cal F}_{\cal P}^1(E)$ is torsion free:
Let $T$ be the torsion subsheaf of ${\cal F}_{\cal P}^1(E)$.
By base change theorem, ${\cal F}_{\cal P}^1(E)$ is locally free
on the open subscheme 
$U:=\{x \in \widehat{X}|H^0(X,E \otimes {\cal P}_x) \ne 0 \}$.
Hence
the proof of (1) implies that $T$ is of dimension 0.
Since $\hat{E}$ satisfies $\WIT_1$ and $T$ satisfies
$\IT_0$, $T$ must be 0.

(3)
${\cal F}_{\cal P}^1(E)$ is $\mu$-stable:
Assume that ${\cal F}_{\cal P}^1(E)$ is not $\mu$-stable.
Let $0 \subset F_1 \subset F_2 \subset \dots \subset F_s= 
{\cal F}_{\cal P}^1(E)$
be the Harder-Narasimhan filtration of ${\cal F}_{\cal P}^1(E)$.
We shall choose the integer $k$ which satisfies 
$\deg(F_i/F_{i-1})>0$ for $i \leq k$ and
$\deg(F_i/F_{i-1}) \leq 0$ for $ i > k$.
We shall prove that
$\widehat{\cal F}^2_{\cal P}(F_k)=0$ and 
$\widehat{\cal F}^0_{\cal P}({\cal F}_{\cal P}^1(E)/F_k)=0$.
Since $\deg(F_i/F_{i-1})>0, i \leq k$,
semi-stability of $F_i/F_{i-1}$ 
implies that $\widehat{\cal F}^2_{\cal P}(F_i/F_{i-1})=0, i \leq k$.
Hence $\widehat{\cal F}^2_{\cal P}(F_k)=0$.
On the other hand, we also see that
$\widehat{\cal F}^0_{\cal P}(F_i/F_{i-1}), i>k$, is of dimension 0.
Since $F_i/F_{i-1}$ is torsion free,
$\widehat{\cal F}^0_{\cal P}(F_i/F_{i-1})=0, i>k$.
Hence we conclude that 
$\widehat{\cal F}^0_{\cal P}({\cal F}_{\cal P}^1(E)/F_k)=0$.

So $F_k$ and ${\cal F}_{\cal P}^1(E)/F_k$ satisfy $\WIT_1$ and
we get an exact sequence 
\begin{equation}
0 \to \widehat{\cal F}^1_{\cal P}(F_k) \to E \to 
\widehat{\cal F}^1_{\cal P}({\cal F}_{\cal P}^1(E)/F_k) \to 0.
\end{equation}
Since $\deg(\widehat{\cal F}^1_{\cal P}(F_k))=\deg(F_k)>0$,
$\mu$-stability of $E$ implies that
$\deg(\widehat{\cal F}^1_{\cal P}(F_k))=1$ and 
$\rk(\widehat{\cal F}^1_{\cal P}(F_k))=\rk(E)$.
Thus $\widehat{\cal F}^1_{\cal P}({\cal F}_{\cal P}^1(E)/F_k)$ 
is of dimension 0.
Then $\widehat{\cal F}^1_{\cal P}({\cal F}_{\cal P}^1(E)/F_k)$ 
satisfies $\IT_0$,
which is a contradiction.
\end{proof}

\subsection{Fourier-Mukai functor on elliptic surfaces}

Let $\pi:X \to C$ be an elliptic surface with a $0$-section $\sigma$
such that every fiber
is integral.
Let $f$ be a fiber of $\pi$.
We identify a compactification of
the relative Jacobian with $\pi:X \to C$ and 
let ${\cal P}$ be a Poincar\'{e} ``line bundle'' on $X \times_C X$,
We regard ${\cal P}$ as a coherent sheaf on $X \times X$.
Let $p_i:X \times X \to X$, $i=1,2$ be two projections.
Then ${\cal P}$ is $p_i$-flat and
parametrizes ``line bundles'' on fibers of $\pi$.  
We consider Fourier-Mukai functor 
$
 {\cal F}_{\cal P}:{\mathbf D}(X) \to {\mathbf D}(X)
$
(see \cite[5]{Br:1}) defined by
\begin{equation}
 {\cal F}_{\cal P}(x):={\mathbf R}p_{2*}({\cal P}\otimes p_1^*(x)), 
 x \in {\mathbf D}(X).  
\end{equation}
We define 
$
 \widehat{\cal F}_{\cal P}:{\mathbf D}(X) \to {\mathbf D}(X)
$
by
\begin{equation}
 \widehat{\cal F}_{\cal P}(x):={\mathbf R}\Hom_{p_{1}}({\cal P},p_2^*(x \otimes K_X)), 
 x \in {\mathbf D}(X).  
\end{equation}
Then $\widehat{\cal F}_{\cal P}[2]$ is the inverse of ${\cal F}_{\cal P}$
(\cite{Br:1}).

Chern classes:
Let $\tau$ be a section of $\pi$.
It is known that ${\cal O}_X(\sigma-\tau)$ satisfies $\WIT_1$ and 
${\cal F}_{\cal P}^1({\cal O}_X(\sigma-\tau))$ is a line bundle on $\tau$.
Replacing ${\cal P}$ by ${\cal P} \otimes p_2^*N$, $N \in \Pic(C)$,
we may assume that $\chi({\cal F}_{\cal P}^1({\cal O}_X))=1$.
Then we see that 
\begin{equation}\label{eq:tau}
 \chi({\cal F}_{\cal P}^1({\cal O}_X(\sigma-\tau)))=\chi({\cal F}_{\cal P}^1({\cal O}_X))-
 (\tau-\sigma,\sigma)=1-(\tau-\sigma,\sigma).
\end{equation}
Indeed, let $\eta$ be the generic point of $C$.
Then we have that
${\cal F}_{\cal P}({\cal O}_{\sigma}(\sigma))_{|\pi^{-1}(\eta)}=
{\cal O}_{|\pi^{-1}(\eta)}$ and
${\cal F}_{\cal P}({\cal O}_{\tau}(\sigma))_{|\pi^{-1}(\eta)}=
{\cal O}_{|\pi^{-1}(\eta)}(\tau-\sigma)$.
Hence $\chi({\cal F}_{\cal P}({\cal O}_{\sigma}(\sigma)))-
\chi({\cal F}_{\cal P}({\cal O}_{\tau}(\sigma)))=((\tau-\sigma)^2)/2$.
Since $\chi({\cal F}_{\cal P}({\cal O}_X(\sigma-\tau))=
\chi({\cal F}_{\cal P}({\cal O}_{\sigma}(\sigma)))-
\chi({\cal F}_{\cal P}({\cal O}_{\tau}(\sigma)))-
\chi({\cal F}_{\cal P}({\cal O}_X))$,
we get \eqref{eq:tau}.

We set
\begin{equation}
 d(\tau):=(\tau-\sigma,\sigma).
\end{equation}
For a coherent sheaf $E$ of 
$(\rk(E),c_1(E),-\ch_2(E))=(r,\sigma-\tau+(l+d(\tau))f,n)$,
\begin{equation}
 (\rk({\cal F}_{\cal P}(E)),c_1({\cal F}_{\cal P}(E)),\chi({\cal F}_{\cal P}(E)))
 =-(0,\tau-\sigma+r\sigma+(n-d(\tau))f,r+l).
\end{equation}

\begin{proof}
For $I_Z(\sigma-\tau+(l+d(\tau))f), Z=\{x_1,\dots,x_{n-d(\tau)}\}$,
$I_Z(\sigma-\tau+(l+d(\tau))f)$ 
satisfies $\WIT_1$ and there is an exact sequence
\begin{equation}
 0 \to \oplus_{i=1}^{n-d(\tau)} {\cal P}_{x_i} \to 
{\cal F}_{\cal P}^1(I_Z(\sigma-\tau+(l+d(\tau))f)) \to
 {\cal F}_{\cal P}^1({\cal O}_X(\sigma-\tau+(l+d(\tau))f)) \to 0.
\end{equation}
Since ${\cal F}_{\cal P}^1({\cal O}_X(\sigma-\tau+(l+d(\tau)))f))\cong 
{\cal F}_{\cal P}^1({\cal O}_X(\sigma-\tau))((l+d(\tau))f)$,
we get our claim for $r=1$ case.
For general cases, we use $E=I_Z(\sigma-\tau+(l+d(\tau))f) 
\oplus {\cal O}_X^{\oplus (r-1)}$.
\end{proof}

The following is an easy consequence of the proof of base change theorem.
\begin{lem}\label{lem:base}
Let $L$ be a coherent sheaf of pure dimension 1 on 
$X$ with $c_1(L)=\tau-\sigma+r \sigma+(n-d(\tau))f$.
Then $L$ satisfies $\WIT_1$ and $\widehat{{\cal F}}_{\cal P}^1(L)$ is torsion free, 
if the following two
conditions are satisfied:
\begin{enumerate}
\item[(1)]
 $\Hom({\cal P}_x,L)=0$, $x \in X$ except finite subset $S$ of $X$.
\item[(2)]
 $\Ext^2({\cal P}_x,L)\cong \Hom(L,{\cal P}_x)^{\vee}=0$ for all $x \in X$. 
\end{enumerate}
\end{lem}

\begin{proof}
We choose a locally free resolution
\begin{equation}
 0 \to V_2 \to V_1 \to V_0 \to {\cal P} \to 0
\end{equation}
of ${\cal P}$ such that $\Ext^j_{p_{X_1}}(V_i,p_2^*L)=0$, $j \geq 1$.
Then ${\mathbf R}\Hom_{p_{X_1}}({\cal P},p_2^*L)$ is quasi-isomorphic to
the complex
\begin{equation}
 0 \to \Hom_{p_{X_1}}(V_0,p_2^*L) \overset{\phi_1}
 \to \Hom_{p_{X_1}}(V_1,p_2^*L) \overset{\phi_2} \to
 \Hom_{p_1}(V_2,p_2^*L) \to 0.
\end{equation}
By our assumption, there is a finite subset $S$ of $X$ such that
$(\phi_1)_x$ is injective for $x \in X \setminus S$ and
$(\phi_2)_x$ is surjective for all $x \in X$.
Hence we get
\begin{enumerate}
\item[(1)]
$\Hom_{p_1}({\cal P},p_2^*L)=\ker \phi_1=0$,
\item[(2)]
$\Ext_{p_1}^2({\cal P},p_2^*L)=\coker \phi_2=0$,
\item[(3)]
$\ker \phi_2$ is a vector bundle and
\item[(4)]
there is an exact sequence
\begin{equation}
 0 \to \Hom_{p_1}(V_0,p_2^*L) \overset{\phi_1} \to
 \ker \phi_2 \to \Ext^1_{p_1}({\cal P},p_2^*L) \to 0.
\end{equation}
\end{enumerate}
Hence $\Ext^1_{p_1}({\cal P},p_2^*L)$ is a vector bundle on $X \setminus S$.
By (4), $\Ext^1_{p_1}({\cal P},p_2^*L)$ is torsion free.
\end{proof}

\begin{lem}\label{lem:H0}
For a purely 1-dimension sheaf $L$ with 
$c_1(L)=\tau-\sigma+r\sigma+(n-d(\tau))f$,
$\Hom_{p_1}({\cal P},p_2^* L)=0$.
\end{lem}

\begin{proof}
By the proof of Lemma \ref{lem:base},
it is sufficient to prove that
$\Hom({\cal P}_x,L)=0$ for some point $x \in X$.
We choose a point $x \in X$ which is not contained in
$\Supp(L)$.
Since $L$ is of pure dimension 1, we get 
$\Hom({\cal P}_x,L)=0$.
\end{proof}

\begin{prop}\label{prop:base-conditions}
Let $L$ be a coherent sheaf of pure dimension 1 on 
$X$ with $c_1(L)=\tau-\sigma+r \sigma+(n-d(\tau))f$ and $\chi(L)>0$.
If $L$ is semi-stable with respect to $\sigma+kf$, $k \gg 0$, then
$L$ satisfies the above two conditions.
\end{prop}

\begin{proof}
By taking account of a Jordan-H\"{o}lder filtration, we may assume that $L$
is stable.
We note that 
\begin{equation}
 \mu(L):=\frac{\chi(L)}{(c_1(L),\sigma+kf)}=
 \frac{\chi(L)}{(k+(\sigma^2))r+n}.
\end{equation} 
Since $\mu({\cal P}_x)=0$, $x \in X$,
obviously (2) holds.
So we shall prove (1).
Let $D=D_{vir}+D_{hol}$ be the decomposition of the scheme-theoretic
support of $L$, where $D_{vir}$ consists of all fiber components and
$D_{hol}$ consists of the other components.
Then we have an exact sequence
\begin{equation}
 0 \to F \to L \to (L_{|D_{hol}})/T \to 0,
\end{equation}
where $T$ is the 0-dimensional submodule of $L_{|D_{hol}}$.   
Then $F$ is a purely 1-dimensional subsheaf of $L$ with $c_1(F)=lf$.
By the stability of $L$, we get 
\begin{equation}
 \mu(F)=\frac{\chi(F)}{l}\leq \frac{\chi(L)}{(k+(\sigma^2))r+n}.
\end{equation}
Since $k$ is sufficiently large (the condition
$k > \max\{((c_1(L),\sigma+k_0 f)\chi(L)-n)/r-(\sigma^2), k_0\}$ 
is sufficient, where $k_0$ satisfies that $\sigma+k_0 f$ is ample),
we get $\chi(F) \leq 0$.
Since $\Hom({\cal P}_x,L_{|D_{hol}}/T)=0$ for all $x \in X$,
we shall prove that
$\Hom({\cal P}_x,F)=0$ except finite numbers of points.
\newline
Proof of the claim:
Let 
\begin{equation}
 0 \subset F_1 \subset F_2 \subset \cdots \subset F_s=F
\end{equation}
be the Harder-Narasimhan filtration of $F$ with respect to $\sigma+kf$.
Then 
\begin{equation}
 \mu(F_1)>\mu(F_2/F_1)>\dots>\mu(F_s/F_{s-1}).
\end{equation}
Since $F_1$ is a subsheaf of $L$, 
we also have the inequality $\chi(F_1) \leq 0$. 
If $\Hom({\cal P}_x,F_1) \ne 0$, then 
$\mu(F_1)=0$ and $F_1$ is $S$-equivalent to
${\cal P}_x \oplus E$ for some $E$.
Hence the choice of $x$ is finite.
Clearly $\Hom({\cal P}_x,F_i/F_{i-1})=0$ for $i \geq 2$.
Hence the claim holds. 
\end{proof}
 
\begin{rem}
Stability is not preserved by the operation 
$L \mapsto L \otimes N$, $N \in \Pic(X)$.
Hence the condition $\chi(L)>0$ is important.
\end{rem}

\begin{lem}\label{lem:wit F}
Let $E$ be a torsion free sheaf of $\rk(E)=r>0$ and
$c_1(E)=\sigma-\tau+(l+d(\tau))f$ on $X$.
Assume that $E$ is semi-stable with respect to $\sigma+k f$, $k \gg 0$.
Then $E$ satisfies $\WIT_1$ and ${\cal F}_{\cal P}^1(E)$ is of pure dimension 1.
\end{lem}

\begin{proof}
We shall first prove that ${\cal P} \otimes p_1^*(E)$ is $p_2$-flat.
Let 
\begin{equation}
 0 \to W_1 \to W_0 \to {\cal P} \to 0
\end{equation}
be a locally free resolution of ${\cal P}$ on $X \times X$.
It is sufficient to prove that
\begin{equation}
 \psi_x:(W_1)_{|x \times X} \otimes E \to (W_0)_{|x \times X}\otimes E
\end{equation}
is injective for all $x \in X$.  
We note that $\rk W_1=\rk W_0$ and ${\cal P}_x \otimes E$ is a 
torsion sheaf on $X$.
Since $E$ is torsion free, $\psi_x$ is injective for all $x \in X$.
Thus ${\cal P} \otimes p_1^*(E)$ is a $p_2$-flat sheaf.

Hence we can use base change theorem.
Since $p_2:X \times_C X \to X$ is relative dimension 1,
$R^2p_{2*}({\cal P} \otimes p_1^*(E))=0$.
Since $E_{|\pi^{-1}(y)}$ is semi-stable for general $y \in C$,
$H^0(X,{\cal P}_x \otimes E)=0$ for a general point $x$ of $X$.
Thus $p_{2*}({\cal P} \otimes p_1^*(E))$ is a torsion sheaf.
By the proof of base change theorem,
locally there is a complex of locally free sheaves 
$V_1 \overset{\phi}\to V_0$ which is quasi-isomorphic to
${\mathbf R}p_{2*}({\cal P}\otimes p_1^*(E))$.
Hence $p_{2*}({\cal P} \otimes p_1^*(E))=\ker \phi=0$,
which means that $E$ satisfies $\WIT_1$.
Also we get ${\mathrm {proj-}}\dim(\coker \phi)=1$.
Hence $R^1p_{2*}({\cal P} \otimes p_1^*(E))$ is of pure dimension 1.
\end{proof}

\begin{cor}\label{cor:cal F}
Let $E$ be a torsion free sheaf on $X$ and assume that
$E_{|\pi^{-1}(y)}$ is a semi-stable vector bundle of degree 0
for a general $y \in C$.
Then $E$ satisfies $\WIT_1$ and 
${\cal F}_{\cal P}^1(E)$ is of pure dimension 1. 
\end{cor}

\begin{proof}
We consider the Harder-Narasimhan filtration of $E$
with respect to $\sigma+kf$, $k\gg 0$.
Applying ${\cal F}_{\cal P}$ to this filtration, we get our corollary
by Lemma \ref{lem:wit F}. 
\end{proof}

\begin{lem}\label{lem:cal F}
Keep the notation as above
and assume that $E$ is semi-stable with respect to 
$\sigma+k f$, $k \gg 0$.
If $r+l>0$, then ${\cal F}_{\cal P}^1(E)$ is semi-stable.
\end{lem}

\begin{proof}
Assume that ${\cal F}_{\cal P}^1(E)$ is not semi-stable.
Then, there is a stable subsheaf $F$ of ${\cal F}_{\cal P}^1(E)$ such that 
$ \mu(F)>\mu({\cal F}_{\cal P}^1(E))=(r+l)/((k+(\sigma^2))r+n)>0$
and $G:={\cal F}_{\cal P}^1(E)/F$ is of pure dimension 1.
Applying $\widehat{\cal F}_{\cal P}$ to the exact sequence
\begin{equation}
 0 \to F \to  {\cal F}_{\cal P}^1(E) \to G \to 0,
\end{equation}
we get an exact sequence
\begin{equation}
 \begin{CD}
  0 @>>>\widehat{\cal F}_{\cal P}^0(F) @>>> 0 @>>>\widehat{\cal F}_{\cal P}^0(G)\\
  @>>> \widehat{\cal F}_{\cal P}^1(F) @>>> E @>>>\widehat{\cal F}_{\cal P}^1(G)\\ 
  @>>> \widehat{\cal F}_{\cal P}^2(F) @>>> 0 @>>>\widehat{\cal F}_{\cal P}^2(G) @>>> 0. 
 \end{CD}
\end{equation}
By Lemma \ref{lem:H0}, $\widehat{\cal F}_{\cal P}^0(G)=0$.
Since $\mu(F)>0$, we also get $\widehat{\cal F}_{\cal P}^2(F)=0$.
Thus $F$ and $G$ satisfies $\WIT_1$ and $\widehat{\cal F}_{\cal P}^1(F)$ is a subsheaf of
$E$.
We set
$c_1(F)=\tau'-\sigma+r'\sigma+(n'-d(\tau'))f$ 
and $\chi(F)=r'+l'$. 
Then 
$(\rk(\widehat{\cal F}_{\cal P}^1(F)),c_1(\widehat{\cal F}_{\cal P}^1(F)),-\ch_2(\widehat{\cal F}_{\cal P}^1(F)))
=(r',\sigma-\tau'+(l'+d(\tau'))f,n')$.
Since $E$ is semi-stable with respect to $\sigma+kf$, $k \gg 0$,
\begin{itemize}
\item[($\star$)]
\begin{enumerate}
 \item
$l'/r'<l/r$,
or
\item
$l'/r'=l/r$ and 
  $-n'/r'\leq -n/r$.
\end{enumerate}
\end{itemize}
On the other hand,
\begin{equation}
 \begin{split}
  \mu({\cal F}_{\cal P}^1(E))-\mu(F)=&
  \frac{r+l}{(k+(\sigma^2))r+n}-
  \frac{r'+l'}{(k+(\sigma^2))r'+n'}\\
  =&\frac{(lr'-l'r)(k+(\sigma^2))+
  (r+l)n'-(r'+l')n}
  {((k+(\sigma^2))r+n)((k+(\sigma^2))r'+n')}.
 \end{split}
\end{equation}
We note that the choice of $r'$ is finite.
In Lemma \ref{lem:bdd},
we shall show that the choice of $n'$ is also finite.
Since $n \geq 0$ (see Lemma \ref{lem:bogomolov}),
there is an integer $N(r,\tau,l,n)$ such that for $k>N(r,\tau,l,n)$,
\begin{enumerate}
\item 
$lr'-l'r>0$ implies 
$(lr'-l'r)(k+(\sigma^2))+(r+l)n'-(r'+l')n>0$
 and 
\item
$lr'-l'r<0$ implies 
$ (lr'-l'r)(k+(\sigma^2))+(r+l)n'-(r'+l')n<0$.
\end{enumerate}
Then $(\star)$ implies that $\mu({\cal F}_{\cal P}^1(E))-\mu(F) \geq 0$,
which is a contradiction.
Therefore ${\cal F}_{\cal P}^1(E)$ is a semi-stable sheaf.
\end{proof}

\begin{lem}\label{lem:bdd}
Keep the notations as above.
Then the choice of $n'$ is finite and 
the number of such $n'$ is bounded in terms of
$(r,\tau,l,n)$.
\end{lem}
\begin{proof}
We fix an ample divisor $\sigma+k_0f$.
Since $F$ is a subsheaf of ${\cal F}_{\cal P}^1(E)$,
\begin{equation}
 0 \leq (c_1(F),\sigma+k_0f) \leq (c_1({\cal F}_{\cal P}^1(E)),\sigma+k_0f).
\end{equation}
Since $(c_1(F),\sigma+k_0f)=(k_0+(\sigma^2))r'+n'$ and 
$r' \leq r$, we get our claim.
\end{proof}

\begin{lem}\label{lem:bogomolov}
Keep the notations as above. Then $n \geq 0$.
\end{lem}
\begin{proof}
By Bogomolov's inequality,
\begin{equation}
  0 \leq c_2(E)-\frac{(r-1)}{2r}(c_1(E)^2)
  =n+\frac{((\sigma-\tau)^2)}{2r}.
\end{equation}
Since $(\sigma^2) \leq 0$, $d(\tau)=(\tau,\sigma)-(\sigma^2)\geq 0$.
Therefore $n \geq -((\sigma-\tau)^2)/2r=d(\tau)/r\geq 0$.
\end{proof}

\begin{lem}
Let $L$ be a purely 1-dimensional sheaf of 
$c_1(L)=\tau-\sigma+r\sigma+(n-d(\tau))f$ and $\chi(L)>0$.
Assume that $L$ is semi-stable with respect to $\sigma+kf$, $k\gg0$.
Then $\widehat{\cal F}_{\cal P}^1(L)$ is a semi-stable sheaf
with respect to $\sigma+kf$, $k\gg0$.
\end{lem}

\begin{proof}
We note that Lemma \ref{lem:base} and Proposition
\ref{prop:base-conditions} imply that 
$L$ satisfies $\WIT_1$ and $\widehat{\cal F}_{\cal P}^1(L)$ is torsion free.
Assume that $\widehat{\cal F}_{\cal P}^1(L)$ 
is not semi-stable with respect to $\sigma+kf$, $k\gg0$.
Then there is a destabilizing subsheaf $F$ of 
$\widehat{\cal F}_{\cal P}^1(L)$ such that
$G:=\widehat{\cal F}_{\cal P}^1(L)/F$ is torsion free.
It is easy to see that $\widehat{\cal F}_{\cal P}^1(L)_{|\pi^{-1}(y)}$ is semi-stable for
general $y \in C$.
Since $k$ is sufficiently large,
$F_{|\pi^{-1}(y)}$ and $G_{|\pi^{-1}(y)}$ are semi-stable vector bundles
of degree 0 for general $y \in C$.
Then Corollary \ref{cor:cal F} implies that
 $F$ and $G$ satisfies $\WIT_1$ and we get an exact sequence
\begin{equation}
 0 \to {\cal F}_{\cal P}^1(F) \to L \to {\cal F}_{\cal P}^1(G) \to 0.
\end{equation}
In the same way as in Lemma \ref{lem:cal F},
we get a contradiction.
Thus $\widehat{\cal F}_{\cal P}^1(L)$ is semi-stable.
\end{proof}

Therefore, we get the following theorem. 
\begin{thm}\label{thm:FM}
Let ${\cal M}(r,c_1,\chi)^{ss}$ be the moduli stack of semi-stable sheaves
$E$ of $(\rk(E),c_1(E),\chi(E))=(r,c_1,\chi)$
with respect to $\sigma+kf$, $k \gg 0$.
Then ${\cal F}_{\cal P}$ gives an isomorphism of moduli stack
\begin{equation}
{\cal M}(r,\sigma-\tau+(l+d(\tau))f,r \chi({\cal O}_X)-n)^{ss} \to 
{\cal M}(0,\tau-\sigma+r\sigma+(n-d(\tau))f,r+l)^{ss}
\end{equation}
if $r+l>0$.
\end{thm}

\begin{rem}\label{rem:general}
If $\gcd (r,l,n)=1$,
then ${\cal M}(r,\sigma-\tau+(l+d(\tau))f,r\chi({\cal O}_X)-n)^{ss}$
consists of stable sheaves. 
\end{rem}

\begin{rem}
The map
${\cal D}':L \mapsto {\cal E}xt_{{\cal O}_X}^1(L,K_X)$ gives an isomorphism
${\cal M}(0,\tau-\sigma+r\sigma+(n-d(\tau))f,r+l)^{ss} \to
{\cal M}(0,\tau-\sigma+r\sigma+(n-d(\tau))f,-(r+l))^{ss}$.
If $r+l<0$, then
${\cal G}_{\cal P}$ induces an isomorphism
\begin{equation}
{\cal M}(r,\sigma-\tau+(-(l+2r)+d(\tau))f,r \chi({\cal O}_X)-n)^{ss} \to 
{\cal M}(0,\tau-\sigma+r\sigma+(n-d(\tau))f,r+l)^{ss}
\end{equation}
\end{rem}

\begin{rem}
Independently,
Jardim and Maciocia \cite{J-M:1} treated
the case where $r+l=0$.
Hern\'{a}ndez Ruip\'{e}rez and Mu\~{n}oz Porras
\cite{H-M:1} studied stable sheaves on elliptic fibrations and
got similar results.
\end{rem}

\section{Moduli spaces of stable sheaves on abelian surfaces}

\subsection{Bogomolov decomposition}

In this section, we assume that $X$ is an abelian surface.
Let $\widehat{X}$ be the dual abelian surface of $X$ and
${\cal P}$ the Poincar\'{e} line bundle on
$\widehat{X} \times X$.
Let $v$ be a Mukai vector with $c_1(v) \in \NS(X)$ and 
$v>0$.
Fix an element $E_0 \in \overline{M}_H(v)$ and
let $\alpha:\overline{M}_H(v) \to X$ be the morphism
such that
\begin{equation}
 \alpha(E):=\det p_{\widehat{X}!}
 ((E-E_0)\otimes({\cal P}-{\cal O}_{\widehat{X} \times X}))
 \in \Pic^0(\widehat{X})=X,
\end{equation}
and $\det:\overline{M}_H(v) \to \widehat{X}$
the morphism sending $E$ to $\det E \otimes \det E_0^{\vee}
\in \widehat{X}$.
We set ${\frak a}_v:=\alpha \times \det$.
Under the assumption $\ell(v)=1$,
 the following was proved in \cite[Thm. 3.1, 3.6]{Y:3}.

\begin{thm}\label{thm:H2}
Let $v \in H^{ev}(X,{\Bbb Z})$
be a primitive Mukai vector such that 
$v>0$ and $c_1(v) \in \NS(X)$.
We assume that $\dim M_H(v)=\langle v^2 \rangle+2 \geq 6$.
Then for an ample line bundle $H$ such that
$\overline{M}_H(v)=M_H(v)$, the following holds.
\begin{enumerate}
\item[(1)]
$\theta_v$ is injective.
\item[(2)]
${\frak a}_v$ is the albanese map.
\item[(3)]
\begin{equation}\label{eq:H2}
\begin{split}
H^2(M_H(v),{\Bbb Z}) &=
 \theta_v(v^{\perp})
\oplus {\frak a}_v^* H^2(X \times \widehat{X},{\Bbb Z}).
\end{split}
\end{equation}
\end{enumerate}
\end{thm}
The proof will be done in section 4.3.
\begin{defn}
We set $K_H(v):={\frak a}_v^{-1}((0,0))$.
\end{defn}
We shall construct an \'{e}tale covering 
such that ${\frak a}_v$ becomes trivial. 
By Theorem \ref{thm:period},
it will become a Bogomolov decomposition of $M_H(v)$.

Let ${\cal F}_{\cal P}:{\mathbf D}(X) \to {\mathbf D}(\widehat{X})$ 
be the Fourier-Mukai
functor defined by ${\cal P}$.
Then 
\begin{equation}
 \alpha(E)=\det {\cal F}_{\cal P}(E) 
 \otimes (\det{\cal F}_{\cal P}(E_0))^{\vee}.
\end{equation}
For a line bundle $L$ on $X$,
we set $\tilde{L}:=\det({\cal F}_{\cal P}(L))$.
Then the following relations hold.

\begin{lem}\label{lem:phi}
 \begin{equation}
  \begin{split}
   & \phi_{\tilde{L}} \circ \phi_L=-\chi(L)1_X,\\
   &  \phi_{L} \circ \phi_{\tilde{L}}=-\chi(L)1_{\widehat{X}}. 
  \end{split}
 \end{equation}
\end{lem}
 
\begin{proof}
For $(x,\hat{x}) \in X \times \widehat{X}$,
the following diagram is commutative (\cite[(3.1)]{Mu:2}).
\begin{equation}\label{eq:tensor}
 \begin{CD}
  {\mathbf D}(X) @>{T_{-x}^*}>>{\mathbf D}(X) 
  @>{\otimes {\cal P}_{\hat{x}}}>> {\mathbf D}(X) \\
  @V{{\cal F}_{\cal P}}VV @V{{\cal F}_{\cal P}}VV
  @V{{\cal F}_{\cal P}}VV \\
  {\mathbf D}(\widehat{X}) @>{\otimes {\cal P}_{{x}}}>>
  {\mathbf D}(\widehat{X}) @>{T_{\hat{x}}^*}>>{\mathbf D}(\widehat{X})
 \end{CD}
\end{equation}
Hence we see that
${\cal F}_{\cal P}(L) \otimes {\cal P}_x={\cal F}_{\cal P}(T^*_{-x}L)
={\cal F}_{\cal P}(L \otimes {\cal P}_{\phi_L(-x)})=
T^*_{\phi_L(-x)}({\cal F}_{\cal P}(L))$.
Since $\rk({\cal F}_{\cal P}(E))=\chi(E)$, by taking determinant of
${\cal F}_{\cal P}(L) \otimes {\cal P}_x$,
we get that 
$\tilde{L} \otimes {\cal P}_{\chi(L)x}=
T_{\phi_L(-x)}^*(\tilde{L})=
\tilde{L} \otimes {\cal P}_{\phi_{\tilde{L}}\circ \phi_L(-x)}$.
Therefore the first equality holds.
Applying \eqref{eq:tensor} again,
we get that
$\widehat{{\cal F}}_{\cal P}(\tilde{L}) \otimes
{\cal P}_{\hat{x}}=
T_{-\phi_{\tilde{L}}(\hat{x})}^*\widehat{{\cal F}}_{\cal P}(\tilde{L})$.
Hence we obtain that
$\chi(\tilde{L})=-\phi_{\det {\widehat{{\cal F}}_{\cal P}(\tilde{L})}}
\circ \phi_{\tilde{L}}$.
Since $c_1({\widehat{{\cal F}}_{\cal P}(\tilde{L})})=c_1(L)$ and
$\chi(\tilde{L})=(c_1(\tilde{L})^2)/2=
(c_1(L)^2)/2=\chi(L)$,
we get the second equality.
\end{proof}
We define a morphism
$\Phi:K_H(v) \times X \times \widehat{X} \to M_H(v)$
by
$\Phi(E,x,y):=T_x^*(E) \otimes {\cal P}_y$.

\begin{lem}\label{lem:alb}
Let $v=r+c_1+a \omega$, $c_1 \in \NS(X)$ be a Mukai vector and
$L$ a line bundle on $X$ such that $c_1(L)=c_1$. Then,
\begin{equation}
\begin{split}
\alpha(T_x^*(E) \otimes {\cal P}_y) &=-ax+\phi_{\tilde{L}}(y)\\
\det(T_x^*(E) \otimes {\cal P}_y) &=\phi_L(x)+ry,
\end{split}
\end{equation}
for a coherent sheaf $E$ of $v(E)=v$.
\end{lem}
\begin{proof}
We shall only prove the first equality.
Applying \eqref{eq:tensor},
we see that 
${\cal F}_{\cal P}(T_x^*(E) \otimes {\cal P}_y)=
T_y^*({\cal F}_{\cal P}(T_x^*(E)))
=T_y^*({\cal F}_{\cal P}(E) \otimes {\cal P}_{-x})=
T_y^*({\cal F}_{\cal P}(E)) \otimes{\cal P}_{-x}$.
Hence we obtain
\begin{equation}
 \det({\cal F}_{\cal P}(T_x^*(E) \otimes {\cal P}_y))=
 T_y^*(\det {\cal F}_{\cal P}(E))\otimes{\cal P}_{-\chi(E)x}.
\end{equation}
Since $c_1( {\cal F}_{\cal P}(E))=c_1( {\cal F}_{\cal P}(L))$,
\begin{equation}
 \alpha(T_x^*(E) \otimes {\cal P}_y)=
 \phi_{\det {\cal F}_{\cal P}(E)}(y)-\chi(E)x=
 \phi_{\tilde{L}}(y)-\chi(E)x.
\end{equation}
By \eqref{eq:RR}, $\chi(E)=-\langle v({\cal O}_X),v(E) \rangle=a$,
and hence we get the first equality.
\end{proof}

Let $v=r+c_1+a \omega$, $c_1 \in \NS(X)$ be a Mukai vector and
$\tau:X \times \widehat{X} \to X \times \widehat{X}$ a homomorphism sending
$(x,y)$ to
$(rx-\phi_{\tilde{L}}(y),-\phi_L(x)-ay)$.
By Lemma \ref{lem:phi} and \ref{lem:alb},
\begin{equation}\label{eq:q-section}
 {\frak a}_v \circ \Phi \circ (1_{K_H(v)} \times \tau)
 (E,x,y)=(nx,ny),
\end{equation}
where $n=\langle v^2 \rangle/2$.
Let $\nu:X \times \widehat{X} \to X \times \widehat{X}$
be the $n$ times map and
we shall consider the fiber product
\begin{equation}
 \begin{CD}
  M_H(v) \times_{X \times \widehat{X}}X \times \widehat{X} @>>> M_H(v)\\
  @VVV @VV{{\frak a}_v}V \\
  X \times \widehat{X} @>{\nu}>> X \times \widehat{X}
 \end{CD}
\end{equation}
Then $\Phi \circ (1_{K_H(v)} \times \tau)$ and the projection
$K_H(v) \times X \times \widehat{X} \to  X \times \widehat{X}$
defines a morphism 
\begin{equation}
 K_H(v) \times X \times \widehat{X} \to
 M_H(v) \times_{X \times \widehat{X}}X \times \widehat{X}.
\end{equation}
We can easily show that this morphism is injective,
and hence it is an isomorphism.
 
\begin{rem}
If $(c_1^2)/2$ and $r$ are relatively prime,
then \cite[Prop. 4.1]{Y:3} implies that 
$M_H(v) \cong \widehat{X} \times \det^{-1}(0)$.
We shall consider the pull-back of 
${\frak a}_v: M_H(v) \to X \times \widehat{X}$ by the morphism
sending $(x,y)$ to $(nx,y)$.
Then we get $M_H(v) \times_{ X \times \widehat{X}} X \times \widehat{X}
\cong K_H(v) \times  X \times \widehat{X}$.
\end{rem}

\begin{defn}
For simplicity,
we also denote the homomorphism
$v^{\perp} \to H^2(M_H(v),{\Bbb Z}) \to H^2(K_H(v),{\Bbb Z})$
by $\theta_v$:
\begin{equation}
 \theta_v(x)=-\frac{1}{\rho}\left[p_{K_H(v)*}
 \left(\ch({\cal E}_{|K_H(v) \times X})
 x^{\vee}\right)\right]_1.
\end{equation}
\end{defn}

\begin{prop}\label{prop:dim-alb} 
Let $v$ be a Mukai vector of $c_1(v) \in \NS(X)$ and 
$v>0$.
Then 
$\dim \im {\frak a}_v \geq 2$ and
${\frak a}_v$ is surjectice if $\langle v^2 \rangle>0$.
\end{prop}

\begin{proof}
The first claim follows from Lemma \ref{lem:alb} and the second claim
follows from \eqref{eq:q-section}.
\end{proof}

\subsection{Fourier-Mukai functor on abelian surfaces}
  
In this subsection, we consider the relation 
of Fourier-Mukai functor on abelian surfaces
to the map ${\frak a}_v:M_H(v) \to X \times \widehat{X}$
(Proposition \ref{Prop:FM-alb}).
For this purpose, we recall some results of Mukai
\cite{Mu:7}.

Let $E$ be a semi-homogeneous sheaf on $X$.
We set 
\begin{equation}
K(E):=\{x \in X|T_x^*(E) \cong E \}.
\end{equation}
Let $L$ be a line bundle on $X$.
Since ${\cal F}_{{\cal P}_X}(T_x^*(E \otimes L))=
{\cal F}_{{\cal P}_X}(E \otimes L) \otimes 
({\cal P}_X)_{|\{x \} \times \widehat{X}}^{\vee}$,
$K(E \otimes L)$ is a finite set for a sufficiently ample 
line bundle $L$,
where ${\cal P}_X$ is the Poincar\'{e} line bundle on $X \times \widehat{X}$.
Replacing $E$ by $E \otimes L$, we assume that $K(E)$ is a finite set. 
Assume that $E$ is simple.
Let $Y$ be the moduli space of simple semi-homogeneous sheaves $F$
such that $\ch(F)=\ch(E)$.
Then we have a morphism $\phi_E:X \to Y$
sending $x \in X$ to $T_x^*(E) \in Y$ and we get an isomorphism
$X/K(E) \cong Y$.
By this identification, the translation $T_{\phi_E(a)}:
Y \to Y, a \in X$ is given by the morphism
$F  \mapsto T_a^*(F)$, $F \in Y$.

We assume that there is a universal family
${\cal E}$ on $Y \times X$.
Let $H$ be an ample line bundle on $Y$.
\begin{lem}
${\cal E}_{|Y \times \{x \}}$ is a stable sheaf on $Y$ 
with respect to $H$.
\end{lem}
\begin{proof}
Let $m:X \times X \to X$ be the multiplication map.
By the universal property of $Y$,
there is a line bundle $L$ on $X$ such that
$m^*(E) \cong (\phi_E \times id_X)^*{\cal E} \otimes p_1^*L$,
where $p_1:X \times X \to X$ is the first projection.
Hence $T_x^*(E) \cong \phi_E^*({\cal E}_{|Y \times \{x \}}) \otimes L$.
Since $T_x^*(E)$ is simple,
\cite[Prop. 6.16]{Mu:1} implies that
$T_x^*(E)$ is stable for all ample line bundles on $X$
(if $\dim E=1$, then $E$ is a stable vector bundle on an abelian subvariety
of $X$).
Since $\phi_E$ is finite, $\phi_E^*H$ is also ample. 
Hence $\phi_E^*({\cal E}_{|Y \times \{x \}})$ is also stable
with respect to $\phi_E^*H$, which implies that
${\cal E}_{|Y \times \{x \}}$ is stable with respect to $H$.
\end{proof}
Then it is not difficult to see that ${\cal E}$ induces 
an isomorphism $ X \to M_H(w)$, $w:=v({\cal E}_{|Y \times \{x \}})$
 and an equivalence of derived categories (cf. \cite[Thm. 2.2]{Mu:8}). 
 
\begin{cor}[Mukai \cite{Mu:7}]\label{cor:FM-const}
 ${\cal E}$ defines Fourier-Mukai transform
 ${\cal F}_{\cal E}:{\mathbf D}(X) \to {\mathbf D}(Y)$.
 \end{cor}

Let $p_X:Y \times X \to X$ (resp. $p_Y:Y \times X \to Y$)
be the projection.
We denote the Poincar\'{e} line bundle on $Y \times \widehat{Y}$
by ${\cal P}_Y$. 
 Let $\lambda:X \times \widehat{X} \to Y$
 be the morphism sending $(x,\hat{x}) \in X \times \widehat{X}$
 to $T_x^*E \otimes ({\cal P}_X)_{\hat{x}} \in Y$.
 Then the translation $T_{\lambda(a,b)}:Y \to Y$, 
$(a,b) \in X \times \widehat{X}$
 is given by the morphism sending
 $F \in Y$ to $T_a^*F \otimes ({\cal P}_X)_b$.
 There is a homomorphism $\mu:X \times \widehat{X} \to \widehat{Y}$
 such that
$$
  (id_Y \times T_x)^*{\cal E} \otimes ({\cal P}_X)_{\hat{x}}
  \cong (T_{\lambda(x,\hat{x})} \times id_X)^*{\cal E} \otimes
  ({\cal P}_Y)_{\mu(x,\hat{x})}, (x,\hat{x}) \in X \times \widehat{X}.
$$
 So we get a homomorphism
 $\lambda \times \mu:X \times \widehat{X} \to Y \times \widehat{Y}$.
 In the same way, we have a homomorphism $\lambda' \times \mu':
 Y \times \widehat{Y} \to X \times \widehat{X}$
 such that
$$
 (T_y \times id_X)^*{\cal E} \otimes ({\cal P}_Y)_{\hat{y}}
 \cong (id_Y \times T_{\lambda'(y,\hat{y})})^*{\cal E} \otimes
 ({\cal P}_X)_{\mu'(y,\hat{y})}, (y,\hat{y}) \in Y \times \widehat{Y}.
$$
  Hence we see that $(\lambda' \times \mu') \circ 
  (\lambda \times \mu)=id_{X \times \widehat{X}}$
  and $(\lambda \times \mu) \circ (\lambda' \times \mu')
  =id_{Y \times \widehat{Y}}$.
  Indeed, by restricting 
  $(id_Y \times T_x)^*{\cal E} \otimes ({\cal P}_X)_{\hat{x}}$ to 
  $Y \times \{a \}$, we get that ${\cal E}_{|Y \times \{a+x \}}=
  {\cal E}_{|Y \times \{a+\lambda'\circ 
  (\lambda \times \mu)(x,\hat{x})\}}$.
  Since $M_H(w) \cong X$, $x=\lambda'\circ 
  (\lambda \times \mu)(x,\hat{x})$.
  Since $p_{X*}({\cal E}nd((id_Y \times T_x)^*{\cal E})) \cong {\cal O}_X$,
 we also get that $\mu'\circ 
  (\lambda \times \mu)(x,\hat{x})=\hat{x}$.
Therefore we have the relation $(\lambda' \times \mu') \circ 
  (\lambda \times \mu)=id_{X \times \widehat{X}}$.
  The other relation also follows from similar arguments.
  Therefore $\lambda \times \mu: X \times \widehat{X} \to Y \times \widehat{Y}$
  is an isomorphism.

\begin{lem}[Mukai \cite{Mu:7}]\label{lem:FM-eq}
For a coherent sheaf $G$ on $X$,
\begin{equation}
 {\cal F}_{\cal E}(T_x^* G \otimes ({\cal P}_X)_{\hat{x}})=
 T_{\lambda(-x,\hat{x})}^* {\cal F}_{\cal E}(G) \otimes
 ({\cal P}_Y)_{\mu(-x,\hat{x})}, (x,\hat{x}) \in X \times \widehat{X}.
\end{equation}
\end{lem}

\begin{proof}
\begin{equation}
 \begin{split}
  {\cal F}_{\cal E}(T_x^* G \otimes ({\cal P}_X)_{\hat{x}}) &=
  {\mathbf R}p_{Y*}({\cal E} \otimes 
  p_X^*(T_x^* G \otimes ({\cal P}_X)_{\hat{x}}))\\
  &={\mathbf R}p_{Y*}(T_{(0,-x)}^*{\cal E} \otimes 
  p_X^*({\cal P}_X)_{\hat{x}}
  \otimes p_X^*G)\\
  &={\mathbf R}p_{Y*}(T_{(\lambda(-x,\hat{x}),0)}^*{\cal E} \otimes
  ({\cal P}_Y)_{\mu(-x,\hat{x})}\otimes p_X^*G)\\
  &=T_{\lambda(-x,\hat{x})}^*{\mathbf R}p_{Y*}({\cal E} \otimes
  p_X^*G) \otimes ({\cal P}_Y)_{\mu(-x,\hat{x})}\\
  &=T_{\lambda(-x,\hat{x})}^* {\cal F}_{\cal E}(G) \otimes
  ({\cal P}_Y)_{\mu(-x,\hat{x})}.
 \end{split}
\end{equation}
\end{proof}

\begin{rem}
Corollary \ref{cor:FM-const} and Lemma \ref{lem:FM-eq}
also hold for any abelian variety.
\end{rem}

Let $K(X)$ be the Grothendieck $K$-group of $X$.
We set 
\begin{equation}
 K(X)_0:=\{\alpha \in K(X)|\text{ $\ch (\alpha)=0$ in $H^*(X,{\Bbb Z})$ }\}.
\end{equation}
For $\alpha \in K(X)$, we define a map
${\frak a}_{\alpha}:K(X)_0+\alpha \to X \times \widehat{X}$
by 
\begin{equation}
 {\frak a}_{\alpha}(\gamma)=(\det({\cal F}_{{\cal P}_X}(\gamma-\alpha)),
 \det(\gamma-\alpha)), \gamma \in K(X)_0+\alpha.
\end{equation}
We define a homomorphism
$\varphi:X \times \widehat{X} \to Y \times \widehat{Y}$
by
$\varphi(x,\hat{x})=
{\frak a}_{\beta}({\cal F}_{\cal E}(I_x \otimes ({\cal P}_X)_{\hat{x}}))$,
$(x,\hat{x}) \in X \times \widehat{X}$,
where $\beta={\cal F}_{\cal E}(I_0)$ and $I_x$ is the ideal sheaf of
$x \in X$.
\begin{lem}
 $\varphi$ is an isomorphism.
\end{lem}

\begin{proof}
Since $I_x \otimes ({\cal P}_X)_{\hat{x}}=
T_{-x}^*(I_0)\otimes ({\cal P}_X)_{\hat{x}}$, $0 \in X$,
Lemma \ref{lem:FM-eq} imlies that
${\cal F}_{\cal E}(I_x \otimes ({\cal P}_X)_{\hat{x}})
\cong T_{\lambda(x,\hat{x})}^*(J) \otimes ({\cal P}_Y)_{\mu(x,\hat{x})}$,
where $J={\cal F}_{\cal E}(I_0)$.
Since $\langle v(J)^2 \rangle=2$,
the homomorphism $Y \times \widehat{Y} \to Y \times \widehat{Y}$ sending
$(y,\hat{y})$ to ${\frak a}(T_y^*(J) \otimes ({\cal P}_Y)_{\hat{y}})$
is an isomorphism (see sect. 4.1), and hence
we get our claim.
\end{proof}
$K(X)_0$ is generated by 
$({\cal P}_X)_{\hat{x}}-{\cal O}_X$, $\hat{x} \in \widehat{X}$
and ${\Bbb C}_{x}-{\Bbb C}_0$, $x \in X$.
Since 
\begin{equation}
 \begin{split}
  I_x \otimes ({\cal P}_X)_{\hat{x}}-I_0 &=
  (I_x \otimes ({\cal P}_X)_{\hat{x}}-({\cal P}_X)_{\hat{x}})
  +({\cal P}_X)_{\hat{x}}-I_0\\
  &=-{\Bbb C}_x+(({\cal P}_X)_{\hat{x}}-{\cal O}_X)+
  ({\cal O}_X-I_0)\\
  &={\Bbb C}_0-{\Bbb C}_x+(({\cal P}_X)_{\hat{x}}-{\cal O}_X),
 \end{split}
\end{equation}
$K(X)_0$ is generated by
$I_x \otimes ({\cal P}_X)_{\hat{x}}-I_0$, $(x,\hat{x}) \in
X \times \widehat{X}$.
Hence we get the following commutative diagram:
\begin{equation}
\begin{CD}
K(X)_0 @>{{\cal F}_{\cal E}}>> K(Y)_0\\
@V{{\frak a}_{\alpha}}VV @VV{{\frak a}_{\beta}}V\\
X \times \widehat{X} @>>{\varphi}> Y \times \widehat{Y}
\end{CD}
\end{equation} 
where $\alpha=I_0$ and $\beta={\cal F}_{\cal E}(\alpha)$.
Therefore we get our main assertion of this subsection.

\begin{prop}\label{Prop:FM-alb}
Assume that ${\cal F}_{\cal E}^i$ induces an isomorphism of
moduli spaces $M_H(v) \to M_{H'}(w)$, where $w=(-1)^i{\cal F}_{\cal E}^H(v)$.
Then
we get the following commutative diagram:
\begin{equation}
\begin{CD}
M_H(v) @>{{\cal F}_{\cal E}^i}>> M_{H'}(w)\\
@V{{\frak a}_v}VV @VV{{\frak a}_w}V\\
X \times \widehat{X} @>>{(-1)^i\varphi}> Y \times \widehat{Y}
\end{CD}
\end{equation} 
In particular, ${\frak a}_v$ is an albanese map if and only
if ${\frak a}_w$ is an albanese map.
\end{prop}

\subsection{Proof of Theorem \ref{thm:deform equiv} and 
\ref{thm:period}}

\subsubsection{Generalized Kummer variety}
 
In this subsection, we shall recall Beauville's results \cite{B:1}
on generalized Kummer varieties.
Then Theorem \ref{thm:period} for $r=1$ follows from his results and
simple calculations.
Let $X$ be an abelian surface.
Let $X^n$ be the $n$-th product of $X$ and $p_i:X^n \to X$
the projection to $i$-th component.
Let $\pi:X^n \to X^{(n)}$
be the quotient map to the $n$-th symmetric product of $X$.
We set $X^{[n]}:=\Hilb_X^n$.
Let $\gamma:X^{[n]} \to X^{(n)}$ be the Hilbert-Chow morphism.
Let $\sigma:X^{(n)} \to X$ be the morphism sending 
$(x_1,x_2,\dots,x_n) \in X^{(n)}$ to $\sum_{i=1}^n x_i \in X$.
Then ${\frak a}:X^{[n]} \to X^{(n)} \to X$ is the albanese map of 
$X^{[n]}$.
If $n=2$, then ${\frak a}^{-1}(0)$ is the
Kummer surface associated to $X$ 
and if $n \geq 3$, then
$K_{n-1}:={\frak a}^{-1}(0)$ is the generalized Kummer variety
constructed by Beauville \cite{B:1}.
$K_{n-1}$ is an irreducible symplectic manifold of
dimension $2(n-1)$.

We assume that $n \geq 3$ and describe $K_{n-1}$ 
up to codimension 2 subscheme.
For integers $i,j,k$, we set
$\Delta^{i,j}:=\{(x_1,x_2,\dots,x_n) \in X^n| x_i=x_j \}$,
$\Delta^{i,j,k}:=\Delta^{i,j} \cap \Delta^{j,k}$.
We set $X_*^n:=X^n \setminus \cup_{i<j<k} \Delta^{i,j,k}$,
$X_*^{[n]}:=X^{[n]} \setminus \cup_{i<j<k}\gamma^{-1}(\pi(\Delta^{i,j,k}))$.
We set 
\begin{equation}
 N:=\{(x_1,x_2,\dots,x_n)| x_1+x_2+\dots+x_n=0 \},
\end{equation}
$N_*:=N \cap X_*^n$, $(K_{n-1})_*:=K_{n-1} \cap X_*^{[n]}$,
$\delta^{i,j}:=\Delta^{i,j} \cap N$.
Since $n \geq 3$,
$\delta^{i,j}$ is connected,
indeed, it is isomorphic to $X^{n-2}$.
Let $\beta:B_{\Delta}(X_*^n) \to X_*^n$ be the blow-up of
$X_*^n$ along $\Delta:=\cup_{i<j} \Delta^{i,j}$
and set $B_{\delta}(N_*)=\beta^{-1}(N_*)$.
Let $E^{i,j}:=\beta^{-1}(\Delta^{i,j})$ 
be the exceptional divisor of $\beta$
and $e^{i,j}:=\beta^{-1}(\Delta^{i,j} \cap N)$.
Let $\Gamma^i \subset B_{\Delta}(X^n_*) \times X$ 
be the graph of the $i$-th projection
$B_{\Delta}(X^n_*) \to X^n_* \to X$.
Then there is a family of ideal sheaves ${\cal I}$ of colength $n$
which fits in an exact sequence
\begin{equation}
0 \to {\cal I} \to {\cal O}_{B_{\Delta}(X^n_*) \times X}
\to \oplus_i {\cal O}_{\Gamma^i} \to 
\oplus_{i<j} {\cal O}_{\Gamma^i \cap E^{i,j}} \to 0,
\end{equation}
and ${\cal I}$ induces a morphism $\pi':B_{\Delta}(X^n_*) \to X^{[n]}_*$.

\begin{equation}
 \begin{CD}
  B_{\delta}(N_*) @>>> B_{\Delta}(X_*^{n}) @>{\beta}>> X^n\\
  @VV{\pi''}V @VV{\pi'}V @VV{\pi}V \\
  (K_{n-1})_* @>>> X_*^{[n]} @>>{\gamma}> X^{(n)}\\
  @.@.@VV{\sigma}V\\
  @.@.X
 \end{CD}
\end{equation}
Since $K_X \cong {\cal O}_X$, by using Grothendieck-Riemann-Roch theorem
for the embedding $\cup_i \Gamma^i \hookrightarrow 
B_{\Delta}(X_*^{n}) \times X$,
we see that
\begin{equation}\label{eq:family}
 \ch({\cal I})
 =1-\sum_{i}\Gamma^i+\sum_{i<j}E^{i,j}\Gamma^i+\cdots,
\end{equation}
where, by Poincar\'{e} duality, we identifies  
$H^*(B_{\Delta}(X_*^{n}),{\Bbb Q})$ with Borel-Moore homology group
$H_*(B_{\Delta}(X_*^{n}),{\Bbb Q})$ (cf. \cite{Fu:1}).
We set $v:=1-n \omega$. 
Then $M_H(v)$ is naturally identified with $X^{[n]} \times \widehat{X}$ 
and ${\frak a} \times 1_{\widehat{X}}:X^{[n]} \times \widehat{X}
\to X \times \widehat{X}$ is identified with the map ${\frak a}_v$
introduced in section 4.1.
Hence $K_H(v)$ is identified with $K_{n-1}$.
We shall consider the homomorphism
$\theta_v:v^{\perp} \to H^2(K_{n-1},{\Bbb Z})$.
Since $v=1-n \omega$,
we get that
\begin{equation}
 v^{\perp}=H^2(X,{\Bbb Z})\oplus {\Bbb Z}(1+n \omega).
\end{equation}
For $\alpha=x+k(1+n \omega),
x \in H^2(X,{\Bbb Z})$,
simple calculations show that
\begin{equation}
 (\pi'')^*\theta_v(\alpha)=
 -\left(\sum_i p_i^* (x)_{|B_{\delta}(N_*)}+k e \right),
\end{equation}
where $e:=\sum_{i<j} e^{i,j}$.

\begin{lem}\label{lem:beauville}
$\theta_v$ is an isomorphism.
In particular, we get identifications
\begin{equation}
 \begin{split}
  H^2(K_{n-1},{\Bbb Z})  & \cong
  H^2(B_{\delta}(N_*),{\Bbb Z})^{{\frak S}_n}\\
  & \cong H^2(X,{\Bbb Z}) \oplus {\Bbb Z}e.
 \end{split}
\end{equation}
\end{lem}
\begin{proof}
By \cite[Prop. 8]{B:1},
$H^2(K_{n-1},{\Bbb Q}) \cong H^2(X,{\Bbb Q}) \oplus {\Bbb Q}\,e$.
Since 
\begin{equation}
 \varphi:H^2(X,{\Bbb Z}) \overset{\theta_v} \to H^2(K_{n-1},{\Bbb Z})
 \to H^2(B_{\delta}(N_*),{\Bbb Z})
\end{equation}
 is injective and
$\im \varphi \subset \beta^*(H^2(N,{\Bbb Z}))$,
we shall prove that the image of 
$f:H^2(X,{\Bbb Z}) \to H^2(N,{\Bbb Z})^{{\frak S}_n}$ is a
primitive submodule of
$H^2(N,{\Bbb Z})$.
Let $\phi:X \times X \to N$ be the morphism such that 
\begin{equation}
 \phi((x,y))=(x,y,0,\dots,0,-x-y) \in N.
\end{equation}
We shall consider the composition $g=\phi^* \circ f
:H^2(X,{\Bbb Z}) \to 
H^2(X \times X,{\Bbb Z})$.
Let $\alpha_i \wedge \alpha_j,i<j$ be the basis of 
$H^2(X,{\Bbb Z})=\wedge^2 H^1(X,{\Bbb Z})$.
Then we see that
\begin{equation}
 g^*(\alpha_i \wedge \alpha_j)=2 p_1^*(\alpha_i \wedge \alpha_j)+
 2p_2^*(\alpha_i \wedge \alpha_j)+(p_1^*\alpha_i \wedge p_2^*\alpha_j
 -p_1^* \alpha_j \wedge p_2^*\alpha_i).
\end{equation}
Hence $\im g$ is a primitive subgroup of $H^2(X \times X,{\Bbb Z})$.
Therefore $\im f$ is primitive.
\end{proof}

Let $\iota:H^2(X,{\Bbb Z}) \to H^2(X^{(n)},{\Bbb Z})$ be the homomorphism 
such that $\pi^*(\iota(x))=\sum_i p_i^* x \in H^2(X^n,{\Bbb Z})$. 
For $x \in H^2(X,{\Bbb Z})$, we get
\begin{equation}
 \begin{split}
  (\pi'')^* \theta_v(x) &=-\sum_i p_i^*(x)_{|B_{\delta}(N_*)}\\
  &=-(\pi^* \iota(x))_{|B_{\delta}(N_*)}\\
  &=-(\pi'')^* (\gamma^* \iota(x))_{|(K_{n-1})_*}.
 \end{split}
\end{equation}
Thus 
\begin{equation}
 \theta_v(x)=-\gamma^*(\iota(x))_{|K_{n-1}}, x \in H^2(X,{\Bbb Z}).
\end{equation}

We shall next prove that $\theta_v$ preserves the bilinear forms.
Since $\langle \alpha^2  \rangle=(x^2)-k^2 (2n)$
for $\alpha=x+k(1+n \omega)$,
 it is sufficient to prove 
\begin{equation}\label{eq:q}
 q_{K_{n-1}}(\theta_v(\alpha))=(x^2)-k^2 (2n).
\end{equation}
We shall choose $l \in H^2(X,{\Bbb Z})$ with $(l^2) \ne 0$.
Let $\overline{B_{\delta}(N_*)}$ be a smooth compactification of
$B_{\delta}(N_*)$ such that there are extensions
$\overline{\beta}:\overline{B_{\delta}(N_*)} \to N$ and
$\overline{\pi''}:\overline{B_{\delta}(N_*)} \to K_{n-1}$
of $\beta$ and $\pi''$ respectively.
 
(1) We shall first consider the relation between 
$q_{K_{n-1}}(\theta_v(x)), x \in H^2(X,{\Bbb Z})$ and 
$q_{K_{n-1}}(\theta_v(l))$.
We set $\rho:=\alpha_1 \wedge \alpha_2 \wedge \alpha_3 \wedge \alpha_4$,
where $\alpha_i$, $1 \leq i \leq 4$ are basis of
$H^1(X,{\Bbb Z})$.
Let 
\begin{equation}
 \eta:\otimes_{i=1}^n p_i^* H^*(X,{\Bbb Z}) \to 
 \otimes_{i=1}^n p_i^* H^{ev}(X,{\Bbb Z})
\end{equation}
be the projection.
We shall first compute $\eta((\sigma \circ \pi)^* (\rho))$.
By direct computations, we see that
\begin{equation}
 \begin{split}
  \eta((\sigma \circ \pi)^* (\rho)) &=
  \eta(\sum_{i,j,k,m}p_i^* \alpha_1
  \wedge p_j^* \alpha_2 \wedge p_k^* \alpha_3 \wedge p_m^*\alpha_4)\\
  &=\sum_i p_i^*(\alpha_1 \wedge \alpha_2 \wedge \alpha_3 \wedge \alpha_4)
  +\sum_{i \ne j}p_i^*(\alpha_1 \wedge \alpha_2)
  \wedge p_j^*(\alpha_3 \wedge \alpha_4).
 \end{split}
\end{equation}
Since $\eta(\sum_{i<j}(\Delta^{i,j}
-p_i^* \rho -p_j^* \rho))=\sum_{i \ne j}p_i^*(\alpha_1 \wedge \alpha_2)
 \wedge p_j^*(\alpha_3 \wedge \alpha_4)$,
\begin{equation}
 \eta((\sigma \circ \pi)^* (\rho))=
 \sum_i p_i^* \rho+\eta(\sum_{i<j}(\Delta^{i,j}
 -p_i^* \rho -p_j^* \rho)).
\end{equation}
We set $\mu:=\sum_i p_i^* \rho+\sum_{i<j}(\Delta^{i,j}
 -p_i^* \rho -p_j^* \rho)$.
Then we see that
\begin{equation*}
 \begin{split}
  \int_{K_{n-1}}(\theta_v(l))^{2n-4}(\theta_v(x))^2 &=
  \frac{1}{n!} \int_{\overline{B_{\delta}(N_*)}}
  (\overline{\pi''}^*\theta_v(l))^{2n-4}
  (\overline{\pi''}^*\theta_v(x))^2\\
  &=\frac{1}{n!} \int_N (p_1^* l+\dots+ p_n^*l)^{2n-4}
  (p_1^*x+\dots +p_n^* x)^2\\
  &=\frac{1}{n!} \int_{X^n} (p_1^* l+\dots +p_n^*l)^{2n-4}
  (p_1^*x+\dots +p_n^* x)^2 \mu\\
  &=\frac{1}{n!}\left\{\frac{n(n-1)}{2}
  \int_{X^{n-1}} (2p_1^* l+p_2^*l+\dots +p_{n-1})^{2n-4}
  (2p_1^*x+p_2^*x+\dots +p_{n-1})^2 \right.\\
  & \quad \left.-n(n-2)\int_{X^{n-1}} (p_1^* l+p_2^*l+\dots +p_{n-1}^*l)^{2n-4}
  (p_1^*x+p_2^*x+\dots +p_{n-1}^* x)^2
  \right\}\\
  &=\frac{(2n-2)!n^2}{n!2^{n-1}}\left(
  \frac{1}{2n-3}(l^2)^{n-2}(x^2)+\frac{2n-4}{2n-3}(l^2)^{n-2}(l,x)^2 
  \right).
 \end{split}
\end{equation*}
In the same way, we see that
\begin{equation}
 \begin{split}
  \int_{K_{n-1}}(\theta_v(l))^{2n-3}(\theta_v(x)) &=
  \frac{(2n-2)!n^2}{2^{n-1}n!}(l^2)^{n-1}(l,x),\\
  \int_{K_{n-1}}(\theta_v(l))^{2n-2}&=
  \frac{(2n-2)!n^2}{2^{n-1}n!}(l^2)^{n}.
 \end{split}
\end{equation}
By \eqref{eq:Bform}, we obtain that
\begin{equation}\label{eq:D}
 q_{K_{n-1}}(\theta_v(x))=\frac{(x^2)}{(l^2)}q_{K_{n-1}}(\theta_v(l)).
\end{equation}

(2) We shall next consider $q_{K_{n-1}}(\theta_v(1+n \omega))$.
For simplicity, we set $\widetilde{e}:=-\theta_v(1+n \omega)$.
Then $(\pi'')^*(\widetilde{e})=e$.
Hence we also denote $\overline{\pi''}^*(\widetilde{e})$ by $e$.
Since $E^{i,j}$ is the exceptional divisor
of $\beta$, $\beta_*(E^{i,j})=0$ and 
$\beta_*((E^{i,j})_{|E^{i,j}})=-\Delta^{i,j}$.
Since $\codim_{X^n}(\cup_{i<j<k} \Delta^{i,j,k})=4$, we get
\begin{equation}
 \begin{split}
  \overline{\beta}_*(e)&=0,\\
  \overline{\beta}_*(e^2)&=-\sum_{i<j}\Delta^{i,j}.
 \end{split}
\end{equation}

Let $\upsilon:X^{n-1}=\Delta^{1,n} \to X$ be the restriction
of $\sigma \circ \pi$ to the diagonal $\Delta^{1,n}$. 
Then,
\begin{equation}
 \upsilon^*(\rho)=\left(2 p_1^* \alpha_1+\sum_{i=2}^{n-1}p_i^* \alpha_1\right)
 \wedge \dots 
 \wedge \left(2 p_1^* \alpha_4+\sum_{i=2}^{n-1}p_i^* \alpha_4 \right).
\end{equation}
We set
\begin{equation}
 \begin{split}
  \mu':=& 2^4 p_1^* \rho+\sum_{i=2}^{n-1}p_i^* \rho+
  \sum_{1<i<j \leq n-1}
  (\Delta^{i,j}-p_i^* \rho-p_j^* \rho)\\
  & \quad +4\sum_{i=2}^{n-1}(\Delta^{1,i}-p^*_1 \rho-p^*_i\rho).
 \end{split}
\end{equation}
In the same way as above, we see that the 
 $\otimes_{i=1}^{n-1}p^*_i H^{ev}(X,{\Bbb Z})$-components of
$\upsilon^* (\rho)$ and $\mu'$ 
are the same.
Hence we see that 

\begin{equation}
 \begin{split}
  \int_{K_{n-1}}(\theta_v(l))^{2n-4}\theta_v(x) \widetilde{e}=& 0,\\
  \int_{K_{n-1}}(\theta_v(l))^{2n-4} \widetilde{e}^2 =&
  \frac{1}{n!} \int_{\overline{B_{\delta}(N_*)}}
  (\overline{\pi''}^*\theta_v(l))^{2n-4}
  (e^2)\\
  =&-\frac{1}{n!} \int_N
  (p_1^* l+\dots+ p_n^*l)^{2n-4}
  (\sum_{i<j} \Delta^{i,j})\\
  =&-\frac{n(n-1)}{2 n!} \int_{X^{n-1}} 
  (2p_1^* l+p_2^*l+\dots +p_{n-1}^*l)^{2n-4}
  \mu'\\
  =&\frac{(2n-2)! n^2}{2^{n-1}n!}\frac{(-2n)}{2n-3}.
 \end{split}
\end{equation}
Thus $\widetilde{e}$ is orthogonal to $H^2(X,{\Bbb Z})$ and
\begin{equation}\label{eq:e}
 q_{K_{n-1}}(\theta_v(1+n \omega))=\frac{-2n}{(l^2)}q_{K_{n-1}}(\theta_v(l)).
\end{equation}
By the identification $\theta_v:v^{\perp} \to H^2(K_{n-1},{\Bbb Z})$, 
$\langle \quad, \quad \rangle$ is a primitive bilinear
form on $H^2(K_{n-1},{\Bbb Z})$.
Then, by the definition of $q_{K_{n-1}}$,
\eqref{eq:D} and \eqref{eq:e} imply that
\eqref{eq:q} holds. 

\begin{prop}[Beauville]\label{prop:beauville}
For a Mukai vector $v$ of $\rk v=1$ and
$\langle v^2 \rangle \geq 6$,
\begin{equation}
 \theta_v:(v^{\perp},\langle\;\;,\;\;\rangle)
 \to (H^2(K_H(v),{\Bbb Z}),B_{K_H(v)})
\end{equation}
is an isometry of Hodge structures.
\end{prop}
\begin{proof}
We set $v=1+\xi+a \omega, \xi \in \NS(X)$
and $n:=\langle v^2 \rangle/2$.
Then we see that
$(1-n \omega)\exp( \xi)=v$.
By \eqref{eq:q},
$\theta_{1-n \omega}$ is an isometry of Hodge structure.
Hence \eqref{eq:isomT} implies our claim.
\end{proof}

Therefore Theorem \ref{thm:period} holds for the rank 1 case.
Obviously, Theorem \ref{thm:deform equiv} also holds. 

\subsubsection{General cases}

\begin{prop}\label{prop:deform}
Let $X_1$ and $X_2$ be abelian (or K3) surfaces, and let
$v_1:=r+\xi_1+a_1 \omega \in H^{ev}(X_1,{\Bbb Z})$
and $v_2:=r+\xi_2+a_2 \omega \in H^{ev}(X_2,{\Bbb Z})$
be primitive Mukai vectors such that
$(1)$ $r>0$, $(2)$ $\ell(v_1)=\ell(v_2)=l$, 
$(3)$ $\langle v_1^2 \rangle=\langle v_2^2 \rangle=2s$,
and $(4)$ $a_1 \equiv a_2 \mod l$.
Then $M_{H_1}(v_1)$ and $M_{H_2}(v_2)$ are deformation equivalent,
where $H_i$, $i=1,2$ are ample divisors on $X_i$ 
such that $\overline{M}_{H_i}(v_i)=M_{H_i}(v_i)$.
Moreover if $X_i$, $i=1,2$ are abelian surfaces, then
${\frak a}_{v_1}:M_{H_1}(v_1) \to X_1 \times \widehat{X_1}$
is deformation equivalent to 
${\frak a}_{v_2}:M_{H_1}(v_2) \to X_2 \times \widehat{X_2}$.
In particular, $K_{H_1}(v_1)$ is an irreducible symplectic manifold
and $\theta_{v_1}$ is an isometry of Hodge structures if and only if 
$K_{H_2}(v_2)$ and $\theta_{v_2}$ have the same properties. 
If $\xi_i$ are ample and $a_1=a_2 \ne 0$, 
then the same assertions also hold for $r=0$.
\end{prop}

\begin{proof}
If $r>0$ and $H_i$ are general, then the results follow from the arguments of
O'Grady \cite{O:1} (cf. \cite[Prop. 1.1]{Y:5}).
Indeed, let ${\cal M}:=\{(X,L)\}$ 
be the moduli space of polarized abelian surfaces
of $(c_1(L)^2)=2n$.
Then ${\cal M}':=\{(X,L) \in {\cal M}|\rho(X) \geq 2\}$ is infinite 
(countable) union of
algebraic subsets (cf. \cite[Exercise 10.10]{L-B:1}).
We can also find a {\it suitable} polarization for
a product of elliptic curves.
However, in order to treat the last statement, we also need the method of
G\"{o}ttsche and Huybrechts \cite{G-H:1}.
Hence we do not use O'Grady's arguments here. 

We first assume that $r>0$.
We note that the stability does not change under the
operation $E \mapsto E(n_i H_i)$.
Hence by \eqref{eq:isomT}, we may replace $v_i$ by 
$v_i\ch(H_i^{\otimes n_i})$, $n_i \gg 0$.
Thus we may assume that $\xi_1$ and $\xi_2$ are ample.
By using Proposition \ref{prop:modification} in Appendix,
we may assume that $H_i$, $i=1,2$ are general with respect to $v_i$
(i.e. $(\natural)$ in sect. 1.2 holds).
By using Proposition \ref{prop:modification} again, 
we may assume that $\rho(X_i) \geq 2$ and $H_i$ are general.
Replacing $v_i$ by $v_i \ch(N_i)$, $N_i \in \Pic(X_i)$,
we may assume that $\xi_i$ is a primitive ample class of
$(\xi_i^2) \geq 4$ (we use Lemma \ref{lem:isomT}).
By the connectedness of the moduli space of polarized K3 (or abelian)
surfaces and Proposition \ref{prop:modification},
we may assume that (1) $X_1=X_2$, (2) $X:=X_1$ has an elliptic fibration
$\pi:X \to C$, (3)
$H_i$, $i=1,2$ are general with respect to $v_i$
and (4)  
$\xi_i/l=\sigma+n_i f$, where $\sigma$ is a section of $\pi$,
$f$ a fiber of $\pi$ and $l=\ell(v_1)(=\ell(v_2))$.
Since $\langle v_1^2 \rangle=\langle v_2^2 \rangle=2s$,
we see that 
$v_1 \exp((a_2-a_1)f/l)=v_2$.
Hence if $H_1=H_2$, then Lemma \ref{lem:isomT} implies our claims.
If $H_1 \ne H_2$, then we take a family of polarized K3 (or abelian) surfaces
$({\cal X},{\cal L}) \to T$ such that
$({\cal X}_{t_0},{\cal L}_{t_0})=(X,\xi_1)$ and 
$\rho({\cal X}_{t_1})=1$ for $t_0,t_1 \in T$.
Applying Proposition \ref{prop:modification} again, 
we can reduce to the case where $H_1=H_2$.
Therefore we get our claim.
If $r=0$, then Proposition \ref{prop:modification} implies 
our assertions.
\end{proof}
{\it Proof of Theorem \ref{thm:deform equiv} and \ref{thm:period}.}
(I) We first consider the case where $\ell(v)=1$. 
Assume that $\rk v>0$.
Let $r$ and $s$ be positive integers. 
We shall find an abelian surface $X$
and a Mukai vector $v=r+\xi+a \omega
\in H^{ev}(X,{\Bbb Z})$ of $\langle v^2 \rangle =2s$
which satisfy the claims of 
Theorem \ref{thm:deform equiv} and \ref{thm:period}.
Let $(X,H)$ be a polarized abelian surface of
$\NS(X)={\Bbb Z}H$ and $(H^2)=2r+2s$.
We set $v=r+H+\omega$.
Since $\langle v^2 \rangle=(H^2)-2r=2s$, 
we shall prove the claims for this $v$.
Let $\widehat{X}$ be the dual abelian surface of $X$ and
${\cal P}$ the Poincar\'{e} line bundle on $X \times \widehat{X}$.
By Proposition \ref{lem:fourier4},
${\cal G}_{\cal P}$ induces an isomorphism
$M_H(v) \to M_{\widehat{H}}(w)$, where $w=1+\widehat{H}+r \widehat{\omega}$.
Then we have a commutative diagram:
\begin{equation}\label{eq:alb}
 \begin{CD}
  M_H(v) @>{{\cal G}_{\cal P}}>> M_H(w)\\
  @V{{\frak a}_v}VV @VV{{\frak a}_w}V\\
  X \times \widehat{X} @= X \times \widehat{X}
 \end{CD}
\end{equation}
where we assume that ${\frak a}_v(E_0)={\frak a}_w({\cal G}^1_{\cal P}(E_0))
=0$.
Hence under ${\cal G}_{\cal P}$, we can identify $K_H(v)$ with
$K_{\widehat{H}}(w)$.
Then, Proposition \ref{prop:comm2} implies that
the following diagram is commutative.
\begin{equation}
 \begin{CD}
  v^{\perp} @>{{\cal G}_{\cal P}^H}>> w^{\perp}\\
  @V{\theta_v}VV @VV{\theta_w}V\\
  H^2(K_H(v),{\Bbb Z}) @= H^2(K_H(w),{\Bbb Z})
 \end{CD}
\end{equation} 
Since $K_{\widehat{H}}(w)$ is an irreducible symplectic manifold and
$\theta_w$ is an isometry of Hodge structures
(Proposition \ref{prop:beauville}),
Theorem \ref{thm:deform equiv} and \ref{thm:period} hold 
for $M_H(v)$, if $\rk v>0$.
Assume that $\rk v=0$.
We set $v=\xi+a \omega$, $\xi \in \NS(X)$.
By our assumption on $v$, $\xi$ is effective (or $\xi=0$).
Since $(\xi^2)=\langle v^2 \rangle>0$,
\cite[Chap. 4 Prop. 5.2]{L-B:1} implies that $\xi$ is ample.
Hence by Proposition \ref{prop:deform},
it is sufficient to prove Theorem \ref{thm:deform equiv} and \ref{thm:period}
for $M_H(H+a \omega)$, where $X$ is an abelian surface of
$\NS(X)={\Bbb Z}H$ and $a \ne 0$.
By Propositions \ref{lem:fourier3} or \ref{lem:fourier4},
we get an isomorphism $M_H(H+a \omega) \cong
M_{\widehat{H}}(|a|+\widehat{H})$.
Since Theorem \ref{thm:deform equiv} and \ref{thm:period} hold if $\rk v>0$,
Theorem \ref{thm:deform equiv} and \ref{thm:period} hold for the case where $\rk v=0$.

(II) We next treat general cases, that is, we shall reduce the general cases
to $\ell(v)=1$ case.
Let $X$ be an abelian surface and  
$v=l(r+c_1)+a \omega$, $c_1 \in \NS(X)$ a primitive Mukai vector such that
$l=\ell(v)$ and $2ls:=\langle v^2 \rangle>0$.

(II-1) We first assume that $\rk v>0$.
We take a positive integer $k$ such that
$n:=rk-(c_1^2)/2>0$.
Since $(l,a)=1$, we may assume that  
$b:=a-kl$ and $r$ are relatively prime.
Let $C$ be an elliptic curve which has an endomorphism 
$\phi:C \to C$ of $\ker \phi \cong {\Bbb Z}/n{\Bbb Z}$.
We set $Y:=C \times C$.
Let $\pi:Y \to C$ be the first projection, $f$ a fiber of $\pi$ and
$\sigma$ the $0$-section of $\pi$.  
Then the graph $\Gamma_{\phi} \subset Y$ satisfies that
$(\Gamma_{\phi},\sigma)=n$.
Hence $\xi_n:=\Gamma_{\phi}-\sigma-nf$ satisfies that
$(\xi_n^2)=-2n$.
We set $w:=l(r+(-\xi_n+f))+b \omega$.
By Remark \ref{rem:general},
$\overline{M}_{\sigma+kf}(w)=M_{\sigma+kf}(w)$, $k \gg 0$.
In particular, it is compact.
Since $\langle w^2 \rangle =2ls$,
by Proposition \ref{prop:deform},
it is sufficient to prove Theorem \ref{thm:deform equiv} and \ref{thm:period}
for $M_{\sigma+kf}(w)$, $k \gg 0$.
By Theorem \ref{thm:FM}, we have an isomorphism
\begin{equation}
 M_{\sigma+kf}(lr+l(-\xi_n+f)+b\omega) \cong 
 M_{\sigma+kf}(l(\xi_n+r \sigma)-bf+l(r+1)\omega), k \gg 0.
\end{equation}
Since $D:=l(\xi_n+r \sigma)-bf$ satsifies $(D^2)=2ls>0$
and $H^0(Y,{\cal O}_Y(-D))=0$, Riemann-Roch theorem implies that
$D$ is effective.
By \cite[Prop. 5.2]{L-B:1}, $D$ is an ample divisor.
Since $(b,r)=1$, $D$ is a primitive ample divisor.
Since we proved our theorem for the case where $\ell(v)=1$,
Theorem \ref{thm:deform equiv} and \ref{thm:period} hold for
 $M_{\sigma+kf}(l(\xi_n+r \sigma)-bf+l(r+1)\omega)$, $k \gg 0$.
Hence by taking account of Proposition \ref{Prop:FM-alb},
we see that Theorem \ref{thm:deform equiv} and \ref{thm:period} hold for
$M_{\sigma+kf}(w)$, $k \gg 0$.

(II-2) We next assume that $\rk v=0$.
In this case, we use an isomorphism
\begin{equation}
 M_{\sigma+kf}(l+(b-l)f-al\omega) \cong 
 M_{\sigma+kf}(l(\sigma+af)+b\omega), k \gg 0
\end{equation}
of moduli spaces on the elliptic surface $\pi:Y=C \times C \to C$.
Then we can reduce our problem to the case where $\rk v>0$.
\qed

{\it Proof of Theorem \ref{thm:H2}.}
By direct computations, we can also show that
Theorem \ref{thm:H2} holds for $r=1$.
Then, by similar method as in the proof of Theorem \ref{thm:period},
we can prove Theorem \ref{thm:H2}:
In order to prove Theorem \ref{thm:H2} (3), we use Proposition
\ref{Prop:FM-alb}.
\qed
\newline
%

\begin{rem}
In the case where $\ell(v)=1$, 
Dekker \cite{D:1} proved Theorem \ref{thm:deform equiv}
by using Fourier-Mukai functor on a product of elliptic curves.
However, there may be some gaps in his proof.
For example, the second line of the proof of
\cite[Thm. 4.9]{D:1} may not be clear and 
the proof of \cite[Thm. 4.11]{D:1} is not true.
Also the arguments in \cite[sect. 5]{D:1} are not correct.
It might be interesting to justify his arguments.
\end{rem}

\begin{cor}\label{cor:NS}
We set $(v^{\perp})_{alg}:=v^{\perp} \cap
(H^0(X,{\Bbb Z}) \oplus \NS(X)\oplus
H^4(X,{\Bbb Z}))$.
Then $\theta_v$ induces an isometry 
$$
(v^{\perp})_{alg} \to \NS(K_H(v)).
$$
\end{cor}
The following example is similar to \cite[5.17]{Mu:6}. 

\begin{ex}\label{ex:ex1}
Let $X$ be an abelian surface with $\NS(X)={\Bbb Z}H$,
$(H^2)=2$.
We set $v=2+H-2 \omega$.
Then $M_H(v)$ is a variety of dimension $12$.
It is easy to see that
$v^{\perp}$ is generated by
$\alpha:=1+\omega$ and $\beta:=H+\omega$.
Since
\begin{alignat}3
 \langle \alpha^2 \rangle=-2, \quad & \langle \alpha,\beta \rangle=-1, 
 \quad & \langle \beta^2 \rangle=2,
\end{alignat} 
$\NS(K_H(v))$ is indecomposable.
Hence $M_H(v)$ is not birationally equivalent to
$\widehat{Y} \times \Hilb_Y^{5}$ for any $Y$.
\end{ex}

\begin{rem}\label{rem:stable}
Combining the dimension counting in \cite{Y:8}, 
we can show the following (the second claim is due to Mukai \cite{Mu:5}).
\begin{itemize}
\item
Assume that $v$ is primitive, $v>0$ and $\langle v^2 \rangle >0$. 
Then for a general $H$,
\begin{enumerate}
\item
$M_H(v)$ contains $\mu$-stable vector bundle, unless $v=rv(L)-\omega$, where 
$L$ is a line bundle on $X$.
\item
$M_H(rv(L)-\omega)$ consists of non-locally free sheaves and is isomorphic
to $X \times \Hilb_{\widehat{X}}^r$.
\end{enumerate}
\end{itemize}
\end{rem}

\subsubsection{Application to Fourier-Mukai functor }

\begin{thm}\label{thm:spl=stable}
Let $v$ be a primitive Mukai vector such that $c_1(v) \in \NS(X)$, $v>0$
and $\langle v^2 \rangle >0$.
Let $S(v)$ be an irreducible open subscheme of the moduli space
of simple sheaves $\Spl(v)$.
Assume that the albanese manifold of $S(v)$ is $\widehat{X} \times X$.
Then $S(v)$ is birationally equivalent to $M_H(v)$, or
a general element $E$ of $S(v)$ fits in an exact sequence
\begin{equation}
 0 \to E_1 \to E \to E_2 \to 0
\end{equation}
where $v(E_1)=l_1w_1$, $v(E_2)=l_2w_2$ and 
$\langle w_1^2 \rangle =\langle w_2^2 \rangle=0$,
$\langle w_1,w_2 \rangle=1$, $(l_1-1)(l_2-1)=0$.

In particular, if $\NS(X) \cong {\Bbb Z}$, then
$S(v)$ is birationally equivalent to $M_H(v)$.
\end{thm}

\begin{proof}
Assume that $S(v)$ does not contain a semi-stable sheaf.
Let $Q(v)$ be an open subscheme of a suitable quot scheme
$\Quot_{{\cal O}_X^{\oplus N_v}/X/{\Bbb C}}$
such that $S(v)$ is a birational quotient of $Q(v)$ by
$GL(N_v)$:
There is an open subscheme $Q(v)'$ of $Q(v)$
and a $GL(N_v)$-invariant morphism $Q(v)' \to S(v)$ such that 
$Q(v)'/GL(N_v)$ is birationally equivalent to $S(v)$.
For a sequence of Mukai vectors
$v_1,v_2,\dots,v_s$,
let $F(v_1,v_2,\dots,v_s)$ be the set of $q \in Q(v)$ such that
the Harder-Narasimhan filtration of ${\cal E}_q$:
\begin{equation}
 0 \subset F_1 \subset F_2 \subset \dots \subset F_s={\cal E}_q
\end{equation}
satisfies $v(F_i/F_{i-1})=v_i$.
Since $Q(v)$ does not contain a semi-stable sheaf,
there is a sequence of $v_1,v_2,\dots,v_s$ such that
$F(v_1,v_2,\dots,v_s)$ is an open dense subscheme of $Q(v)$.
Let $\zeta:F(v_1,v_2,\dots,v_s) \to \prod_i \overline{M}_H(v_i)$
be a morphism sending $q \in F(v_1,v_2,\dots,v_s)$ to
$([F_1],[F_2/F_1],\dots,[F_s/F_{s-1}]) \in \prod_i \overline{M}_H(v_i)$,
where $[F_i/F_{i-1}]$ is the $S$-equivalence class of $F_i/F_{i-1}$.
By \cite[sect. 5.2]{Y:8}, we see that $\zeta$ is dominant.
Composing $\zeta$ with $\prod_i{\frak a}_{v_i}$,
we get a morphism $F(v_1,v_2,\dots,v_s) \to (\widehat{X} \times X)^s$.
Obviously, this map is $GL(N_v)$-invariant.
Hence we get a morphism $\eta:S(v) \to (\widehat{X} \times X)^s$.
By Proposition \ref{prop:dim-alb}, $\dim \im (\eta) \geq 2s$.
Since the albanese manifold of $S(v)$ is $\widehat{X} \times X$,
by the universal property of the albanese map,
we get that $s \leq 2$.
Moreover by the second claim of Proposition \ref{prop:dim-alb},
we get that $\langle v_1^2 \rangle=\langle v_2^2 \rangle=0$.
Hence a general member $E$ of $S(v)$ fits in an exact sequence
\begin{equation}
 0 \to E_1 \to E \to E_2 \to 0,
\end{equation}
where $v(E_i)=v_i$.
We set $v_i:=l_iw_i$, where $w_i$ are primitive.
Then $\langle v^2 \rangle=2l_1l_2n$, where $n=\langle w_1,w_2 \rangle$.
Since $E$ is simple, $\Hom(E_2,E_1)=0$.
Hence $n=\langle w_1,w_2 \rangle>0$.
We denote the moduli stack of semi-stable sheaves $E$ of $v(E)=v_i$ by
${\cal M}_H(v_i)^{ss}$.
Then we get that $\dim S(v)-1=
\dim {\cal M}_H(v_1)^{ss}+\dim {\cal M}_H(v_2)^{ss}
+\langle v_1,v_2 \rangle$.
Since $\dim {\cal M}_H(v_i)^{ss}=l_i$ (see \cite[Lem. 1.8]{Y:8}), 
we have $l_1+l_2+l_1l_2n=2l_1l_2n+1$,
which implies that $n=1$ and $(l_1-1)(l_2-1)=0$.
Assume that $\NS(X)={\Bbb Z}H$.
We set $w_i:=r_i+d_iH+a_i \omega$, $i=1,2$.
If $\langle w_i \rangle=d_i^2(H^2)-2r_ia_i=0$, then
we get that $\langle w_1,w_2 \rangle=-(r_2d_1-r_1d_2)^2(H^2)/2 r_1r_2<0$.
Hence $S(v)$ is birationally equivalent to $M_H(v)$. 
\end{proof}

\begin{cor}\label{cor:wit->birat}
Let $X$ be an abelian surface of $\NS(X)={\Bbb Z}H$. 
Let $v$ be a primitive Mukai vector such that $c_1(v) \in \NS(X)$, $v>0$
and $\langle v^2 \rangle >0$.
Let ${\cal F}_{\cal P}:{\mathbf D}(X) \to {\mathbf D}(Y)$ 
(resp. ${\cal G}_{\cal P}:{\mathbf D}(X) \to {\mathbf D}(Y)_{op}$) 
be a Fourier-Mukai functor
associated to a universal family on $X \times Y$.
Assume that $\WIT_i$ holds for some $E \in M_H(v)$.
Then ${\cal F}_{\cal P}$ (resp. ${\cal G}_{\cal P}$)
induces a birational map
$M_H(v) \cdots \to M_{\widehat{H}}(w)$, where
$w=(-1)^i{\cal F}^H_{\cal P}(v)$ (resp. $w=(-1)^i{\cal G}^H_{\cal P}(v)$).
\end{cor}

\begin{proof}
By Theorem \ref{thm:H2},
$\Alb(M_H(v)) \cong \widehat{X} \times X \cong \widehat{Y} \times Y$.
Hence we can apply Theorem \ref{thm:spl=stable} to
${\cal F}_{\cal P}^i(M_H(v))$.
\end{proof}

\begin{guess}\label{conj:1}
Assume that ${\cal P}$ is the Poincar\'{e} line bundle on 
$X \times \widehat{X}$. 
If $\NS(X)={\Bbb Z}H$, then ${\cal F}_{\cal P}$ or
${\cal G}_{\cal P}$ induces a birational map
$M_H(v) \cdots \to M_{\widehat{H}}(w)$,
where $w=\pm{\cal F}_{\cal P}^H(v)$ or $w=\pm{\cal G}_{\cal P}^H(v)$.
\end{guess}

For an evidence, we can show the following theorem.

\begin{thm}\label{thm:evidence}
Assume that $\NS(X)={\Bbb Z}H$ and 
${\cal P}$ is the Poincar\'{e} line bundle on $X \times 
\widehat{X}$.
Let $r+dH+a \omega$ be a primitive Mukai vector such that
$r>0$ and $d \geq 0$.  
\begin{enumerate}
\item
If $a \leq 0$ and $d>0$, then
${\cal F}_{\cal P}$ induces a birational map
$M_H(r+dH+a \omega) \cdots \to 
M_{\widehat{H}}(-a+d\widehat{H}-r \widehat{\omega})$.
\item
If $d=0$ or $0<a \leq 4$, then 
${\cal G}_{\cal P}$ induces a birational map
$M_H(r+dH+a \omega) \cdots \to 
M_{\widehat{H}}(a+d \widehat{H}+r \widehat{\omega})$.
\end{enumerate}
\end{thm}
  
The proof will be done in section 6.

\subsection{Non-K\"{a}hler symplectic manifold}

For another application of Theorem \ref{thm:period},
we shall give an example of non-K\"{a}hler compact symplectic manifold,
which is a symplectic analogue of Hironaka's example \cite{Hi:1}.
Let $X$ be an abelian surface of $\NS(X)={\Bbb Z}H$ and assume that 
$v:=r+c_1+a \omega$ satisfies that $\ell(v)=1$ and
$\langle v^2 \rangle=2r$.
We set
\begin{equation}
 M_H(v)_s:=\{E \in M_H(v)|\text{ $E$ is not locally free} \}.
\end{equation}
We shall describe $M_H(v)_s$.
We shall first prove that $M_H(v)_s$ is smooth.
Let $E$ be a point of $M_H(v)_s$ and $x$ the pinch point of $E$.
Then $E\otimes {\cal O}_{X,x} \cong I_x 
\oplus {\cal O}_{X,x}^{\oplus (r-1)}$, where ${\cal O}_{X,x}$ is the stalk
of ${\cal O}_X$ at $x$. 
Since $\Ext^2(E,E) \cong {\Bbb C}$ and 
$E^{\vee \vee}$ is simple,
by using the local-global spectral sequence,
we see that 
\begin{equation}
\Ext^1(E,E) \to H^0(X,{\cal E}xt^1(E,E))
\end{equation}
is surjective.
In the local deformation space
of $E\otimes {\cal O}_{X,x}$,
the locus of non-locally free sheaves is smooth of codimension $r-1$.
Hence $M_H(v)_s$ is smooth of codimension $r-1$.
  
Bu our assumption,
$\langle (v+\omega)^2 \rangle=0$.
Hence $Y:=M_H(v+\omega)$ is an abelian surface and
consists of semi-homogeneous vector bundles.
For simplicity,
we assume that there is a universal bundle
${\cal F}$ on $Y \times X$
(e.g. $v=r+d H+(d^2k-1) \omega$, $(H^2)=2rk$ and $(r, dk)=1$).
We shall prove that ${\Bbb P}:={\Bbb P}({\cal F})$ is isomorphic to $M_H(v)_s$.
Let $\varpi:{\Bbb P} \to Y \times X$ be the projection and
$\varpi^*{\cal F} \to {\cal O}_{\Bbb P}(\lambda)$ the universal quotient
line bundle.
Let $\Gamma$ be the graph of the projection 
${\Bbb P} \to X$.
Then there is a surjective homomorphism
\begin{equation}
 \phi:\varpi^*{\cal F} \boxtimes {\cal O}_X \to {\cal O}_{\Gamma}(\lambda).
\end{equation}
Since ${\cal O}_{\Gamma}(\lambda)$ is flat over ${\Bbb P}$,
$\ker \phi$ is a flat family of torsion free sheaves of
$v(\ker \phi_t)=v$, $t \in {\Bbb P}$.
Thus we get a morphism ${\Bbb P} \to M_H(v)_s$.
It is easy to see that this morphism is proper and bijective.
Hence it is isomorphic.

Let us consider the fiber product
$K_H(v)_s:=M_H(v)_s \times_{M_H(v)} K_H(v)$.
We shall prove that $K_H(v)_s$ is disjoint union of
$r^2$ copies of ${\Bbb P}^{r-1}$.
Let $E_0$ be an element of ${\Bbb P}$ such that
$\varpi(E_0)=(0,0)$.
We may assume that ${\frak a}(E_0)=(0,0)$.
Then the restriction of ${\frak a}$ to ${\Bbb P}$ 
factors ${\Bbb P} \to Y \times X \to \widehat{X} \times X$, where
$Y \times X \to \widehat{X} \times X$ is the map sending
$(y,x)$ to $(\det(y),\alpha(y)-x)$.
By \cite{Mu:1}, $\# \ker(\det)=r^2$
(see Lemma \ref{lem:semi-hom}).
Hence $K_H(v)_s$ is $r^2$ copies of ${\Bbb P}^{r-1}$.

Let $P_1,P_2,\dots,P_{r^2}$ be the $r^2$ copies of ${\Bbb P}^{r-1}$.
Assume that $r \geq 3$.
Let $\widetilde{K_H(v)} \to K_H(v)$ be the blow-up of 
$K_H(v)$ along $P_1$.
Then the exceptional divisor ${\Bbb P}(\Omega_{P_1}^1)$ has 
a natural morphism $\psi:{\Bbb P}(\Omega_{P_1}^1) \to \check{P_1}$,
where $\check{P_1}$ is the dual projective space of $P_1$.
Then we can contract fibers of $\psi$ and get a symplectic manifold
$K_H(v)'$.
That is, $K_H(v)'$ is Mukai's elementary transform of  
$K_H(v)$ along $P_1$ (\cite{Mu:3}, \cite[2.5]{H:2}).

\begin{prop}
$K_H(v)'$ is not K\"{a}hler.
\end{prop}

\begin{proof}
By our assumption on $X$ and Corollary \ref{cor:NS},
$H^{1,1}(K_H(v))_{\Bbb Q}={\Bbb Q}^{\oplus 2}$.
By using $q_{K_H(v)}$,
we have an isomorphism $H^{1,1}(K_H(v))_{\Bbb Q} \cong
H^{2r-3,2r-3}(K_H(v))_{\Bbb Q}$.
Let $l_1,l_2,\dots,l_{r^2}$ be lines on
$P_1,P_2,\dots,P_{r^2}$ respectively.
For a line bundle $A$ of $c_1(A)=c_1(v)$,
we set $M_H(v,A):=\{E \in M_H(v)| \det E=A \}$.
Let $N_H(v,A)$ be the Uhlenbeck's compactification of
the moduli space of $\mu$-stable vector bundles $E$ of
Mukai vector $v$ and $\det E=A$.
By Li \cite{Li:1}, $N_H(v,\det E_0)$ is a projective scheme 
and there is a contraction $f:M_H(v,\det E_0) \to N_H(v,\det E_0)$.
By this contraction,
each $P_i$ are contracted to points on $N_H(v,\det E_0)$.
Hence each $l_i$ are perpendicular to 
$f^*(\Pic(N_H(v,\det E_0)))$.
Let $L$ be the pull-back of an ample line bundle on $N_H(v,\det E_0)$.
Then $(L,C) \ne 0$ for any curve which does not contained in $\cup_i P_i$.
Therefore ${\Bbb Q}\, l_1={\Bbb Q}\,l_2=\dots={\Bbb Q}\,l_{r^2}$.
Let $P_1' \subset K_H(v)'$ be the center of the elementary transformation
and $l_1'$ a line on $P_1'$.
Assume that $K_H(v)'$ is also K\"{a}hler.
Since $H^{1,1}(K_H(v)')_{\Bbb Q} \cong H^{1,1}(K_H(v))_{\Bbb Q}$,
$H^{2r-3,2r-3}(K_H(v)')_{\Bbb Q}$ is also of dimension 2.
Since $L$ is trivial on a neighborhood of $\cup_i P_i$, we get 
$(L,l_1')=(L,l_2)=\dots=(L,l_{r^2})=0$,
where we denote the extension of $L_{|K_H(v)\setminus \cup_i P_i}$
to $K_H(v)'$ by the same $L$.
Therefore we also get that
${\Bbb Q}\,l_1'={\Bbb Q}\, l_2=\dots={\Bbb Q}\, l_{r^2}$.
Since K\"{a}hler form has positive intersection
with effective 1-cycles,
we get that
\begin{equation}\label{eq:l}
 \begin{split} 
  &{\Bbb Q}^+l_1={\Bbb Q}^+l_2=\dots={\Bbb Q}^+l_{r^2},\\
  &{\Bbb Q}^+l_1'={\Bbb Q}^+l_2=\dots={\Bbb Q}^+l_{r^2}.
 \end{split}
\end{equation}
On the other hand, since $K_H(v)'$ is the elementary transform
of $K_H(v)$ along $P_1$,
for an ample line bundle $B$ on $K_H(v)$,
$(B,l_1)>0$ and $(B,l_1')<0$.
By \eqref{eq:l}, this is impossible.
Therefore $K_H(v)'$ is not K\"{a}hler.
\end{proof}
 
\begin{lem}\label{lem:semi-hom}
Keep the notations as above. Then $\# \ker (\det:Y \to \widehat{X}) =r^2$.
\end{lem}

\begin{proof}
Let $E$ be an element of $Y$.
We set 
\begin{equation}
 \begin{split}
  K(E)&:=\{x \in X| T^*_x E \cong E \},\\
  K(\det E)&:=\{x \in X| T^*_x \det E \cong \det E \},\\
  \Sigma(E)&:=\{y \in \widehat{X}|E \otimes {\cal P}_y \cong E \},
 \end{split}
\end{equation}
where ${\cal P}$ is the Poincar\'{e} line bundle on $X \times \widehat{X}$.
Then $Y \cong X/K(E)$ and $\ker \det =K(\det E)/K(E)$.
By \cite[Cor. 7.8]{Mu:1},
\begin{equation}
 \frac{\# K(\det E)}{\# K(E)}=\frac{\# X_r}{\# \Sigma(E)},
\end{equation}
where $X_r$ is the set of $r$-torsion points of $X$. 
Since $\# X_r=r^4$ and $\#\Sigma(E)=r^2$ (\cite[Prop. 7.1]{Mu:1}),
we obtain our lemma.
\end{proof}

\subsection{Intermediate jacobian}

For an irreducible symplectic manifold $M$,
$H^3(M,{\cal O}_{M})=0$.
In this subsection, we shall compute the intermediate jacobian
$J_2(K_H(v)):=H^2(K_H(v),\Omega^1_{K_H(v)})/H^3(K_H(v),{\Bbb Z})$.
We set $H^{odd}(X,{\Bbb Z}):=H^1(X,{\Bbb Z}) \oplus H^3(X,{\Bbb Z})$
 and define weight 3 Hodge structure by
\begin{equation}
 \begin{cases}
  H^{0,3}(H^{odd}(X,{\Bbb C}))=0\\
  H^{1,2}(H^{odd}(X,{\Bbb C}))=H^{0,1}(X) \oplus H^{1,2}(X)\\
  H^{2,1}(H^{odd}(X,{\Bbb C}))=H^{1,0}(X) \oplus H^{2,1}(X)\\
  H^{3,0}(H^{odd}(X,{\Bbb C}))=0.
 \end{cases}
\end{equation}
We define a homomorphism
\begin{equation}
 j_v:H^{odd}(X,{\Bbb Z}) \to 
 H^3(K_H(v),{\Bbb Z})_f
\end{equation}
 by
\begin{equation}
 j_v(x):=-\frac{1}{\rho}
 \left[p_{K_H(v)*}\left(
 \ch({\cal E}_{|K_H(v) \times X})x\right)\right]_{3/2},
\end{equation}
where $H^3(K_H(v),{\Bbb Z})_f$ is the torsion free quotient of 
$H^3(K_H(v),{\Bbb Z})$.
Since $H^1(K_H(v),{\Bbb Z})=0$, it is easy to see that
$j_v$ does not depend on the choice of a quasi-universal family.
\begin{prop}
Let $v$ be a primitive Mukai vector
such that $v>0$, $c_1(v) \in \NS(X)$ and $\langle v^2 \rangle \geq 6$.
Then
$j_v:H^1(X,{\Bbb Z}) \oplus H^3(X,{\Bbb Z}) \to 
 H^3(K_H(v),{\Bbb Z})_f$ is an isomorphism preserving Hodge structures.
In particular, $J_2(K_H(v)) \cong \widehat{X} \times X$.
\end{prop} 
In the same way as in the proof of Theorem \ref{thm:period},
we shall first prove our claim for
the case where $\rk (v)=1$.
Since $\codim_N(N \setminus N_*)=4$,
$H^3(N_*,{\Bbb Z}) \cong H^3(N,{\Bbb Z})$.
Since G\"{o}ttsche (see \cite[p. 50]{Go:1}) proved that $b_3(K_H(v))=8$,
it is sufficient to prove that 
$\im j$ is a direct summand of $H^3(K_H(v),{\Bbb Z})$.  
Since $\codim_{K_H(v)}(K_H(v) \setminus K_H(v)_*)=2$,
$H^3(K_H(v),{\Bbb Z}) \to H^3(K_H(v)_*,{\Bbb Z})$ is injective.
Since $H^3(K_H(v)_*,{\Bbb Z})_f \to 
H^3(B_{\delta}(N_*),{\Bbb Z})^{{\frak S}_n}$
is injective, we shall regard $\im j$ as a submodule of 
$H^3(B_{\delta}(N_*),{\Bbb Z})^{{\frak S}_n}$.
By \eqref{eq:family},
we see that
\begin{equation}
 j_v(x_1+x_3)=\sum_i p_i^*(x_3)_{|B_{\delta}(N_*)}
 -\sum_{i<j}e^{i,j}p_i^*(x_1)_{|B_{\delta}(N_*)}.
\end{equation}
By similar way as in the proof of Lemma \ref{lem:beauville}, 
we see that $\im j_v$ is a direct summand of 
$H^3(B_{\delta}(N_*),{\Bbb Z})^{{\frak S}_n}$.
Thus we get our claim for $\rk v=1$ case.
For general cases, by the following lemma, we can argue as in the proof of 
Theorem \ref{thm:period}.

\begin{lem}
Let ${\cal P}$ be the Poincar\'{e} line bundle on
$X \times \widehat{X}$. Then
\begin{enumerate}
\item[(1)]
${\cal G}_{\cal P}^H$ is defined over ${\Bbb Z}$.
\item[(2)]
Assume that ${\cal G}_{\cal P}$ induces an isomorphism
$K_H(v) \to K_{\widehat{H}}({\cal G}_{\cal P}^H(v))$.
Then the following diagram is commutative.
\begin{equation}
 \begin{CD}
  H^{odd}(X,{\Bbb Z}) @>{(-1)^{i+1}{\cal G}_{\cal P}^H} >>
  H^{odd}(\widehat{X},{\Bbb Z})\\
  @V{j_v}VV @VV{j_{{\cal G}_{\cal P}^H(v)}}V\\
  H^3(K_H(v),{\Bbb Z}) @=H^3(K_{\widehat{H}}({\cal F}_{\cal P}^H(v)),{\Bbb Z}) 
 \end{CD}
\end{equation}
\end{enumerate}
The same assertion also holds for relative Fourier-Mukai functor on
a product of elliptic curves in section 3.2.
\end{lem}
\begin{proof}
We first treat original Fourier-Mukai functor.
By Lemma \ref{lem:homology},
${\cal G}_{\cal P}^H$ is defined over ${\Bbb Z}$.
The second assertion follows from the same computation
in Proposition \ref{prop:comm}.

Let $X=C_1 \times C_2$ be a product of
two elliptic curves $C_1,C_2$.
Let ${\cal P}$ be a universal family on
relative jacobian on $i:X \times_{C_1} X \hookrightarrow X \times X$.  
By a direct computation, we see that 
$\ch (i_*({\cal P})) \in H^*(X \times X,{\Bbb Z})$.
Hence the first claim holds.
The second assertion follows from the same computation
in Proposition \ref{prop:comm}. 
\end{proof}

\subsection{The case where $\langle v^2 \rangle=4$}

In this subsection, we shall treat the remaining case,
that is, $\langle v^2 \rangle =4$.
In this case, $K_H(v)$ is a K3 surface.
We assume that $H$ is general with respect to $v$.
We shall determine this K3 surface.
Let $v=r+\xi+a \omega$,
$\xi \in \NS(X)$ be a Mukai vector of
$\langle v^2 \rangle =4$.
Then we see that $\ell(v)=1,2$.
Replacing $v$ by $v \ch(H^{\otimes m})$, $m \gg 0$,
we may assume that $\xi$ belongs to the ample cone.
Let $\iota:X \to X$ be the $(-1)$-involution of $X$
and $x_1,x_2,\dots,x_{16}$ the fixed points of $\iota$.
Let $\pi:\widetilde{X} \to X$ be the blow-ups of $X$ at 
$x_1,x_2,\dots,x_{16}$ and $E_1,E_2,\dots,E_{16}$
the exceptional divisors of $\pi$. 
Let $q_1:X \to X/\iota$ be the quotient map.
Then the morphism $q_1 \circ \pi:\widetilde{X} \to X/\iota$ factors through
the quotient $\widetilde{X}/\iota$ of $\widetilde{X}$ by $\iota$ 
: $\widetilde{X} \overset{q_2} \to \widetilde{X}/\iota
\overset{\varpi} \to X/\iota$.
Let $\Km(X):=\widetilde{X}/\iota$ be the Kummer surface associated to $X$
and $\varpi:\widetilde{X}/\iota \to X/\iota$ the minimal resolution of 
$X/\iota$.
We set $C_i:=q_2(E_i)$, $i=1,2,\dots,16$.
 
Since the notion of stability only depends on the numerical equivalence 
class of $H$, replacing $H$ by $H \otimes P$, $P \in \Pic^0(X)$,
we may assume that $H$ is symmetric,
that is, $\iota^*H=H$.
Then $H$ has a $\iota$-linearization.
Hence $H^{\otimes 2}$ descend to an ample line bundle $L$ on
$X/\iota$.
Then $L_m:=\varpi^*(L^{\otimes m})(-\sum_{i=1}^{16}C_i)$,
$m \gg0$ is an ample line bundle on $\Km(X)$.
We shall fix a sufficiently large integer $m$.
We would like to relate $K_H(v)$ to $M_{L_m}(w)$
for some $w \in H^*(\Km(X),{\Bbb Z})$.

\subsubsection{The case where $\ell(v)=1$}

Let $w=r+c_1+a \omega\in H^{ev}(\Km(X),{\Bbb Z})$
be an isotropic Mukai vector with $c_1 \in \NS(X)$.
We shall look for some conditions on $w$ such that 
$K_H(v) \cong M_{L_m}(w)$.
By Mukai \cite{Mu:4}, $M_{L_m}(w)$ is not empty,
if $\rk(w)>0$.
%
%
Assume that $\pi_*(q_2^* c_1)=\xi$.

\begin{lem}
$M_{L_m}(w)$ consists of $\mu$-stable 
sheaves, if $\rk w>0$.
\end{lem}
\begin{proof}
%
Let $E$ be an element of $M_{L_m}(w)$.
Then $F:=E^{\vee \vee}$ is a $\mu$-semi-stable vector bundle of
$v(F)=w+k \omega, k \geq 0$.
$q^*_2(F)$ is a $\mu$-semi-stable vector bundle on $X$ with
respect to $q_2^*(L_m)=\pi^*(H^{\otimes 2m})(-2\sum_{i=1}^{16}E_i)$.
Since $m$ is sufficiently large, 
$\pi_*(q_2^*(F))^{\vee \vee}$ is $\mu$-semi-stable with respect to $H$.
Since $\ell(w)=1$ and 
$H$ is a general ample line bundle, 
$\pi_*(q_2^*(F))^{\vee \vee}$ is $\mu$-stable.
By the choice of $m$, $q_2^*(F)$ is also $\mu$-stable.
Hence $F$ is a $\mu$-stable vector bundle,
which implies that $M_{L_m}(w)$ consists of $\mu$-stable sheaves.
\end{proof}
Since $\dim M_{L_m}(w)=\langle w^2 \rangle+2=2$,
every member of $M_{L_m}(w)$ is locally free.
Moreover general members $F$ of $M_{L_m}(w)$ are 
rigid on each $(-2)$-curves $C_i$.
\begin{lem}
We set
$N(w,i):=\{F \in M_{L_m}(w)| \Ext^1_{{\cal O}_C}(F_{|C_i},F_{|C_i}) \ne 0\}$.
Then $N(w,i)$ is not empty if and only if $r|\deg(F_{|C_i})$.
Moreover if $N(w,i)$ is not empty, then $N(w,i)$ is a rational curve.
\end{lem}
\begin{proof}
We assume that $\Ext^1_{{\cal O}_C}(F_{|C_i},F_{|C_i}) \ne 0$.
We set $F_{|C_i}=\oplus_{j=1}^k{\cal O}_{C_i}(a_j)^{\oplus n_j}$,
$a_1<a_2<\dots<a_k$.
Let $F':=\ker (F \to {\cal O}_{C_i}(a_1)^{\oplus n_1})$
be the elementary transformation of $F$ along
${\cal O}_{C_i}(a_1)^{\oplus n_1}$.
Since $v({\cal O}_{C_i}(a_1))=C_i+\chi({\cal O}_{C_i}(a_1))\omega
=C_i+(a_1+1)\omega$, we get
\begin{equation}
v(F')=v(F)-n_1(C_i-(a_1+1)\omega).
\end{equation}
Hence we see that
$\langle v(F')^2\rangle=-2 n_1(\sum_{j \geq 2}n_j(a_j-a_1-1)) \leq 0$.
By the choice of $L_m$,
$F'$ is also $\mu$-stable (cf. \cite[Prop. 2.3]{Y:2}).
Hence $-2 \leq -2 n_1(\sum_{j \geq 2}n_j(a_j-a_1-1))$.
Since $F$ is not rigid (i.e. $\Ext^1(F,F) \ne 0$),
$\sum_{j \geq 2}n_j(a_j-a_1-1) > 0$.
Thus $n_1=\sum_{j\geq 2} n_j(a_j-a_1-1)=1$.
Therefore we get that
\begin{equation}
 F_{|C_i} \cong {\cal O}_{C_i}(a_1)
 \oplus {\cal O}_{C_i}(a_1+1)^{\oplus (r-2)} \oplus  {\cal O}_{C_i}(a_1+2).
\end{equation}
In this case, $\langle v(F')^2 \rangle =-2$, and hence
$F'$ is a unique stable vector bundle
of $v(F')=v(F)-n_1(C_i-(a_1+1)\omega)$.
It is not difficult to see that the choice of inverse transformations
is parametrized by ${\Bbb P}^1$. Therefore $N(w,i)$ is a rational curve.
\end{proof}

We shall consider the pull-back $q^*_2(F)$ of a general member $F$.
Since $F_{|C_i}$, $1 \leq i \leq 16$ are rigid,
replacing $q_2^*(F)$ by
$q_2^*(F)(\sum_{i=1}^{16} s_i E_i)$ for some integers $s_i$,
we may assume that $q^*_2(F)_{|E_i} \cong 
{\cal O}_{E_i}(-1)^{\oplus k_i}  
\oplus {\cal O}_{E_i}^{\oplus (r-k_i)}$.
Let $\phi:q^*_2(F) \to \oplus_{i=1}^{16}{\cal O}_{E_i}(-1)^{\oplus k_i}$
be the quotient map
induced by the quotients $q^*_2(F)_{|E_i} \to {\cal O}_{E_i}(-1)^{\oplus k_i}$.
Then 
\begin{enumerate}
\item
$G:=\ker \phi$ is the elementary transformation of $q_2^*(F)$
along $\oplus_{i=1}^{16}{\cal O}_{E_i}(-1)^{\oplus k_i}$,
\item
$G$ satisfies that $G_{|E_i} \cong {\cal O}_{E_i}^{\oplus r}$.
\end{enumerate}
Hence $\pi_*(G)$ is a stable vector bundle on $X$.
So we get a rational map
$f:M_{L_m}(w) \cdots \to M_H(v)$, where $v=v(\pi_*(G))$.
Since $M_{L_m}(w)$ is a K3 surface, the image of $M_{L_m}(w)$
belongs to a fiber of ${\frak a}$.
Since $q_2^*(F)$ is a stable (and hence a simple vector bundle)
and $\iota$ has fixed points,
$\iota$-linearization on $F$ is uniquely determined by
$q_2^*(F)$.
Hence $f$ is generically injective.  
By a simple calculation, we get that 
\begin{equation}\label{eq:v(G)}
 \begin{split}
  \langle v(G)^2 \rangle &=2r c_2(G)-(r-1)(c_1(G)^2)\\
  &=4r c_2(F)-2(r-1)(c_1(F)^2)-\sum_{i=1}^{16}k_i(r-k_i)\\
  &=2(\langle w^2 \rangle+2r^2)-\sum_{i=1}^{16}k_i(r-k_i)\\
  &=4r^2-\sum_{i=1}^{16}k_i(r-k_i).
 \end{split}
\end{equation}
Hence if $\langle v(G)^2 \rangle=4r^2-\sum_{i=1}^{16}k_i(r-k_i)=4$, 
then the fiber of ${\frak a}$
is isomorphic to $M_{L_m}(w)$.

Conversely for a Mukai vector $v=r+d N+a \omega \in H^{ev}(X,{\Bbb Z})$
such that $(a)$ $N$ is a $(1,n)$-polarization,
$(b)$ $(r,d)=1$ and $(c)$ $\langle v^2 \rangle =d^2 (N^2)-2r a=4$,
we shall look for such a vector $w \in H^{ev}(\Km(X),{\Bbb Z})$.
We shall divide the problem into two cases.
We also treat the case where $r=0$.

Case (I). We first assume that 
$r$ is even.
In this case, $d$ must be odd.
By the condition $(c)$, $(N^2)=2n$ is divisible by 4.
Thus $n$ is an even integer.
In this case, replacing $N$ by $N \otimes P$,
$P \in \Pic^0(X)$,
we may assume that $N$ has an $\iota$-linearization which acts 
trivially on the fibers of $N$ at exactly 4 points (cf. [L-B, Rem. 7.7]).
Replacing the indices, we assume that the 4 points are 
$x_1,x_2,x_3,x_4$.
We set 
\begin{equation}
 \begin{cases}
  N_1:=\pi^*(N^{\otimes d})(\frac{r-2}{2}\sum_{i=1}^4 E_i+
  \frac{r}{2}\sum_{i \geq 5}E_i),\\
  N_2:=N_1(-r E_1).
 \end{cases}
\end{equation}
Then for suitable linearizations,
$N_1$ and $N_2$ descend to line bundles 
$\xi_1$ and $\xi_2$ on $\Km(X)$ respectively.
By simple calculations, we get that
\begin{equation}
 \begin{split}
  (\xi_1^2) &=d^2 \frac{(N^2)}{2}-2r^2+2r-2\\
  &=r(a-2r+2),\\
  (\xi_2^2) &= r(a-2r+1).
 \end{split}
\end{equation}
(I-1) We first assume that $r>0$ and set
\begin{equation}
 w:=
  \begin{cases}    
   r+\xi_1+\frac{a-2r+2}{2}\omega, \text{ if $a$ is even},\\
   r+\xi_2+\frac{a-2r+1}{2}\omega, \text{ if $a$ is odd}.
  \end{cases}
\end{equation}
Then we get that $\langle w^2 \rangle =0$.
Let $F$ be a general stable vector bundle of $v(F)=w$.
By the choice of $\xi_1$ and $\xi_2$, 
the restriction of $q_2^*(F)$ or $q_2^*(F)(E_1)$
to $E_i$ is isomorphic to
${\cal O}_{E_i}(-1)^{\oplus k_i} \oplus {\cal O}_{E_i}^{\oplus (r-k_i)}$,
where $k_i=(r-2)/2$ for $1 \leq i \leq 4$ and
$k_i=r/2$ for $i \geq 5$.
Then by \eqref{eq:v(G)},
we get that $\langle v(\pi_*(G))^2 \rangle=4$.
Since $\rk(\pi_*(G))=r$ and $c_1(\pi_*(G))=dN$,
$v(\pi_*(G))$ must be equal to $v$.
Therefore $K_H(v)$ is isomorphic to $M_{L_m}(w)$.
\newline
(I-2) We next assume that $r=0$.
In this case, $d=1$ and $(N^2)=4$.
Then $\dim H^0(X,N)=2$. 
We set
$w:=\xi_1+b \omega$, $b \ne 0$.
Since $(\xi_1^2)=0$ and $(\xi_1,L_m)>0$,
we get $\dim H^0(\Km(X),\xi_1) \geq 2$.
Since $\dim  H^0(X,N) \geq H^0(\Km(X),\xi_1)$,
$\dim H^0(\Km(X),\xi_1) =2$.
If $|\xi_1|$ does not have a fixed component,
it defines an elliptic fibration
$\Km(X) \to {\Bbb P}^1$.
In this case, we see that $M_{L_m}(w)$ is not empty and
is isomorphic to a compactification of the relative jacobian. 
For an elliptic curve $D \in |\xi_1|$,
$q_2^{-1}(D) \to D$ is a double cover branched at
$q_2(C_i) \cap D$, $1 \leq i \leq 4$. 
Then we have a birational map
$M_{L_m}(w) \cdots \to K_H(v)$ by sending a line bundle 
$F$ ($\in M_{L_m}(w)$) on $D$  
to $\pi_*(q_2^*(F)) \in K_H(v)$. 
Therefore $K_H(v)$ is isomorphic to $M_{L_m}(w)$. 
If $|\xi_1|$ have a fixed component, then
by a classification of $N$,
\begin{enumerate}
\item
$X$ is a product of two elliptic curves $E_1$, $E_2$  
and
\item
 $N={\cal O}_{E_1}(P_1) \boxtimes{\cal O}_{E_2}(P_2+(-1)^*(P_2))$,
where $P_1 \in E_1$ is a 2-torson point and $P_2 \in E_2$. 
\end{enumerate}
In this case, the second projection
$p_2:X \to E_2$ induces an elliptic fibration
$e:\Km(X) \to E_2/\pm1={\Bbb P}^1$.
Then the section $\{P_1\} \times E_2 \subset X$ of $p_2$ induces a
section $\sigma$ of $e$.
Let $f$ be a fiber of $e$.
Then we see that $\xi_1=\sigma+f$.
Let $y$ be a point of $E_2/\pm1$ such that $e$ and $E_2 \to E_2/\pm1$
are smooth over $y$.
Assume that $M_{L_m}(w) \ne \emptyset$.
Let $F$ be an element of $M_{L_m}(w)$ such that $\Supp(F)=\sigma \cup
e^{-1}(y)$.
Then it is easy to see that $\pi_*(q^*_2(F))$ belongs to
$K_H(v)$.
Hence we have a birational map $M_{L_m}(w) \cdots \to K_H(v)$.
Therefore $K_H(v)$ is isomorphic to $M_{L_m}(w)$.  
For the non-emptyness of $M_{L_m}(w)$, we quote the following easy lemma.

\begin{lem}
Let $\lambda$ be a rational number such that
$b/\lambda \not \in {\Bbb Z}$.
Let $n$ be an integer of $b/\lambda<n <b/\lambda+1$.
Let $F$ be a coherent sheaf of $v(F)=w$.
Then $F$ is stable with respect to 
$\sigma+(\lambda+1)f+L'$, $L' \in \langle \sigma,
f \rangle^{\perp}$, if and only if $F$ fits in a non-trivial extension
\begin{equation}
 0 \to {\cal O}_{\sigma}(b-n-1) \to F \to I \to 0,
\end{equation}
where $I$ is a line bundle of degree $n$ on a fiber of $e$.
In particular, $M_{\sigma+(\lambda+1)f+L'}(w)$ is isomorphic 
to a compactification of relative jacobian.
\end{lem}

\vspace{1pc}

Case (II). We assume that $r$ is odd.
Replacing $v$ by $v \ch(N)$, we may assume that $d$ is even.
We set $N_1:=\pi^*(N^{\otimes d})(\frac{r-1}{2}\sum_{i=1}^{16}E_i)$.
Then for a suitable linearization,
$N_1$ descends to a line bundle $\xi$ on $\Km(X)$.
By a simple calculation, we get that 
$(\xi^2)=r(a-2r+4)$.
Since $d$ is even and $r$ is odd, condition $(c)$ implies that
$a$ is an even integer.
We set $w:=r+\xi+\{(a-2r+4)/2\}\omega$.
Then we get that $\langle w^2 \rangle=0$.
In the same way as above, we see that $v(\pi_*(G))=v$,
which implies that $K_H(v) \cong M_{L_m}(w)$.

\begin{rem}
By the choice of $k_i$,
if $r>2$, then $N(w,i)$ is empty.
Thus $f$ is a morphism.
If $r=2$, then $N(w,i)$, $1 \leq i \leq 4$
is not empty and these closed subset correspond to
the closed subset 
\begin{equation}
 N(v,i):=\{G \in K_H(v)|
 \text{$G$ is not locally free at $x_i$}\}.
\end{equation}
\end{rem}   

\subsubsection{The case where $\ell(v)=2$}
In this case, we may assume that $v=2(r+dN)+a \omega$,
where $N$ is a primitive ample line bundle and $d^2(N^2)=ra+1$.
We use the same notations and methods as in the case where $\ell(v)=1$.
For this purpose, we prepare a lemma.
\begin{lem}\label{lem:v4}
Keep the notations as above.
Assume that $r \geq 2$.
If $H$ is general, then every $\mu$-semi-stable sheaf $E$ of
$v(E)=v$ is a $\mu$-stable vector bundle.
\end{lem}

\begin{proof}
Let $E$ be a $\mu$-semi-stable sheaf of $v(E)=v$.
Assume that $E$ is not $\mu$-stable.
Since $H$ is general,
there is an exact sequence
\begin{equation}
 0 \to E_1 \to E \to E_2 \to 0,
\end{equation}
where $E_i$, $i=1,2$ are $\mu$-stable sheaves of 
$v(E_i)=r+dN+a_i \omega$.
Since $d^2(N^2)=ra+1$ and $a_1+a_2=a$,
we see that $\langle v(E_1)^2 \rangle=
r(a_2-a_1)+1$ and $\langle v(E_2)^2 \rangle=r(a_1-a_2)+1$.
Since $\langle v(E_1)^2 \rangle, \langle v(E_2)^2 \rangle \geq 0$,
$\langle v(E_1)^2 \rangle=0$ or $\langle v(E_2)^2 \rangle=0$.
Then we get $r(a_1-a_2)=\pm1$.
Since $r>1$, this is impossible.
Therefore $E$ is $\mu$-stable.
If $E$ is not locally free, then 
$F:=E^{\vee \vee}$ is a $\mu$-stable locally free sheaf of 
$v(F)=v+b \omega$, $b>0$.
Since $r>1$, we see that
$\langle v(F)^2 \rangle=\langle v^2 \rangle-4rb<0$.
Hence $E$ must be locally free. 
\end{proof}

We first assume that $r>1$.
We set $N_1:=\pi^*(N^{\otimes 2d})(-r\sum_{i=1}^{16}E_i+2E_1)$.
Then for a suitable linearization on $N_1$,
$N_1$ descend to a line bundle $\xi_1$ on $\Km(X)$.
We set $w:=2r+\xi_1+\frac{a+1-4r}{2}\omega$.
Then we see that $\langle w^2 \rangle=0$.
Let $F$ be a $\mu$-semi-stable vector bundle of $v(F)=w$
with respect to $L_m$.
By the choice of $\xi_1$ and $\lambda$,
$\pi_*(q_2^*(F))^{\vee \vee}$ is a $\mu$-semi-stable vector bundle of
$v(\pi_*(q_2^*(F))^{\vee \vee})=v+b \omega$, $b \geq 0$.
By Lemma \ref{lem:v4} and 
the choice of $L_m$,
$q_2^*(F)$ is a $\mu$-stable vector bundle
with respect to $q_2^*(L_m)$ and $\pi_*(q_2^*(F))^{\vee \vee}$
 is a $\mu$-stable vector bundle of $v(\pi_*(q_2^*(F))^{\vee \vee})=v$.
Hence $K_H(v)$ is isomorphic to $M_{L_m}(w)$.

If $r=1$, then $M_H(v)$ is isomorphic to
$M_H(v')$, $v'=2-\omega$.
In this case, Mukai \cite[Cor. 4.5]{Mu:5} (see Theorem \ref{thm:00/08/17})
showed that
Fourier-Mukai transform ${\cal G}_{\cal P}$
by Poincar\'{e} line bundle
${\cal P}$ on $\widehat{X} \times X$ gives an isomorphism
$M_H(v') \to \Hilb_{\widehat{X}}^2 \times X$.
Therefore $K_H(v)$ is isomorphic to $\Km(\widehat{X})$.    
By the isomorphism induced by ${\cal G}_{\cal P}$,
a generel member of $M_H(v')$ fits in an exact sequence
\begin{equation}\label{eq:v=(2,0,-1)}
 0 \to E \to {\cal P}_x \oplus {\cal P}_y \to {\Bbb C}_z \to 0
\end{equation}
where $x,y \in \widehat{X}$ and $z \in X$.
By using \eqref{eq:v=(2,0,-1)},
we give another description of $K_H(v)$.
Assume that $y=-x$ and $z=0$.
Then $E$ belongs to $K_H(v)$.
Since $\iota^*({\cal P}_x)\cong {\cal P}_y$,
$G:=\pi^*({\cal P}_x \oplus {\cal P}_y)$ has a $\iota$-linearization.
We note that the action of $\iota$ on $G_{|E_i}=
\pi^*({\cal P}_x)_{|E_i} \oplus \pi^*({\cal P}_y)_{|E_i}$
is given by $(a,b) \mapsto (b,a)$.
Let $\psi_i:G_{E_i} \to {\cal O}_{E_i}$ be the homomorphism
sending $(a,b)$ to $a-b$.
Let $\psi:G \to \oplus_{i=1}^{16}{\cal O}_{E_i}$ be the 
composition of the restriction $G \to \oplus_{i=1}^{16}G_{|E_i}$
with $\oplus \psi_i:\oplus_{i=1}^{16}G_{|E_i} \to 
\oplus_{i=1}^{16}{\cal O}_{E_i}$.
Then $\ker(\psi)$
is $\iota$-linearlized locally free sheaf such that
the action is trivial on $\cup_i E_i$.
Hence there is a locally free sheaf $F$ on $\Km(X)$ such that
$q_2^*(F)=\ker (\psi)$. 
By our construction of $\ker(\psi)$,
$\ker(\psi)$ is stable with respect to $q_2^*(L_m)$.
Hence $F$ is stable with respect to $L_m$.
Moreover we see that $\ch(q_2^*(F))=2-\sum_i E_i-8\omega$,
and hence $v(F)=2-(\sum_i C_i)/2-2\omega=:w$.
Since $w$ is isotropic, $M_{L_m}(w)$ is a K3 surface.
Therefore $K_H(v)$ is isomorphic to $M_{L_m}(w)$.

Combining all together, we obtain the following theorem.
\begin{thm}
Let $v\in H^{ev}(X,{\Bbb Z})$ be a Mukai vector
such that $v>0$, $c_1(v) \in \NS(X)$ and $\langle v^2 \rangle =4$.
Assume that $H$ is general with respect to $v$.
Let $\Km(X)$ be the Kummer surface associated to $X$.
Then there is an isotropic Mukai vector 
$w \in H^{ev}(\Km(X),{\Bbb Z})$ and an ample line bundle $H'$
on $\Km(X)$ such that $K_H(v)$ is isomorphic to $M_{H'}(w)$.
\end{thm}


\section{Appendix: Stable sheaves on a family of abelian (or K3) surfaces}

In this appendix, we shall prove the following proposition.
\begin{prop}\label{prop:modification}
Let $T$ be a connected smooth curve and
$({\cal X}, {\cal L})$ a pair of a smooth family of
abelian surfaces (resp. K3 surfaces)
$p_T:{\cal X} \to T$ and a relatively ample line bundle
${\cal L}$.
Assume that ${\cal X}_{t_1}$ 
is an abelian surface (resp. K3 surface)
of $\NS({\cal X}_{t_1})= {\Bbb Z} {\cal L}_{t_1}$, 
for a point $t_1 \in T$
(in particular ${\cal L}$ is primitive).
Let $v=r+d{\cal L}+a \omega \in R^*p_{T*}{\Bbb Z}$
be a prmitive Mukai vector.
Then, 
\begin{enumerate}
\item[(1)]
there is an algebraic space
${\frak M}(v) \to T$ which is smooth, proper
and ${\frak M}(v)_t=M_{H_t}(v_t), t \in T$
for a general ample divisor $H_t$ on ${\cal X}_t$.
Moreover if we choose a finite subset $T_0 \subset T$ and
choose any ample divisor $H_t'$ on ${\cal X}_t$, $t \in T_0$,
such that
$\overline{M}_{H_t'}(v_t)=M_{H_t'}(v_t)$,
then we can construct ${\frak M}(v)$ so that
${\frak M}(v)_t=M_{H_t'}(v_t)$ for $t \in T_0$.
\item[(2)]
For the family ${\frak M}(v) \to T$, the homomorphism
$\theta_{v_t}:v_t^{\perp} \to H^2({\frak M}(v)_t,{\Bbb Z}), t \in T$ 
gives a homomorphism between local systems
$\{v_t^{\perp} \}_{t \in T} \to \{H^2({\frak M}(v)_t,{\Bbb Z})\}_{ t \in T}$. 
\item[(3)]
If $p_T:{\cal X} \to T$ is a family of abelian surfaces
with a section,
then we have a family of morphisms ${\frak a}_{v_t}:
{\frak M}(v)_t \to ({\cal X} \times_T \Pic^0_{{\cal X}/T})_t$,
$t \in T$.
\end{enumerate}
\end{prop}
\begin{proof}
The proof is similar to \cite[Prop. 3.3]{Y:3}.
But we repeat the proof
since our assumption on $v$ is weaker than that in \cite[Prop. 3.3]{Y:3}.

 Let $g:\Pic_{{\cal X}/T} \to T$ be the relative Picard scheme.
We denote the connected component of $\Pic_{{\cal X}/T}$ containing
the section of $g$ which corresponds to the family $d{\cal L}$ by 
$\Pic^{\xi}_{{\cal X}/T}$.
Since $\Pic^0_{{\cal X}/T} \cong \Pic^{\xi}_{{\cal X}/T}$,
$\Pic^{\xi}_{{\cal X}/T} \to T$ is a smooth morphism.
Let $h:\overline{{\frak M}}_{{\cal X}/T}(v) \to T$
be the moduli scheme parametrizing $S$-equivalence classes
of ${\cal L}_t$-semi-stable sheaves $E$ on ${\cal X}_t$ with
$v(E)=v_t$ \cite{Ma:1}.
Let $D$ be the closed subset of 
$\overline{{\frak M}}_{{\cal X}/T}(v)$
consisting of properly ${\cal L}_t$-semi-stable sheaves on ${\cal X}_t$.
Since $h$ is a proper morphism, $h(D)$ is a closed subset of $T$.
Since $h(D)$ does not contain $t_1$ and $T$ is an irreducible curve,
$h(D)$ is a finite point set.
Since our problem is local, we may assume that
$h(D)=\{t_0 \}$.
Let $H$ be a general ample divisor on ${\cal X}_{t_0}$.
Let $s:\Spl_{{\cal X}/T}(v) \to T$ be the moduli space of 
simple sheaves $E$ on ${\cal X}_t, t \in T$ with
$v(E)=v_t$ \cite[Thm. 7.4]{A-K:1}.
Let $U_1$ be the closed subset of $s^{-1}(T \setminus \{t_0 \})$
consisting of simple sheaves on ${\cal X}_t$, $t \in T \setminus \{t_0 \}$ 
which are not stable with respect to ${\cal L}_t$ 
and $\overline{U_1}$
the closure of $U_1$ in $\Spl_{{\cal X}/T}(v)$.
Let $U_2$ be the closed subset of $s^{-1}(t_0)$
consisting of simple sheaves which are not semi-stable with respect to $H$.
In Lemma \ref{lem:valuative}, we shall prove that 
$\overline{U_1} \cap s^{-1}(t_0)$ is a subset of
$U_2$.
We set ${\frak M}(v):=\Spl_{{\cal X}/T}(v) 
\setminus (\overline{U_1} \cup U_2)$.
Then ${\frak M}(v)$ is an open subspace of $\Spl_{{\cal X}/T}(v)$
which is of finite type and contains all 
$H$-stable sheaves on ${\cal X}_{t_0}$. 
By using valuative criterion of separatedness and properness,
we get that $s:{\frak M}(v) \to T$ is a proper morphism.
Indeed, since ${\frak M}(v) \times_T (T \setminus \{t_0\}) 
\to T \setminus \{t_0\}$
is proper, it is sufficient to check these properties near the fibre 
${\cal X}_{t_0}$.   
The separatedness follows from base change theorem and stability
with respect to $H$ (cf. \cite[Lem. 7.8]{A-K:1}),
and the properness follows from Lemma \ref{lem:valuative}
below
and the projectivity of ${\frak M}(v)_{t_0}$.
Since $\Pic^{\xi}_{{\cal X}/T} \to T$ is a smooth morphism,
\cite[Thm. 1.17]{Mu:3} implies that $s:{\frak M}(v) \to T$
is a smooth morphism.
Thus (1) holds.

We next prove the second claim.
By the proof of \cite[Thm. A.5]{Mu:4}, there is a quasi-universal family
${\cal E}$ on ${\frak M}(v) \times_T {\cal X}$.
Since ${\cal L}$ is relatively ample, there is a locally free resolution
of ${\cal E}$.
Hence we can define Chern character of ${\cal E}$.
Then $\theta_{v_t}(x)=-(1/\rho)[p_{{\frak M}(v)_t*}((\ch {\cal E})_t
(\sqrt{\td_{{\cal X}/T}})_t x^{\vee})]_1$, $x \in v_t^{\perp}$,
where $\td_{{\cal X}/T}$ is the relative Todd class of ${\cal X} \to T$
and $\rho$ is the similitude of ${\cal E}$.
Therefore $\theta_{v_t}$ is a homomorphism between local systems.

Assume that $p_T$ is a family of abelian surfaces with a section
$\sigma$.
Then $R^4p_{T*}{\Bbb Z} \cong {\Bbb Z}\sigma$ 
and we have a flat family of
coherent sheaves $V_1$ and $V_2$ on ${\cal X}$ such that 
$v((V_1)_t-(V_2)_t)=v_t$, $t \in T$.
Moreover there is a universal line bundle on 
${\cal X} \times_T \Pic^0_{{\cal X}/T}$. 
Hence by using the formal difference $V_1-V_2$, we can define a morphism
${\frak a}:{\frak M}(v) \to
{\cal X} \times_T \Pic^0_{{\cal X}/T}$.
\end{proof}

\begin{lem}\label{lem:line}
Keep the notation above.
Let $R$ be a discrete valuation ring over ${\cal O}_T$.
Assume that ${\cal O}_T \to R$ is injective.
Let $L$ be a line bundle over ${\cal X} \otimes_{{\cal O}_T}R$.
Then $c_1(L)=l c_1({\cal L}\otimes_{{\cal O}_T}R)$ for some $l \in {\Bbb Z}$.
\end{lem}

\begin{proof}
We set ${\cal L}_R:={\cal L}\otimes_{{\cal O}_T}R$.
Let $P(x)$ be the Hilbert polynomial of $L$ with respect to 
${\cal L}_{R}$.
Let $E$ be a locally free sheaf on ${\cal X}$ such that 
there is a surjective homomorphism $E \otimes_{{\cal O}_T}R \to L$,
and we shall consider the quot scheme
$Q:=\Quot_{E/{\cal X}/T}^{P(x)}$.
Then $L$ defines a morphism $\tau:\Spec(R) \to Q$
such that $L=(\tau \times_T 1_{{\cal X}})^* {\cal Q}$,
where ${\cal Q}$ is the universal quotient.  
Let $Q_0$ be the connected component of $Q$
which contains the image of $\Spec(R)$.
Since $\Spec(R) \to T$ is dominant,
${\frak q}:Q_0 \to T$ is dominant, and hence surjective.
Since $\NS({\cal X}_{t_1})={\Bbb Z}{\cal L}_{t_1}$, we get that
$c_1({\cal Q}_{q_1})=l c_1({\cal L}_{q_1})$, $l \in {\Bbb Z}$,
where $q_1 \in {\frak q}^{-1}(t_1)$. 
Hence $c_1({\cal Q})=l c_1({\cal L} \otimes_{{\cal O}_T}{\cal O}_{Q_0})$ 
as an element of 
$R^2p_{Q_0*}{\Bbb Z}$, where $p_{Q_0}:{\cal X} \times_T Q_0 \to Q_0$ is the 
projection.
Therefore we get our lemma.
\end{proof}

\begin{lem}\label{lem:valuative}
Keep the notation as above.
Let $R$ be a discrete valuation ring,
$K$ the quotient field of $R$, and $k$ the residue field of $R$.
Let $\Spec(R) \to T$ be a dominant morphism such that 
$\Spec(k) \to T$ defines the point $t_0$.
\begin{enumerate}
\item[(1)]
Let $E$ be a $R$-flat coherent sheaf on $X_R$ such that $E_K:=E \otimes_R K$
is not semi-stable with respect to ${\cal L}_K$.
Then $E_{t_0}$ is not semi-stable with respect to any ample divisor on
${\cal X}_{t_0}$. 
\item[(2)]
For a ${\cal L}_K$-stable sheaf $E_K$ on $X_K$,
there is a $R$-flat coherent sheaf $E$ on $X_R$ such that $E \otimes_R K=E_K$
and $E \otimes_R k$ is a $H$-stable sheaf.
\end{enumerate}
\end{lem}

\begin{proof}
(1)
We first assume that $\rk E>0$.
Since $E_K$ is not semi-stable,
there is a quotient sheaf $E_K \to F_K$ such that
\begin{enumerate}
\item
\begin{equation}\label{eq:not-ss1}
  \frac{(c_1(E_K),c_1({\cal L}_K))}{\rk E_K}>
  \frac{(c_1(F_K),c_1({\cal L}_K))}{\rk F_K},
\end{equation}
or
\item
\begin{equation}\label{eq:not-ss2}
 \frac{(c_1(E_K),c_1({\cal L}_K))}{\rk E_K}=
 \frac{(c_1(F_K),c_1({\cal L}_K))}{\rk F_K},\; 
 \frac{\chi(E_K)}{\rk E_K}>
 \frac{\chi(F_K)}{\rk F_K}.
\end{equation}
\end{enumerate}
Let $E \to F$ be a flat extension of $E_K \to F_K$.
By Lemma \ref{lem:line},
\eqref{eq:not-ss1} and \eqref{eq:not-ss2} implies that
$c_1(E)/\rk E-c_1(F)/\rk(F)=l c_1({\cal L}_R), l \in {\Bbb Q}_{\geq 0}$.
Since $(c_1({\cal L}_{t_0}),H) >0$ for any ample divisor $H$ 
on ${\cal X}_{t_0}$,
similar relations to \eqref{eq:not-ss1}, \eqref{eq:not-ss2}
hold for $F \otimes_R k$.
Thus $E \otimes_R k$ is not semi-stable with respect to 
any ample divisor on ${\cal X}_{t_0}$. 

If $\rk E=0$, then by using the inequality
\begin{equation}
 \frac{\chi(E_K)}{(c_1(E_K),c_1({\cal L}_K))}>
 \frac{\chi(F_K)}{(c_1(F_K),c_1({\cal L}_K))},
\end{equation}
we can prove our claim.
 
(2)
The proof is very similar to \cite[Lem. 3.4]{Y:3}.
In \cite{Y:3}, we only consider the case where $\mu$-semi-stable sheaves 
are $\mu$-stable.
By using a similar method as in the proof of (1),
we can easily modify the arguments in the second paragraph of  
the proof of \cite[Lem. 3.4]{Y:3}.
\end{proof}

\section{Appendix: Evidence of Conjecture \ref{conj:1}}

Throughout of this section, we assume that $X$ is an
abelian surface with $\NS(X)={\Bbb Z}H$,
where $H$ is the ample generator.
We consider Fourier-Mukai transform ${\cal F}_{\cal P}$
(resp. ${\cal G}_{\cal P}$) induced by the Poincar\'{e} line bundle
${\cal P}$ on $X \times \widehat{X}$.
In particular, we shall prove Theorem \ref{thm:evidence}. 

\subsection{Proof of Theorem  \ref{thm:evidence} (i)}

\begin{lem}\label{lem:BN1-2}
Let $v=r+dH+a \omega$ be a Mukai vector such that
$r>0$, $d \geq 0$, $(r,d)=1$ and $a \leq 0$.
Then $H^0(X,E)=0$ for a general $E \in M_H(v)$. 
\end{lem}

\begin{proof}
We shall prove our claim by induction on $r$.

(I) Assume that $r=1$.
Then $M_H(v)\cong \widehat{X} \times \Hilb_X^{\langle v^2 \rangle/2}$.
For $I_Z \in \Hilb_X^{\langle v^2 \rangle/2}$ and $L \in \Pic^0(X)$,
we get $\chi(I_Z(dH) \otimes L)=a \leq 0$.
Hence 
$H^0(X,I_Z(dH) \otimes L)=0$ for general $I_Z$ and $L$.

(II) Let $(r_1,d_1)$ be a pair of integers such that
$d_1r-d r_1=1$ and $0<r_1<r$.
We set $(r_2,d_2):=(r-r_1,d-d_1)$.
We may assume that $d_1 >0$ and $d-d_1 \geq 0$.
We shall choose Mukai vectors 
$v_i:=r_i+d_i H+a_i \omega$, $i=1,2$ such that
$a_i \leq 0$, $i=1,2$.
We shall choose $E_i \in M_H(v_i)$, $i=1,2$ such that
$H^0(X,E_i)=0$.
If $r_2>1$, then we may assume that $E_2$ is locally free.
In this case, we have $d_2>0$. Hence
$-\chi(E_1,E_2)=
d_1d_2(H^2)-r_2a_1-r_1a_2>0$.
Then there is a non-trivial extension
\begin{equation}
0 \to E_2 \to E \to E_1 \to 0.
\end{equation}
By \cite[Lem. 2.1]{Y:5},
$E$ is a $\mu$-stable sheaf.
Therefore our claim holds.
If $r_2=1$, then 
$-\chi({\cal O}_X(d_2H),E_1)=d_1d_2(H^2)-a_1-r_1d_2^2(H^2)/2
=d_2(d_1r_2-r_1d_2/2)(H^2)-a_1>0$, unless 
$d_2=a_1=0$.
Assume that $a_1 >0$ or $d_2>0$.
Then we have a non-trivial extension
\begin{equation}
0 \to {\cal O}_X(d_2H) \to E' \to E_1 \to 0.
\end{equation}
 \cite[Lem. 2.1]{Y:5} implies that $E'$ is a $\mu$-stable sheaf of 
$v(E')=v+b \omega$, 
where $b=\chi({\cal O}_X(d_2H))-\chi(E_2)\geq \chi({\cal O}_X(d_2H))$.
We choose points $x_1,x_2,\dots,x_b \in X$ such that
$H^0(X,I_Z(d_2H))=0$, where $Z=\{x_1,x_2,\dots,x_b \}$.
If we take a suitable surjection
$\phi:E' \to \oplus_{i=1}^b {\Bbb C}_{x_i}$,
then $E:=\ker \phi$ fits in an exact sequence
\begin{equation}
0 \to I_Z(d_2H) \to E \to E_1 \to 0.
\end{equation}
 Since $v(E)=v$, $E$ is a desired element of $M_H(v)$.
Next we assume that $a_1=d_2=0$.
Then we get $d=d_1=1$.
By Proposition \ref{lem:fourier3}, the claim holds.
\end{proof}

By using Lemma \ref{lem:BN1-2}, we get the following proposition.
\begin{prop}\label{prop:BN1-2}
Let $v=l(r+dH)+a \omega$ be a primitive Mukai vector of
$r,l>0$, $d \geq 0$, $(r,d)=1$ and $a \leq 0$.
Then $H^0(X,E)=0$ for a general $E \in M_H(v)$. 
\end{prop}
\begin{proof}
We choose integers $a_1,a_2,\dots,a_l$ such that
$\sum_{i=1}^l a_i=a$ and $a_i \leq 0$ for $1 \leq i \leq l$.
We set $v_i:=r+dH+a_i \omega$.
By Lemma \ref{lem:BN1-2}, 
we can choose elements $E_i \in M_H(v_i)$, $1 \leq i \leq l$
such that $H^0(X,E_i)=0$.
We set $E:=\oplus_{i=1}^l E_i$.
Then $E$ is $\mu$-semi-stable and $H^0(X,E)=0$.
Since $\langle v^2 \rangle \geq 2l^2$,
the proof of \cite[Lem. 4.4]{Y:5} implies that our proposition holds.
\end{proof}
The following is the same as Theorem \ref{thm:evidence} (i).
\begin{thm}
 Let $v=r+dH+a \omega$ be a primitive Mukai vector such that $d>0$ and
 $a \leq 0$.
 Then ${\cal F}_{\cal P}$ induces a birational map
 $M_H(v) \cdots \to M_{\widehat{H}}(-{\cal F}_{\cal P}^H(v))$.
\end{thm}
\begin{proof}
By Proposition \ref{prop:BN1-2},
${\cal F}_{\cal P}^0(E)=0$
for a general $E \in M_H(v)$.
Since ${\cal F}_{\cal P}^2(E)=0$,
$\WIT_1$ holds for a general $E$.
Then by Corollary \ref{cor:wit->birat}, we get our claim.
\end{proof}

\subsection{Proof of Theorem \ref{thm:evidence} (ii)}

If $d=0$, then by Theorem \ref{thm:00/08/17}, 
we get Theorem \ref{thm:evidence} for this case. 
Next we shall consider the case where $d,a>0$.
We devide the proof into several cases.

We first treat the case where $d \equiv 1 \mod r$.

\begin{lem}\label{lem:rk=1}
(1) Assume that $\chi(I_Z(H))>0$.
Then for $n \geq 1$,
$H^1(X,I_Z(nH) \otimes {\cal P}_x) =0$
except finitely many points $x \in \widehat{X}$.
(2) If $\chi(I_Z(H))=0$ and $n \geq 1$, then 
$H^1(X,I_Z(nH) \otimes {\cal P}_x) =0$
for a general $x \in \widehat{X}$.
\end{lem}

\begin{proof}
We first prove (1) by induction on $n$.
By Proposition \ref{lem:fourier4}, the assertion holds for $n=1$.
For $n>1$,
we choose a curve $C$ of $C \equiv H$ (numerical equiv.).
Since $H^0(X,{\cal O}_X(H) \otimes {\cal P}_y) \to
H^0(X,{\cal O}_Z(H) \otimes {\cal P}_y)$ is surjective for 
some $y \in \widehat{X}$ because of the claim for $n=1$,
we may assume that $C$ does not meet $Z$.
Then we have an exact sequence
\begin{equation}
 0 \to I_Z((n-1)H) \to I_Z(nH) \otimes {\cal P}_y
 \to {\cal O}_C(nH) \otimes {\cal P}_y \to 0.
\end{equation}
Since $\deg({\cal O}_C(nH))=n(2g(C)-2)> 2g(C)-2$,
$H^1(X,{\cal O}_C(nH) \otimes {\cal P}_x)=0$ for all $x \in \widehat{X}$.
Hence if $H^1(X,I_Z((n-1)H) \otimes {\cal P}_x)=0$ except
finitely many points $x \in \widehat{X}$,
then the claim also holds for $I_Z(nH)$.

For the claim (2), we use Remark \ref{rem:r=0} 
instead of Proposition \ref{lem:fourier4}.
 \end{proof}

\begin{lem}\label{lem:special} 
Assume that $v=r+H+a \omega$ satisfies
$r>1$ and $2r>(H^2)-2ra \geq 0$.
Then there is a locally free sheaf $E \in M_H(v)$ such that
$E$ fits in an exact sequence
\begin{equation}\label{eq:special}
0 \to \oplus_{i=1}^{r-1}{\cal P}_{x_i} \to
E \to I_Z(H) \to 0,
\end{equation}
where $x_i \in \widehat{X}$.
\end{lem}

\begin{proof}
By our assumption,
$M_{\widehat{H}}(a+\widehat{H}+r \widehat{\omega}) \ne \emptyset$.
We choose an element $F \in M_{\widehat{H}}(a+\widehat{H}+r \widehat{\omega})$.
By Proposition \ref{lem:fourier4}, $\WIT_2$ holds for $F$ with respect to 
${\cal G}_{\cal P}$ and ${\cal G}_{\cal P}^2(F)$ belongs to $M_H(v)$.
Since $r>1$, we may assume that $E:={\cal G}_{\cal P}^2(F)$ is locally free.
We show that $E$ satisfies our claim.
Let $x_1,x_2,\dots,x_{r-1}$ be distinct points of $\widehat{X}$.
Let $G$ be an element of $M_{\widehat{H}}(a+\widehat{H}+\widehat{\omega})$ 
which fits in an exact sequence
\begin{equation}
0 \to G \to F \to  \oplus_{i=1}^{r-1}{\Bbb C}_{-x_i} \to 0.
\end{equation}
Applying ${\cal G}_{\cal P}$, we get an exact sequence
\begin{equation}
0 \to {\cal G}_{\cal P}^2(\oplus_{i=1}^{r-1}{\Bbb C}_{-x_i}) \to
{\cal G}_{\cal P}^2(F) \to {\cal G}_{\cal P}^2(G) \to 0
\end{equation}
and 
${\cal G}_{\cal P}^2(G)$ belongs to $M_H(1+H+a \omega)$.
We set $I_Z(H):={\cal G}_{\cal P}^2(G)$.
Since ${\cal G}_{\cal P}^2(\oplus_{i=1}^{r-1}{\Bbb C}_{-x_i}) \cong
\oplus_{i=1}^{r-1}{\cal P}_{x_i}$,
we get our claim.
\end{proof}

\begin{lem}\label{lem:small}
Assume that $a \geq 0$.
If there is an element $E \in M_H(r+dH+(a+b) \omega)$, $b>0$ such that
$H^1(X,E)=0$, then
there is an element $F \in M_H(r+dH+a\omega)$ such that
$H^1(X,F)=0$.
\end{lem}

\begin{proof}
We note that $\dim H^0(X,E)=a+b$.
For general points $x_1, x_2,\dots,x_b \in X$ and a general surjection
$f:E \to \oplus_{i=1}^b {\Bbb C}_{x_i}$,
$H^0(X,E) \to H^0(X,\oplus_{i=1}^b {\Bbb C}_{x_i})$ is surjective.
Hence $H^1(X,\ker f)=0$.
Since $\ker f$ is $\mu$-semi-stable,
Remark \ref{rem:stable} implies that 
there is an element $F \in M_H(r+dH+a\omega)$ of
$H^1(X,F)=0$.
\end{proof}

\begin{thm}\label{thm:d equiv 1} 
Assume that $r>1$, $d \equiv 1 \mod r$, $d>0$ and $a>0$. 
Then for a general element of $M_H(r+dH+a \omega)$,
$\WIT_2$ folds with respect to ${\cal G}_{\cal P}$.
In particular we have a birational map
$M_H(r+dH+a \omega) \cdots \to 
M_{\widehat{H}}(a+d \widehat{H}+r \widehat{\omega})$.
\end{thm}
 
\begin{proof}
We note that $\Hom(E,{\cal P}_x)=0$ for all $x \in \widehat{X}$ and
$E \in M_H(r+dH+a \omega)$.
Hence ${\cal G}_{\cal P}^1(E)$ is torsion free.
If $\Ext^1(E,{\cal P}_x) \cong
H^1(X,E \otimes {\cal P}_{-x})^{\vee}=0$ for some $x \in \widehat{X}$,
then ${\cal G}_{\cal P}^1(E)=0$.
Thus $\WIT_1$ holds for a general $E$. 
So we shall show that $H^1(X,E \otimes {\cal P}_{-x})^{\vee}=0$
for some $x \in \widehat{X}$.
By Lemma \ref{lem:small}, it is sufficient to show our claim
under the assumption $2r>d^2(H^2)-2ra \geq 0$.
By Lemma \ref{lem:rk=1}, our claim follows from Lemma \ref{lem:special}. 
\end{proof}

We next treat the case where $d \equiv -1 \mod r$.

\begin{lem}\label{lem:special2}
Let $E$ be a stable sheaf in Lemma \ref{lem:special}.
Then $A:=E(-H)$ satisfies that
$H^1(X,A^{\vee}(nH) \otimes {\cal P}_x)=
\Ext^1(A(-nH),{\cal P}_x)=0$, $n \geq 0$ for a general $x
\in \widehat{X}$.
\end{lem}

\begin{proof}
We note that there is a curve $C$ of $C \equiv H$
which does not meet $Z$.
We first show that $\Ext^2(A_{|C},{\cal P}_x)=
H^0(X,A_{|C} \otimes {\cal P}_{x}^{\vee})^{\vee}=0$ 
for a general $x \in \widehat{X}$.
Since $\chi({\cal O}_C)=-(H^2)/2<0$,
for a general $x \in \widehat{X}$, we have
$H^0(X,{\cal O}_C \otimes {\cal P}_x^{\vee})=0$.
Hence $H^0(X,I_{Z|C}\otimes {\cal P}_x^{\vee})=0$ 
for a general $x\in \widehat{X}$.
Therefore $H^0(X,A_{|C}\otimes {\cal P}_x^{\vee})=0$ 
for a general $x\in \widehat{X}$.
If $n>0$, then 
we have $H^0(X,A(-nH)_{|C}\otimes {\cal P}_x^{\vee})=0$ 
for all $x \in \widehat{X}$.
We consider an exact sequence induced by \eqref{eq:special}:
\begin{equation}
 0 \to {\cal O}_X \to A^{\vee} 
 \overset{f}{\to} \oplus_{i=1}^{r-1}{\cal P}_{x_i}^{\vee}(H).
\end{equation}
Let $J$ be the image of $f$.
Then we have a filtration 
$0 \subset J_1 \subset J_2 \subset \dots \subset J_{r-1}=J$
such that $J_i/J_{i-1}=I_{Z_i} \otimes {\cal P}_{x_i}^{\vee}(H)$.
Since $\sum_i \deg(Z_i)=\deg Z \leq (H^2)/2$,
$\chi(I_{Z_i} \otimes {\cal P}_{x_i}^{\vee}(H))=
 (H^2)/2-\deg(Z_i) \geq 0$.
Hence Proposition \ref{lem:fourier4} implies that
$H^1(X,J_i/J_{i-1} \otimes {\cal P}_x)=0$ for a general $x \in \widehat{X}$.
Therefore $H^1(X,J \otimes {\cal P}_x)=0$ for a general $x \in \widehat{X}$,
which implies that
$H^1(X,A^{\vee}\otimes {\cal P}_x)=0$ for a general $x \in \widehat{X}$.
Since $\Ext^2(A(-(n-1)H)_{|C}, {\cal P}_x)=0$ for $n \geq 1$ and 
a general $x\in \widehat{X}$,
by using an exact sequence
\begin{equation}
 0 \to A(-nC) \to A(-(n-1)C) \to A(-(n-1)C)_{|C} \to 0
\end{equation}
we get that
$\Ext^1(A(-nC),{\cal P}_x)=0$ for 
a general $x \in \widehat{X}$.
\end{proof}

By Lemma \ref{lem:special2},
we get the following theorem.

\begin{thm}\label{thm:d equiv -1} 
Assume that $r>1$, $d \equiv -1 \mod r$, $d>0$ and $a>0$. 
Then for a general element of $M_H(r+dH+a \omega)$,
$\WIT_2$ folds with respect to ${\cal G}_{\cal P}$.
In particular we have a birational map
$M_H(r+dH+a \omega) \cdots \to 
M_{\widehat{H}}(a+d \widehat{H}+r \widehat{\omega})$.
\end{thm}

\begin{proof}
By Lemma \ref{lem:small}, we may assume that
$2r >d^2(H^2)-2ar \geq 0$.
By Lemma \ref{lem:special2},
there is a locally free sheaf $E \in M_H(r+dH+a \omega)$
such that $H^1(X,E)=0$.
Therefore we get our claims.
\end{proof}

In order to complete the proof of Theorem \ref{thm:evidence} (ii),
we need two more results. 

\begin{thm}\label{thm:(r,0,a)}
${\cal G}_{\cal P}$ induices a birational map
$M_H(r+rdH+a \omega) \cdots \to
M_{\widehat{H}}(a+rd\widehat{H}+r \widehat{\omega})$,
if $r, d, a> 0$ and $(r,a)=1$.
\end{thm}

\begin{proof}
Obviously $H^1(X,{\cal O}_X(dH)^{\oplus r})=0$.
Since $a \leq rd^2(H^2)/2$, by Lemma \ref{lem:small},
we get our claim.
\end{proof}

\begin{prop}\label{prop:(4,2,a)}
${\cal G}_{\cal P}$ induces a birational map
$M_H(4+(4n+2)H+a \omega) \cdots \to
M_{\widehat{H}}(a+(4n+2)\widehat{H}+4 \widehat{\omega})$,
if $n \geq 0$ and $a$ is an odd positive integer.
\end{prop}

\begin{proof}
We set $v=4+(4n+2)H+\{(2n+1)^2(H^2)/2\} \omega$.
Then $\langle v^2 \rangle =0$.
Let $E$ be a semi-stable sheaf of $v(E)=v$.
Since $E$ is semi-homogeneous,
$H^1(X,E)=0$.
Since $a \leq (2n+1)^2(H^2)/2$, by Lemma \ref{lem:small},
we get our claim.
\end{proof}

Combining all together (Theorem \ref{thm:d equiv 1}, \ref{thm:d equiv -1},
\ref{thm:(r,0,a)} and Proposition \ref{prop:(4,2,a)}), 
we get Theorem \ref{thm:evidence} (ii).


\section{Appendix: 
Isomorphisms of moduli spaces induced by more general functor
${\cal F}_{\cal E}$}

In this section, we treat more general cases than
section 3.1.
Let $(X,H)$ be a polarized abelian (or K3) surface of
$(H^2)=2r_0k$, where $r_0$ and $k$ are positive integers of
$(r_0,k)=1$.
We assume that $\NS(X)={\Bbb Z}H$.
We set $v_0:=r_0+d_0 c_1(H)+d_0^2k \omega_X$,
where $d_0$ is an integer of $(r_0,d_0)=1$ and
$\omega_X$ is the fundamental class of $X$.
Then $\langle v_0^2 \rangle =0$.
So $Y:=M_H(v_0)$ is an abelian (or K3) surface.
Since $X$ and $Y$ are isogenous, $\NS(Y) \cong {\Bbb Z}$. 
Since $(r_0,d_0^2k)=1$,
there is a universal family ${\cal E}$ on $X \times Y$.
We assume that 
\begin{itemize}
\item
${\cal E}$ is locally free.\footnote{
If ${\cal E}$ is not locally free,
then Fourier-Mukai functor is the same as reflection
functor \cite{Mu:4}. This case was treated in \cite{Mr:1} and \cite{Y:5}.}
\end{itemize}
Bridgeland \cite{Br:2}, Mukai \cite{Mu:7}
(see Corollary \ref{cor:FM-const}),\cite{Mu:8} 
and Orlov \cite{O:1} proved that
${\cal E}$ satisfies conditions \eqref{eq:str-simple1} and 
\eqref{eq:str-simple2}.
Thus ${\cal E}$ defines Fourier-Mukai functor 
${\cal F}_{\cal E}:{\mathbf D}(X) \to {\mathbf D}(Y)$. 
Let ${\cal F}_{\cal E}^H:H^{ev}(X,{\Bbb Z}) \to H^{ev}(Y,{\Bbb Z})$
be the induced isometry.
We set
\begin{equation}
 \widehat{H}:=\det(p_{Y!}({\cal E} \otimes {\cal O}_H(kr_0-2kd_0)))^{\vee}.   
\end{equation}
We claim that $\widehat{H}$ is a primitive ample line bundle of 
$(c_1(\widehat{H})^2)=(c_1(H)^2)$.
\newline
Proof of the claim:
By Grothendieck-Riemann-Roch theorem,
we have 
\begin{equation}
 c_1(\widehat{H})=\theta_{v_0}(v({\cal O}_H(kr_0-2kd_0))^{\vee}).
\end{equation}
By direct computations, we see that 
$v_0$ and $v({\cal O}_H(kr_0-2kd_0))^{\vee}=-(c_1(H)+2kd_0 \omega_X)$
generate $(v_0)^{\perp} \cap ({\Bbb Z} \oplus \NS(X) 
\oplus {\Bbb Z} \omega_X)$.
Since 
\begin{equation}\label{eq:00/08/21}
 \theta_{v_0}:(v_0)^{\perp}/{\Bbb Z}v_0
 \to H^2(Y,{\Bbb Z})
\end{equation}
is an isomorphism,
$c_1(\widehat{H})$ is primitive.
By Li \cite{Li:1} (or \cite{BBH:2}) and the following lemma, $\widehat{H}$ is ample.
Hence $\widehat{H}$ is a primitive ample line bundle.
\begin{lem}\label{rem:00/08/21}
Let $X$ be an abelian surface and
$w=r+c_1+a \omega_X$, $c_1 \in \NS(X)$ 
a primitive isotropic Mukai vector of $r>0$.
Let $L$ be an ample line bundle such that $r|(c_1(L),c_1)$.
Let $D$ be a divisor on $C \in |L|$ of degree $(C^2)/2-(c_1,C)/r$.
Then 
\begin{equation}
\widehat{L}:=\det(p_{Y!}({\cal E} \otimes {\cal O}_C(D)))^{\vee}  
\end{equation}
is an ample line bundle on $Y$.
\end{lem}

\begin{proof}
Let $E$ be a quasi-homogeneous vector bundle of $v(E)=w$.
Mukai \cite{Mu:1} showed that 
$E$ is stable with respect to any ample line bundle on $X$.
Let $m:X \times X \to X$ be the multiplication map.
Then $m^* E$ is a family of quasi-homogeneous vector bundles
of Mukai vector $w$. 
Hence we get a morphism $\phi_E:X \to M_H(w)$.
Replacing $E$ by $E \otimes L^{\otimes n}$, $n \gg 0$,
we may assume that $c_1(E)$ is ample.
Then this map is finite.
In order to prove the ampleness of 
$\widehat{L}$, it is sufficient to prove the ampleness of 
$\phi_E^*(\widehat{L})=
\det(p_{2!}(m^* E \otimes p_{1}^*({\cal O}_C(D))))^{\vee}$,
where $p_i:X \times X \to X$ is the $i$-th projection. 
By using Grothendieck-Riemann-Roch theorem, we have
$c_1(\phi_E^*(\widehat{L}))=-[p_{2*}(m^*(w)p_{1}^*(C-\frac{(c_1,C)}{r} \omega_X))]_1$.
A simple calculation shows that
$c_1(\phi_E^*(\widehat{L}))=\frac{(c_1,C)}{r}c_1-\frac{(c_1^2)}{2r}C$.
Since $X$ is an abelian surface,
$(c_1(\phi_E^*(\widehat{L}))^2)=a^2(c_1(L)^2)>0$
and $(c_1(\phi_E^*(\widehat{L})),c_1)=(c_1(L),c_1)(c_1^2)/2r>0$ imply that
$\phi_E^*(\widehat{L})$ is ample.
\end{proof}
 
Since $-v({\cal F}_{\cal E}({\cal O}_H(kr_0-2kd_0)))=
c_1(\widehat{H})+\lambda \omega_Y, \lambda \in {\Bbb Z}$
and ${\cal F}_{\cal E}^H$ is an isometry of
Mukai lattice,
\begin{equation}
 \begin{split}
  (c_1(\widehat{H})^2)&=
  \langle v({\cal F}_{\cal E}({\cal O}_H(kr_0-2kd_0)))^2 \rangle\\
  &=\langle v({\cal O}_H(kr_0-2kd_0))^2 \rangle
  =(c_1(H)^2 ).
 \end{split}
\end{equation}
\begin{lem}\label{lem:image}
Let $d_1$ and $l$ be integers which satisfy
$d_1(kd_0)-lr_0=1$.
Then replacing ${\cal E}$ by ${\cal E}\otimes p_Y^* N$, $N \in \Pic(Y)$,
we get 
\begin{equation}
 \begin{cases}
  {\cal F}_{\cal E}^H(1)=d_0^2 k+d_0 l c_1(\widehat{H})+l^2 r_0 \omega_Y\\
  {\cal F}_{\cal E}^H(c_1(H))=2d_0 k r_0+(2d_0 kd_1-1)c_1(\widehat{H})+
  (2d_0 k^2 d_1^2-2d_1 k) \omega_Y\\
  {\cal F}_{\cal E}^H(\omega_X)=r_0+d_1 c_1(\widehat{H})+d_1^2 k \omega_Y.
 \end{cases}
\end{equation}
\end{lem}
\begin{proof}
We set 
\begin{equation}
  \begin{cases}
    c_1({\cal F}^H_{\cal E}(1)) =a c_1(\widehat{H}),\\
    c_1({{\cal F}^H_{\cal E}(c_1(H))}) =b c_1(\widehat{H}),\\
    c_1({\cal F}^H_{\cal E}(\omega_X)) =c c_1(\widehat{H}),
  \end{cases}
\end{equation}
where $a,b,c \in {\Bbb Z}$.
Since $\widehat{{\cal F}}_{\cal E}({\Bbb C}_y)={\cal E}_y^{\vee}$ 
for $y \in Y$,
we have $\widehat{{\cal F}}^H_{\cal E}(\omega_Y)=v_0^{\vee}$.
This implies that ${\cal F}_{\cal E}^H(v_0^{\vee})=\omega_Y$. 
Hence we get the relation
\begin{equation}
 r_0a-d_0 b+d_0^2kc=0.
\end{equation}
Since $(r_0,d_0)=1$,
$b \equiv d_0kc \mod r_0$.
By the definition of $\widehat{H}$,
$-b+2kd_0c=1$.
Hence we get $kd_0c \equiv 1 \mod r_0$.
By the definition of $d_1$, we get $c \equiv d_1 \mod r_0$.
Replacing ${\cal E}$ by
${\cal E} \otimes \widehat{H}^{\otimes((d_1-c)/r_0)}$,
we may assume that
\begin{equation}
  \begin{cases}
    c_1({\cal F}^H_{\cal E}(c_1(H)))=(2kd_0d_1-1) c_1(\widehat{H}),\\
    c_1({\cal F}^H_{\cal E}(\omega_X))=d_1 c_1(\widehat{H}).
  \end{cases}
\end{equation}
Since ${\cal F}^H_{\cal E}$ is an isometry,
$\langle {\cal F}^H_{\cal E}(\omega_X)^2 \rangle=
\langle \omega_X^2 \rangle=0$.
Hence we get
\begin{equation}
{\cal F}^H_{\cal E}(\omega_X)=r_0+d_1 c_1(\widehat{H})+d_1^2k \omega_Y.
\end{equation}
Since ${\cal E}$ is a universal family of stable sheaves of
Mukai vector $v_0$,
we get the following relations:
\begin{equation}
 \begin{cases}
  {\cal F}^H_{\cal E}(v_0^{\vee})=\omega_Y\\
  {\cal F}^H_{\cal E}(\omega_X)=r_0+d_1 c_1(\widehat{H})+d_1^2k \omega_Y\\
  {\cal F}^H_{\cal E}(-c_1(H)+2kd_0 \omega_X)=x+c_1(\widehat{H})+y \omega_Y,
 \end{cases}
\end{equation}
where $x,y \in {\Bbb Z}$.
Since ${\cal F}^H_{\cal E}$ is an isometry,
we see that $x=0$ and $y=2kd_1$.
Hence we get our lemma.
\end{proof}
\begin{lem}
We set $w_0:=r_0+d_1 c_1(\widehat{H})+d_1^2 k \omega_Y$.
Then $X \cong M_{\widehat{H}}(w_0)$ and 
${\cal E}$ is a universal family on 
$M_{\widehat{H}}(w_0) \times Y$.
In particular, ${\cal E}_{|\{x \} \times Y}$, $x \in X$
 is a $\mu$-stable vector bundle.
\end{lem}
\begin{proof}
By \eqref{eq:str-simple1}, 
${\cal E}_{|\{x \} \times Y}$, $x \in X$ is a simple vector bundle on $Y$.
Since $\NS(Y)={\Bbb Z}\widehat{H}$ and $w_0=v({\cal E}_{|\{x \} \times Y})$ 
is an isotropic Mukai vector, 
${\cal E}_{|\{x \} \times Y}$
is a stable vector bundle on $Y$ (\cite[Prop. 3.14]{Mu:4}).
Since $(r,d_1)=1$, ${\cal E}_{|\{x \} \times Y}$
is $\mu$-stable.
Hence there is a morphism $\psi:X \to M_{\widehat{H}}(w_0)$.
By \eqref{eq:str-simple1}, this morphism is injective.
Since $M_{\widehat{H}}(w_0)$ is a smooth projective surface,
$\psi$ is an isomorphism.
\end{proof}
In order to generalize Proposition \ref{lem:fourier4},
and \ref{lem:fourier3}, let us introduce some notations.
Let $G$ be a locally free sheaf on $X$.
For a torsion free sheaf $E$ on $X$, we define
\begin{equation}
\begin{split}
\rk_G(E):& =\rk(E \otimes G^{\vee}),\\
\deg_G(E):&=\deg (E \otimes G^{\vee}),\\
\mu_G(E):& =\frac{\deg_G(E)}{\rk_G(E)}.
\end{split}
\end{equation}
For $x \in {\mathbf D}(X)$ (resp. $v(x) \in H^{ev}(X,{\Bbb Z})$),
we can also define $\rk_G(x)$ and $\deg_G(x)$
(resp. $\rk_G(v(x))$ and $\deg_G(v(x))$).
Then we see that 
\begin{equation}
\begin{split}
\mu_G(E) & =\frac{\deg(E) \rk (G)-\deg(G) \rk(E)}{\rk(G)\rk(E)}\\
&=\mu(E)-\mu(G).
\end{split}
\end{equation}  
The following lemma is obvious.
\begin{lem}
\begin{enumerate}
\item[(1)]
$E$ is $\mu$-stable if and only if
\begin{equation}
\mu_G(F)<\mu_G(E)
\end{equation}
for any subsheaf $F \subset E$ of $\rk(F)<\rk(E)$.
\item[(2)]
Let $E,F$ be $\mu$-semi-stable sheaves of $\mu_G(E)>\mu_G(F)$.
Then $\Hom(E,F)=0$.
\end{enumerate}
\end{lem}
Assume that $\deg_G(E)=1$. Then it is easy to see that
$E$ is $\mu$-stable if and only if
$\deg_G(F) \leq 0$ for any subsheaf $F \subset E$ of $\rk(F)<\rk(E)$.  
\begin{lem}\label{lem:deg}
We choose points $s \in X$ and $t \in  Y$. 
We set 
\begin{equation}
\begin{split}
G_1 &:={\cal E}_{|X \times \{t \}}^{\vee},\\
G_2 &:={\cal E}_{|\{s\} \times Y}.
\end{split}
\end{equation}
Then for a Mukai vector $v$,
\begin{equation}
\deg_{G_1}(v)=-\deg_{G_2}({\cal F}_{\cal E}^H(v))=
\deg_{G_2^{\vee}}({\cal F}_{\cal E}^H(v)^{\vee}).
\end{equation}
\end{lem}
\begin{proof}
We set 
\begin{equation}
\begin{split}
v&=r+d c_1(H)+a \omega_X,\\
{\cal F}_{\cal E}^H(v)&=r'+d' c_1(\widehat{H})+a'\omega_Y.
\end{split}
\end{equation}
It is sufficient to prove that
$r'd_1-d'r_0=dr_0+rd_0$.
This follows from the following relations which come from
Lemma \ref{lem:image}: 
\begin{equation}
 \begin{cases}
  r'=r(d_0^2 k)+d(2d_0 r_0 k)+a r_0,\\
  d'=r(d_0 l)+d(2d_0 d_1 k-1)+ad_1.
 \end{cases}
\end{equation}
\end{proof}
Due to this lemma, we can use the same 
arguments as in Propositions \ref{lem:fourier4}
and \ref{lem:fourier3}.
Hence we get the following theorem.
\begin{thm}\label{thm:fourier}
Keep the notations as above.
Let $v:=r+d c_1(H)+a \omega_X$ be a Mukai vector
of $d r_0+r d_0=1$.
\begin{enumerate}
\item
If $-\langle v,v_0^{\vee} \rangle>0$, then
${\cal G}_{\cal E}$ induces an isomorphism
$M_H(v) \to M_{\widehat{H}}({\cal F}_{\cal E}^H(v)^{\vee})$.
\item
If $\langle v,v_0^{\vee} \rangle>0$, then ${\cal F}_{\cal E}$
induces an isomorphism
$M_H(v) \to M_{\widehat{H}}(-{\cal F}_{\cal E}^H(v))$.
\end{enumerate}
\end{thm}
\begin{ex}\label{ex:1}
Keep the notations as above.
We set 
\begin{equation}
 \begin{cases}
  d_0=-(r_0-1),\\
  r=d=1,\\
  k=-n+s r_0,\\
  a=(r_0^2-1)s-r_0n,
 \end{cases}
\end{equation}
where $s>0$, $s r_0>n > 0$ and $(r_0,n)=1$.
Then $\langle v^2 \rangle=2s$ and
$\langle v,v_0^{\vee} \rangle=n>0$.
Applying Theorem \ref{thm:fourier},
we get an isomorphism
$M_H(v) \to M_{\widehat{H}}(-{\cal F}^H_{\cal E}(v))$.
In particular,
we get a more direct proof of \cite[Thm. 0.2]{Y:5}.
\end{ex}

\begin{ex}\label{ex:birat} 
We shall explain an example of two different K3 surfaces
$X,Y$ such that $\Hilb_X^n$ and $\Hilb_Y^n$ are isomorphic or birational
for some $n$. 
Let $(X,H)$ be any polarized K3 surface of $(H^2)=12$.
We set $v_0:=2-H+3 \omega_X$ and
$Y:=M_H(v_0)$.
Then $Y$ is a K3 surface of 
$H^2(Y,{\Bbb Z}) \cong v_0^{\perp}/{\Bbb Z}v_0$.
Under the identification
$-\theta_{v_0}: v_0^{\perp}/{\Bbb Z}v_0 \to H^2(Y,{\Bbb Z})$,
Li \cite{Li:1} implies that
 $\widehat{H}:=H-6\omega_X$ is a nef and big divisor and
$\widehat{H}$ is ample if and only if $Y$ consists of 
$\mu$-stable vector bundles.
In general, $(Y,\widehat{H}) \not \cong (X,H)$ (\cite[Rem. 10.13]{Mu:9}).

Proof of the claim:
Let ${\cal M}_{12}$ be the moduli space of polarized K3 surfaces of degree 12
and $\overline{{\cal M}}_{12}$ the moduli space of 
quasi-polarized K3 surfaces of degree 12.
Then the correspondence $(X, H) \mapsto (Y,\widehat{H})$
gives a morphism ${\cal M}_{12} \to \overline{{\cal M}}_{12}$. 
Since $\overline{{\cal M}}_{12}$ is quasi-projective (and hence Hausdorff),
we shall prove that $(X,H) \not \cong (Y,\widehat{H})$
for an elliptic K3 surface $\pi:X \to {\Bbb P}^1$
of $\rho(X)=2$.
We may assume that there is a section $\sigma$ 
of $\pi$ and $H=\sigma+7f$, where $f$ is a fiber of $\pi$.
Assume that $(X,H) \cong (Y,\widehat{H})$.
Then there is an isometry 
$\psi:\NS(X) \to \NS(Y)$ such that
$\psi(H)=\widehat{H}$.
By direct computations, we see that
$v_0^{\perp}$ is generated by
$e_1:=1-3f$, $e_2:=2f-\omega_X$ and $v_0$.
Hence $\NS(Y) ={\Bbb Z}e_1\oplus {\Bbb Z} e_2$,
where $(e_1,e_1)=(e_2,e_2)=0$ and $(e_1,e_2)=1$.
Since $\widehat{H}=H-6\omega_X \equiv 2-3\omega_X \mod$ ${\Bbb Z} v_0$,
$\widehat{H}=2e_1+3e_2$.
We set $\psi(\sigma+f)=ae_1+be_2$ and
$\psi(f)=ce_1+de_2$. 
Since $\psi$ is an isometry, $ab=cd=0$.
By $\psi(H)=\widehat{H}$, we get
$ae_1+be_2+6(ce_1+de_2)=2e_1+3e_2$.
Thus $a+6c=2$ and $b+6d=3$.
Since $ab=0$, this is impossible.
Therefore $(X,H) \not \cong (Y,\widehat{H})$.

We shall choose a polarized K3 surface $(X,H)$ such that
$\rho(X)=1$ and $(X,H) \not \cong (Y,\widehat{H})$.
Since $\rho(X)=\rho(Y)=1$, $X \ne Y$. 
We shall consider Fourier-Mukai functor defined by $v_0$.
We shall choose $d_1=-1$ and $l=1$.
Let ${\cal E}$ be a universal vector bundle
on $X \times Y$. 
Then, by Lemma \ref{lem:image}, we may assume that 
\begin{equation}
 \begin{cases}
  {\cal F}_{\cal E}^H(1)=3-\widehat{H}+2 \omega_Y\\
  {\cal F}_{\cal E}^H(H)=-12+5\widehat{H}-12 \omega_Y\\
  {\cal F}_{\cal E}^H(\omega_X)=2-\widehat{H}+3 \omega_Y.
 \end{cases}
\end{equation}
By ${\cal G}_{\cal E}^H=D \circ {\cal F}_{\cal E}^H$,
we get an isomorphism $M_H(1+H+5\omega_X) \cong 
M_{\widehat{H}}(1+\widehat{H}+5\omega_Y)$.
Thus $\Hilb_X^2 \cong \Hilb_Y^2$.
Since ${\cal F}_{\cal E}^H(1+H+4 \omega_X)=-(1-2\omega_Y)$,
we also get an isomorphism
$\Hilb_X^3 \cong \Hilb_Y^3$.

For $n=4$, we still have a birational map
$\Hilb_X^4 \cdots \to \Hilb_Y^4$.
In this case, $\Hilb_X^4 \not \cong \Hilb_Y^4$.
We first construct a birational map $\Hilb_X^4 \cdots \to \Hilb_Y^4$ and
next we show that $\Hilb_X^4 \not \cong \Hilb_Y^4$.
\newline
(1) Construction of birational map:
By Theorem \ref{thm:fourier},
${\cal F}_{\cal E}^H$ induces an isomorphism
$M_H(1+H+3 \omega_X) \cong M_{\widehat{H}}(3-\widehat{H}+\omega_Y)$.
We set $v=1+H+3 \omega_X$ and $w:=-{\cal F}_{\cal E}^H(v)$.
By taking the dual of $E \in M_{\widehat{H}}(w)$, we get a birational map
$M_{\widehat{H}}(w) \cdots \to M_{\widehat{H}}(w^{\vee})$.
We set $\hat{v}:=1+\widehat{H}+3 \omega_Y$.
Since $\langle (w^{\vee}+v({\cal O}_X))^2 \rangle<-2$,
\cite[Rem. 2.3]{Y:5} implies that $H^1(Y,E)=0$ for all $E \in M_{\widehat{H}}(w)$.
Then ${\cal G}_{I_{\Delta}}$ induces an isomorphism
$M_{\widehat{H}}(w^{\vee}) \to M_{\widehat{H}}(\hat{v})$,
where $I_{\Delta}$ is the ideal sheaf of the diagonal
$\Delta \subset Y \times Y$ (this is due to Markman \cite{Mr:1}).
Hence we get a desired birational map
\begin{equation}
 M_H(v) \to M_{\widehat{H}}(w) \cdots \to M_{\widehat{H}}(w^{\vee}) 
 \to M_{\widehat{H}}(\hat{v}).
\end{equation}
(2) We shall next prove that $M_H(v) \not \cong M_{\widehat{H}}(\hat{v})$.
We first describe the ample cone of $M_H(v)$ and 
$M_{\widehat{H}}(\hat{v})$.
We claim that
\begin{enumerate}
\item[(a)]
\begin{equation}
 {\Amp}(M_H(v))=
  {\Bbb Q}_{> 0}\theta_v(-(H+12 \omega_X))+{\Bbb Q}_{> 0}\theta_v(4+H).
\end{equation}
\item[(b)]
${\Bbb Q}_{> 0}\theta_v(-(H+12 \omega_X))$ gives the
Hilbert-Chow morphism $M_H(v) \to S^4X$.
\item[(c)]
By the identification $M_H(v) \cong M_{\widehat{H}}(-{\cal F}_{\cal E}^H(v))$,
${\Bbb Q}_{> 0} \theta_v(4+H)$ corresponds to
the morphism $M_{\widehat{H}}(-{\cal F}_{\cal E}^H(v)) \to 
N_{\widehat{H}}(-{\cal F}_{\cal E}^H(v))$, where 
$N_{\widehat{H}}(-{\cal F}_{\cal E}^H(v))$ is the Uhlenbeck compactification of
the moduli space of $\mu$-stable bundles. 
\end{enumerate}
Proof of (a), (b), (c):
Let $Z$ be the closed subscheme of $M_H(v)$
such that 
${\cal F}_{\cal E}(Z)$ is the set of non-locally free sheaves.
We shall show that $Z \cong {\cal F}_{\cal E}(Z)$ 
is a ${\Bbb P}^2$-bundle over $X \times Y$.
For simplicity, we only give a set-theoretic description of 
${\cal F}_{\cal E}(Z)$ here.
For more detail, see the argument in section 4.4.
\newline
Description of ${\cal F}_{\cal E}(Z)$:
For $E \in {\cal F}_{\cal E}(Z)$, there is an exact sequence
\begin{equation}
 0 \to E \to E^{\vee \vee} \to {\Bbb C}_y \to 0,
\end{equation}
where $E^{\vee \vee} \in M_{\widehat{H}}(3-\widehat{H}+2 \omega_Y)$
and $y \in Y$.
Since ${\cal G}_{\cal E}$ induces an isomorphism
$X=M_H(1+H+6\omega_X) \cong M_{\widehat{H}}(3+2\widehat{H}+8\omega_Y)$
and
$M_{\widehat{H}}(3+2\widehat{H}+8\omega_Y) \cong 
M_{\widehat{H}}(3-\widehat{H}+2 \omega_Y)$,
${\cal F}_{\cal E}(Z)$ is a ${\Bbb P}^2$-bundle over $X \times Y$.
 
We return to the proof of (a), (b), (c).
By Li \cite{Li:1},
${\Bbb Q}_{> 0} \theta_{w}(-(\widehat{H}-4\omega_Y))$ gives 
the contraction $M_{\widehat{H}}(w) \to 
N_{\widehat{H}}(w)$ 
which contracts all fibers of ${\cal F}_{\cal E}^H(Z) \to X \times Y$.
Hence $\theta_v(4+H)=\theta_w(-(\widehat{H}-4\omega_Y))$ 
(cf. Proposition \ref{prop:comm})
gives a boundary of the ample cone.
Since the other boundary corresponds to the Hilbert-Chow
morphism and it is given by ${\Bbb Q}_{> 0} \theta_v(-(H+12 \omega_X))$, 
we get our claim.

In the same way as above, we see that
\begin{enumerate}
\item[(${\mathrm a}'$)]
\begin{equation}
 {\Amp}(M_{\widehat{H}}(\hat{v}))=
  {\Bbb Q}_{> 0}\theta_v(-(\widehat{H}+12 \omega_Y))
  +{\Bbb Q}_{> 0}\theta_v(4+\widehat{H}).
\end{equation}
\item[(${\mathrm b}'$)]
${\Bbb Q}_{> 0}\theta_v(-(\widehat{H}+12 \omega_Y))$ gives the
Hilbert-Chow morphism $M_{\widehat{H}}(\hat{v}) \to S^4Y$.
\item[(${\mathrm c}'$)]
By the identification $M_{\widehat{H}}(\hat{v}) \cong 
M_{\widehat{H}}(w^{\vee})$,
${\Bbb Q}_{> 0} \theta_{\hat{v}}(4+\widehat{H})$ corresponds to
the morphism $M_{\widehat{H}}(w^{\vee}) \to 
N_{\widehat{H}}(w^{\vee})$.
\end{enumerate}
Now we can show that $M_H(v) \not \cong M_{\widehat{H}}(\hat{v})$:
Assume that $M_H(v) \cong M_{\widehat{H}}(\hat{v})$.
Then the isomorphism induces an isometry 
$\zeta:H^2(M_H(v),{\Bbb Z}) \to H^2(M_{\widehat{H}}(\hat{v}),{\Bbb Z})$
such that $\zeta(E_v)=E_{\hat{v}}$,
where $E_v$ and $E_{\hat{v}}$ are exceptional divisors of
Hilbert-Chow morphisms.
Hence we get an Hodge isometry
 $H^2(X,{\Bbb Z})=E_v^{\perp} \to E_{\hat{v}}^{\perp}=
H^2(Y,{\Bbb Z})$.
Then Torelli's theorem implies that $X \cong Y$, which is a contradiction. 
Therefore $M_H(v) \not \cong M_{\widehat{H}}(\hat{v})$. 
\qed

As a final remark, we shall consider the relation between 
$M_H(v)$ and $M_{\widehat{H}}(\hat{v})$.
Let $\widehat{Z}$ be the closed subscheme of $M_{\widehat{H}}(\hat{v})$
such that ${\cal G}_{I_{\Delta}}(\widehat{Z})$ is the set of
non-locally free sheaves. 
Then $\widehat{Z}$ is a ${\Bbb P}^2$-bundle over $X \times Y$.
We claim that $M_H(v)$ is the elementary transform of 
$M_{\widehat{H}}(\hat{v})$ along $\widehat{Z}$.
\newline
Proof of the claim:
Let $\widetilde{M_{\widehat{H}}(\hat{v})} \to M_{\widehat{H}}(\hat{v})$ 
be the blow-up of $M_{\widehat{H}}(\hat{v})$ along $\widehat{Z}$ and
$E$ the exceptional divisor.  
Let $\xi:\widetilde{M_{\widehat{H}}(\hat{v})} \to M_{\widehat{H}}(\hat{v})'$ 
be the contraction to another direction, that is, 
$M_{\widehat{H}}(\hat{v})'$ is the elementary transform of 
$M_{\widehat{H}}(\hat{v})$ along $\widehat{Z}$. 
Assume that 
$\theta_{\hat{v}}(24+7 \widehat{H}+12\omega_Y)_{|f}={\cal O}_f(-n)$,
where $f$ is a fiber of $\widehat{Z} \to X \times Y$.
Then $\theta_{\hat{v}}(24+7 \widehat{H}+12\omega_Y)(-nE)_{|g}$
is trivial for all fibers $g$ of $\xi$.
Hence there is a line bundle ${\cal L}'$ on $M_{\widehat{H}}(\hat{v})'$
such that $\xi^* {\cal L}'=
\theta_{\hat{v}}(m(4+\widehat{H})+(24+7 \widehat{H}+12\omega_Y))
(-nE)$.
If we choose a sufficiently large $m$, then
${\cal L}'$ is nef and big.
Hence base point free theorem implies that
there is a morphism $M_{\widehat{H}}(\hat{v})' \to 
\Proj(\oplus_k H^0(M_{\widehat{H}}(\hat{v})',{{\cal L}'}^{\otimes k}))$.
We set
${\cal L}:=\theta_v(m(4+H)-(H+12\omega_X))$.
Since $m$ is sufficiently large, ${\cal L}$ is ample.
By the canonical identification
$\Pic(M_{\widehat{H}}(\hat{v})')=\Pic(M_{\widehat{H}}(\hat{v}))
=\Pic(M_H(v))$,
we can identify ${\cal L}'$ with ${\cal L}$. 
Since $H^0(M_{\widehat{H}}(\hat{v})',{{\cal L}'}^{\otimes k})=
H^0(M_H(v),{\cal L}^{\otimes k})$ and ${\cal L}$ is ample,
we get a morphism $\phi:M_{\widehat{H}}(\hat{v})' \to M_H(v)$
which is an isomorphism. 
\end{ex}

\begin{rem}
 By using Example \ref{ex:1}, we can construct more examples of
such pairs $X, Y$:
In the notations of Example \ref{ex:1},
if $r_0, k>1$, then we can show that
$M_H(v_0) \not \cong X$ for a general $X$.
If we choose $n=\pm 1$, then we have $\Hilb_X^{s+1} \cong 
\Hilb_Y^{s+1}$.
\end{rem}

\section{Appendix: Deformation types of moduli spaces of stable sheaves on K3 surfaces}

In this section, we shall prove the following theorem. 
\begin{thm}\label{thm:K3-deform}
Let $v=l(r+c_1)+a \omega$, $c_1 \in \NS(X)$ 
be a primitive Mukai vector of $l=\ell(v)$.
Assume that $\rk v>0$ or $c_1$ is ample.
Then $M_H(v)$ is deformation 
equivalent to $\Hilb_X^{\langle v^2 \rangle/2+1}$. 
\end{thm}

\subsection{Proof of the theorem}

The proof is very similar to that of Theorem \ref{thm:deform equiv}.
However there may be a $(-2)$-curve on a K3 surface.
Thus for a divisor $D$,
$(D^2)>0$ does not imply ampleness of $D$.  
Hence we need a modification of the proof of Theorem \ref{thm:deform equiv}.

\subsubsection{The case where $\rk v>0$}
We first assume that $\rk v>0$.
By \cite{Mu:4} and \cite{Y:8}, we may assume that $\langle v^2 \rangle>0$.
Let $X$ be a K3 surface and 
$v=l(r+c_1)+a \omega$, $c_1 \in \NS(X)$ a primitive Mukai vector such that
$l=\ell(v)$ and $2ls:=\langle v^2 \rangle>0$.
We take a positive integer $k$ such that
$n:=rk-(c_1^2)/2>0$.
Since $(l,a)=1$, we may assume that  
$b+lr:=a-kl$ and $r$ are relatively prime.
Let $L$ be an even lattice whose intersection matrix
is given by
\begin{equation}
\begin{pmatrix}
0 & 1 & 0\\
1 & 0 & 0\\
0 & 0 & -2n
\end{pmatrix}.
\end{equation} 
By Nikulin \cite{N:1},
it has a primitive embedding into the K3 lattice.
By the surjectivity of period map,
there is a K3 surface $Y$ whose Picard lattice is isometric to
$L$.
Since $L$ contain a hyperbolic lattice,
there is an elliptic fibration $\pi:Y \to {\Bbb P}^1$ which has a section
$\sigma$.
Let $f$ be a fiber of $\pi$.  
If there is a reducible fiber, then we see that
$n=1$. 
Assume that $n>r^2$.
Then every fiber of $\pi$ is irreducible.
Let $\xi_n$ be a divisor such that $(\xi_n)=-2n$ and 
$(\xi_n,f)=(\xi_n,\sigma)=0$.

We set $w:=l(r+(-\xi_n+f))+(b+lr) \omega$.
By Remark \ref{rem:general},
$\overline{M}_{\sigma+kf}(w)=M_{\sigma+kf}(w)$, $k \gg 0$.
Since $\langle w^2 \rangle =2ls$,
by Proposition \ref{prop:deform},
it is sufficient to prove Theorem \ref{thm:K3-deform} for $M_{\sigma+kf}(w)$, 
$k \gg 0$.
By Theorem \ref{thm:FM}, we have an isomorphism
\begin{equation}
 M_{\sigma+kf}(lr+l(-\xi_n+f)+(b+lr)\omega) \cong 
 M_{\sigma+kf}(l(\xi_n+r \sigma)-bf+l(r+1)\omega), k \gg 0.
\end{equation}
Since $D:=l(\xi_n+r \sigma)-bf$ satsifies $(D^2)=2ls>0$
and $H^0(Y,{\cal O}_Y(-D))=0$, Riemann-Roch theorem implies that
$D$ is effective.
Assume that $D$ is ample.
Since $(b,r)=1$, $D$ is a primitive ample divisor.
By \cite[Thm. 0.2]{Y:5} (or Proposition \ref{prop:deform}, Theorem
\ref{thm:fourier} and Example 3.1),
our theorem holds for the case where $\ell(v)=1$.
Thus $M_{\sigma+kf}(D+l(r+1) \omega)$ and hence
$M_{\sigma+kf}(w), k \gg 0$ is deformation equivalent to 
$\Hilb_Y^{\langle w^2 \rangle/2+1}$.

Ampleness of $D$:
We shall prove that
$D$ is ample unless $r|n-1$ and $-b \leq l(n+1)/r+rl$. 
Let $C:=\alpha \sigma+\lambda f+\beta \xi_n$, 
$\alpha,\beta,\lambda \in {\Bbb Z}$ be a $(-2)$-curve in $Y$.  
Since $(\alpha \sigma+\lambda f+\beta \xi_n)^2=-2$,
$\lambda=(n \beta^2+\alpha^2-1)/\alpha$ is an integer.
Since every fiber is irreducible, we get $\alpha >0$. 
By our assumption,
$\langle v^2 \rangle=(D^2)=-2l^2(n+r^2)-2lbr>0$.
Thus $-b>l(n+r^2)/r$.
Then we see that
\begin{equation}
\begin{split}
(D,C)&=
-2nl \beta-b \alpha+\lambda lr-2 \alpha lr\\
&>l(-2n \beta+\lambda r-2 \alpha r+\alpha(n+r^2)/r)\\
&=l(n(-2 \beta+ \beta^2r/ \alpha+\alpha/r)+\alpha r-r/\alpha-2 \alpha r
+r \alpha)\\
&=l(n\alpha/r(\beta r/\alpha-1)^2-r/\alpha)\\
&=l(n(\beta r-\alpha)^2-r^2)/\alpha r>0,
\end{split}
\end{equation} 
unless $\alpha =r \beta$.
If $\alpha =r \beta$, then $\lambda=(n \beta^2+r^2\beta^2-1)/r \beta$.
This implies that $\beta=1$ and $r|n-1$.
Hence if $r \not |n-1$, then $D$ is ample.
Assume that $r|n-1$.
If $-b>l(n+1)/r+rl$, then we see that
$(D,r\sigma+\lambda f+\xi_n)>0$.
Therefore, $D$ is ample. Thus our claim holds.

We shall next treat the remaining case (i.e. $-b \leq l(n+1)/r+rl$).
In this case, $(D,r\sigma+\lambda f+\xi_n) \leq 0$.
Thus $D$ is not ample.
By deforming the pair $(Y,D)$, we shall reduce the problem to 
the case where $D$ is ample.
Let $({\cal Y}, {\cal L}) \to T$ be a family of polarized K3 surfaces
of $({\cal L}^2_t)=2k$ such that the period map of
polarized K3 surfaces is submersive for every point of $T$. 
Replacing $T$ by a suitable covering, we assume that 
$({\cal Y}, {\cal L}) \to T$ has a section.
Let $\zeta$ be a numerical equivalence class of $(\zeta,{\cal L}_t)=(D,\sigma+
kf)$ and
$(\zeta^2)=2sl$.
Let $\varpi:\Pic_{{\cal Y}/T}^{\zeta} \to T$ 
be the relative Picard scheme over $T$ and ${\cal P}$ be the universal family
of line bundles.
For a point $x$ of $\Pic_{{\cal Y}/T}^{\zeta}$, we choose a sufficiently small
open neighbourhood $U$ of $x$.
Then by the description of the period domain,
$U \to  T$ is an immersion of $\codim_T U=1$.
Assume that $\varpi(x)$ corresponds to 
$(Y,\sigma+kf)$ and $x$ corresponds to $D$.  
Then in a neighborhood of $x$,
there is a point $y$ such that $\rho({\cal Y}_{\varpi(y)})=2$
and ${\cal P}_y$ is ample.
Indeed, it is easy to see that
$r \sigma+\xi_n+(r+(n-1)/r)f$, $D$ and $\sigma+kf$
are linearly independent.
Hence if $\rho({\cal Y}_{\varpi(y)})=2$, then 
$r \sigma+\xi_n+(r+(n-1)/r)f$ does not belong to $\Pic({\cal Y}_{\varpi(y)})$.
Thus ${\cal P}_y$ must be ample.
Moreover we may assume that 
$\overline{M}_{{\cal L}_{\varpi(y)}}(c_1({\cal P}_y)+l(r+1) \omega_y)=
M_{{\cal L}_{\varpi(y)}}(c_1({\cal P}_y)+l(r+1) \omega_y)$.

By using the family ${\cal Y} \times_T \Pic_{{\cal Y}/T}^{\zeta} \to
\Pic_{{\cal Y}/T}^{\zeta}$, we see that
$M_{\sigma+kf}(D+l(r+1) \omega)$ is deformation equivalent to
$M_{{\cal L}_{\varpi(y)}}(c_1({\cal P}_y)+l(r+1) \omega_y)$,
where $\omega_y$ is the fundamental class of ${\cal Y}_{\varpi(y)}$.
Since ${\cal P}_y$ is ample and ${\cal P}_y$ is primitive, 
$M_{{\cal L}_{\varpi(y)}}(c_1({\cal P}_y)+l(r+1) \omega_y)$ is deformation 
equivalent to $\Hilb_Y^{ls+1}$.
Therefore our theorem also holds for this case.
\qed

\subsubsection{The case where $\rk v=0$}
We next treat the case where $\rk v=0$ and $c_1(v)$ is ample.
Let $Y$ be the same K3 surface as above.
Then we have an isomorphism
\begin{equation}
 M_{\sigma+kf}(l-l\xi_n+(a-l)f+(1+b)l\omega) \cong 
 M_{\sigma+kf}(l(\xi_n+\sigma-bf)+a\omega), k \gg 0.
\end{equation}
In the same way as above, we see that
$D:=\xi_n+\sigma-bf$ is ample for $-b>n+2$.
If $-b=n+2$, then $(D,C)>0$ except
for the $(-2)$-curve $C=\sigma+\xi_n+nf$.
Hence $D$ and $\sigma+kf$ deform to ample divisors.
We note that $(D^2)=2(-b-1-n)$.
Since we can take an arbitrary positive $-b-1-n$,
our theorem also holds for this case.
\qed

\subsection{A remark on Theorem \ref{thm:K3-deform}}
\subsubsection{A remark on Theorem \ref{thm:K3-deform}}
In the notation of 5.1,
assume that $r |((c_1^2)/2+1)$ and $\langle v^2 \rangle < 2 l^2$
(this case corresponds to $(D,\alpha \sigma+\lambda f+ \beta \xi_n)<0$).
We note that $v_0:=r+c_1+\{((c_1^2)/2+1)/r\}\omega$ is a $(-2)$-vector. 
We set $v:=l v_0-b\omega$.
Since $0<\langle v^2 \rangle < 2 l^2$,
we get $l<b \rk v_0 <2l$.
Then $w:=(b \rk v_0-l)v_0^{\vee}-b \omega$ satisfies
$\langle w^2 \rangle>2\ell(w)^2$.
For a general polarization $H$,
let $E_0$ be the stable vector bundle of $v(E_0)=v_0$.
We set $X_1=X_2=X$ and $p_i:X_1 \times X_2 \to X_i$,
$i=1,2$ the projections.  
We set
\begin{equation}\label{eq:E_0}
{\cal E}:=\ker(ev:E_0 \boxtimes E_0^{\vee} \to {\cal O}_{\Delta}),
\end{equation}
where $\Delta \subset X_1 \times X_2$ is the diagonal and
$ev$ is the evaluation map.
Then ${\cal E}_{|\{x_1 \}\times X_2}$, $x_1 \in X_1$ 
(resp. ${\cal E}_{|X_1\times \{x_2\}}$, $x_2 \in X_2$) is  
a stable sheaf of $v({\cal E}_{|\{x_1 \}\times X_2})=v_0$
(resp. of $v({\cal E}_{|X_1\times \{x_2\}})=v_0^{\vee}$).

We consider a functor
${\cal H}_{\cal E}:{\bf D}(X_1) \to {\bf D}(X_2)_{op}$ 
which is the composition of 
reflection by $v_0$ with the taking dual functor:
\begin{equation}
 {\cal H}_{\cal E}(x):={\bf R}\Hom_{p_2}(p_1^*(x),{\cal E}),
 x \in {\bf D}(X_1).
\end{equation}
Then we have the following theorem.
\begin{thm}\label{thm:00/08/16}
Assume that $l<b \rk v_0 <2l$.
Then ${\cal H}_{\cal E}$ induces an isomorphism
$M_H(lv_0-b \omega) \to M_H((b \rk v_0-l)v_0^{\vee}-b \omega)$
for a general ample divisor $H$.
\end{thm}

\subsubsection{Proof of Theorem \ref{thm:00/08/16}}
Let $X$ be a K3 surface or an abelian surface.
We shall prove a generalization of Theorem \ref{thm:00/08/16}
which is also a generalization of \cite[Cor. 4.5]{Mu:5}.
In order to state our theorem (Theorem \ref{thm:00/08/17}),
we prepare some notations.

Let $v_1:=r_1+c_1+a_1 \omega_X, r_1>0, c_1 \in \NS(X)$ be a 
primitive isotropic Mukai vector on $X$.
For a general ample divisor $H$, we set $Y:=M_H(v_1)$.
Assume that there is a universal family ${\cal E}$ on
$X \times Y$.
We set $w_1:=v({\cal E}_{|X \times \{y \}})=
r_1+\widehat{c}_1+\widehat{a}_1 \omega_Y$, $y \in Y$. 
We consider a functor
${\cal H}_{\cal E}:
{\bf D}(X) \to {\bf D}(Y)_{op}$ defined by
\begin{equation}
 {\cal H}_{\cal E}(x):={\bf R}\Hom_{p_Y}(p_X^*(x),{\cal E}),
 x \in {\bf D}(X),
\end{equation}
where $p_X:X \times Y \to X$ (resp. $p_Y:X \times Y \to Y$) be the projection.
Then ${\cal H}_{\cal E}$ gives an equivalence of categories and
the inverse is given by
\begin{equation}
 \widehat{{\cal H}}_{\cal E}(y):={\bf R}\Hom_{p_X}(p_Y^*(y),{\cal E}),
 y \in {\bf D}(Y)_{op}.
\end{equation}
${\cal H}_{\cal E}$ induces an isometry $H^{ev}(X,{\Bbb Z}) \to 
H^{ev}(Y,{\Bbb Z})$.
We denote it by ${\cal H}_{\cal E}^H$.
We also denote the inverse $({{\cal H}_{\cal E}^H})^{-1}$ by 
$\widehat{\cal H}_{\cal E}^H$. 
\begin{rem}
If ${\cal E}$ is locally free, then
${\cal H}_{\cal E}={\cal G}_{{\cal E}^{\vee}}$.
\end{rem}

We have an isomorphism
$\NS(X) \otimes {\Bbb Q} \to v_1^{\perp} \cap \omega_X^{\perp}$
by sending $D \in \NS(X) \otimes {\Bbb Q}$ to 
$D+\frac{1}{r}(D,c_1) \omega_X \in v_1^{\perp} \cap \omega_X^{\perp}$.
Since ${\cal H}_{\cal E}^H$ induces an isomorphism
$v_1^{\perp} \cap \omega_X^{\perp} \to w_1^{\perp} \cap \omega_Y^{\perp}$,
we have an isomorphism
$\delta:\NS(X) \otimes {\Bbb Q} \to \NS(Y) \otimes {\Bbb Q}$.
This map is nothing but $-\theta_{v_1}$ in \eqref{eq:00/08/21}. 
For a ${\Bbb Q}$-line bundle $L \in \Pic(X) \otimes {\Bbb Q}$,
we choose a ${\Bbb Q}$-line bundle $\widehat{L}$ on $Y$ such that
$\delta(c_1(L))=c_1(\widehat{L})$.
Let $H$ be an ample line bundle on $X$.
Then Lemma \ref{rem:00/08/21} implies that 
$\widehat{H}$ is ample, if $Y$ consists of $\mu$-stable vector bundles.
By the proof of \cite[Lem. 2.1]{Y:8}, 
$Y$ consists of $\mu$-stable vector bundles
unless ${\cal E}$ is given by \eqref{eq:E_0}.
In this case, a direct computation (or \cite{Li:1}) shows that
$\widehat{H}$ is ample.

Assume that 
\begin{itemize}
\item[$(\star)$]
\begin{enumerate}
\item
$H$ is an ample divisor on $X$ 
which is general with respect to $v$ and $\widehat{H}$
is an ample divisor on $Y$ which is general with respect to 
$-{\cal H}_{\cal E}^H(v)$.
\item
${\cal E}_{|\{x \} \times Y}$ is stable with respect to $\widehat{H}$.
\end{enumerate} 
\end{itemize}
\begin{rem}\label{rem:00/08/17}

\begin{itemize}
\item
If $\NS(X)={\Bbb Z}$, then the assumption $(\star)$ holds.
\item
If $X$ is abelian or $Y$ consists of non-locally free sheaves, then
the assumption $(\star)$ holds for a general $H$.
\end{itemize}
\end{rem}
{\bf Problem.} Is ${\cal E}_{|\{x \} \times Y}$ always stable with respect to
$\widehat{H}$?

\vspace{1pc}

For a coherent sheaf $E$ on $X$
(resp. $F$ on $Y$), we set $\deg(E):=(c_1(E),H)$
(resp. $\deg(F):=(c_1(F),\widehat{H})$).
We consider twisted degree 
$\deg_{G_1}(E)$ and $\deg_{G_2}(F)$,
where $G_1:={\cal E}_{|X \times \{y \}}$ and
$G_2:={\cal E}_{|\{x \} \times Y}$.
Then 
\begin{lem}
$\deg_{G_1}(v)=\deg_{G_2}({\cal H}_{\cal E}^H(v))$. 
\end{lem}
\begin{proof}

\begin{equation}\label{eq:deg-preserve}
 \begin{split}
  \deg_{G_1}(v)&=
  ((\rk v_1)c_1(v)-(\rk v)c_1(v_1),H)\\
  &=\langle (\rk v_1)v-(\rk v)v_1,H+\frac{1}{r}(H,c_1) \omega_X \rangle\\
  &=\langle (\rk v_1)v,H+\frac{1}{r}(H,c_1) \omega_X \rangle\\
  &=\langle (\rk w_1){\cal H}_{\cal E}^H(v),\widehat{H}+
   \frac{1}{r}(\widehat{H},\widehat{c}_1)\omega_Y \rangle
   \quad \quad \quad \quad (\because \text{${\cal H}_{\cal E}^H$ 
   is an isometry})\\
  &=\langle (\rk w_1){\cal H}_{\cal E}^H(v)-(\rk{\cal H}_{\cal E}^H(v))w_1,
   \widehat{H}+\frac{1}{r}(\widehat{H},\widehat{c}_1)
   \omega_Y \rangle\\
  &=((\rk w_1)c_1({\cal H}_{\cal E}^H(v))-(\rk{\cal H}_{\cal E}^H(v))c_1(w_1),
   \widehat{H})\\
  &=\deg_{G_2}({\cal H}_{\cal E}^H(v)).
 \end{split}
\end{equation}

\end{proof}

\begin{defn}
For a Mukai vector
$v=l v_1-a \omega_X$, $l, a \in \frac{1}{\ell(v_1)}{\Bbb Z}$,
we set
$l(v):=l$ and $a(v):=a$.
\end{defn}
Since ${\cal H}_{\cal E}^H(v_1)=\omega_Y$ and
$\widehat{{\cal H}}_{\cal E}^H(w_1)=\omega_X$,
we get
\begin{equation}
{\cal H}_{\cal E}^H(l v_1-a \omega_X)=l \omega_Y-a w_1.
\end{equation}
We can now state our theorem.

\begin{thm}\label{thm:00/08/17}
Assume that $l,a>0$.
Then ${\cal H}_{\cal E}$ induces an isomorphism
$M_H(lv_1-a \omega_X) \to M_{\widehat{H}}(aw_1-l \omega_Y)$
for a general ample divisor $H$ on $X$ which satisfies $(\star)$.
\end{thm}
By \ref{rem:00/08/17}, Theorem \ref{thm:00/08/16} follows.
In order to prove Theorem \ref{thm:00/08/17},
we prepare some lemmas.

\begin{lem}\label{lem:G^2-2}
Assume that $a >0$.
Then $\Hom({\cal E}_{|X \times \{y \}},E)=0$
for all $y \in Y$ and $E \in M_H(v)$.
\end{lem}

\begin{proof}
Since ${\cal E}_{|X \times \{y \}}$ 
and $E$ are semi-stable, it is sufficient to show that
$-a({\cal E}_{|X \times \{y \}})/l({\cal E}_{|X \times \{y \}})> -a/l$.
Since $v({\cal E}_{|X \times \{y \}})=v_1$,
$-a({\cal E}_{|X \times \{y \}})/l({\cal E}_{|X \times \{y \}})-( -a/l)=a/l>0$.
\end{proof} 

\begin{lem}\label{lem:G^0-2}
For a $\mu$-semi-stable sheaf $E$ of $v(E)=v$, 
there is a finite subset $S \subset Y$
such that 
$\Hom(E,{\cal E}_{|X \times \{y \}})=0$ for all $y \in X \setminus S$.
\end{lem}

\begin{proof}
Considering Jordan-H\"{o}lder filtration of $E$ with respect to
$\mu$-stability,
we may assume that $E$ is $\mu$-stable.
If ${\cal E}_{|X \times \{y \}}$ is locally free, then 
obviously the claim holds.
Hence we assume that ${\cal E}_{|X \times \{y \}}$ is not locally free.
Under the notation \eqref{eq:E_0},
if $E^{\vee \vee} \ne E_0$, then clearly $\Hom(E,E_0)=0$.
Hence $\Hom(E,{\cal E}_{|X \times \{y \}})=0$ for all $y \in Y$.
If $E^{\vee \vee}=E_0$, then
$\Hom(E,{\cal E}_{|X \times \{y \}})=0$ for 
$y \in Y \setminus \Supp(E^{\vee \vee}/E)$.
\end{proof}
{\it Proof of Theorem \ref{thm:00/08/17}.}
By the symmetry of the condition, it is sufficient to show that
$\WIT_1$ holds for $E \in M_H(v)$ and
${\cal H}^1_{\cal E}(E)$ is stable with respect to $\widehat{H}$.
By Lemma \ref{lem:G^2-2} and
\ref{lem:G^0-2}, $\WIT_1$ holds and 
${\cal H}^1_{\cal E}(E)$ is torsion free.
We shall show that $E$ is semi-stable.

(I) ${\cal H}^1_{\cal E}(E)$ is $\mu$-semi-stable:
Assume that ${\cal H}^1_{\cal E}(E)$ is not $\mu$-semi-stable.
Let $0 \subset F_1 \subset F_2 \subset \dots \subset 
F_s= {\cal H}^1_{\cal E}(E)$
be the Harder-Narasimhan filtration of ${\cal H}^1_{\cal E}(E)$
with respect to $\mu$-semi-stability.
We shall choose the integer $k$ which satisfies 
$\deg_{G_2}(F_i/F_{i-1})\geq 0, i \leq k$ and
$\deg_{G_2}(F_i/F_{i-1}) < 0, i > k$.
We claim that
$\widehat{\cal H}^0_{\cal E}(F_k)=0$ 
and $\widehat{\cal H}^2_{\cal E}({\cal H}_{\cal E}^1(E)/F_k)=0$.
Indeed $\deg_{G_2}(F_i/F_{i-1}) \geq 0, i \leq k$
and the $\mu$-semi-stability of $F_i/F_{i-1}$ 
imply that $\widehat{\cal H}^0_{\cal E}(F_i/F_{i-1}), i \leq k$
is of dimension 0.
Since $\widehat{\cal H}^0_{\cal E}(F_i/F_{i-1})$ is torsion free,
$\widehat{\cal H}^0_{\cal E}(F_i/F_{i-1})=0, i \leq k$.
Hence $\widehat{\cal H}^0_{\cal E}(F_k)=0$.
On the other hand, we also see that
$\widehat{\cal H}^2_{\cal E}(F_i/F_{i-1})=0, i>k$.
Hence we conclude that $\widehat{\cal H}^2_{\cal E}
({\cal H}_{\cal E}^1(E)/F_k)=0$.

So $F_k$ and ${\cal H}_{\cal E}^1(E)/F_k$ satisfy $\WIT_1$ and
we get an exact sequence 
\begin{equation}
0 \to \widehat{\cal H}^1_{\cal E}({\cal H}_{\cal E}^1(E)/F_k) \to E \to 
\widehat{\cal H}^1_{\cal E}(F_k) \to 0.
\end{equation}
By \eqref{eq:deg-preserve},
 $\deg_{G_1}(\widehat{\cal H}^1_{\cal E}(F_k))=-\deg_{G_2}(F_k)<0$.
This means that
$E$ is not $\mu$-semi-stable with respect to $H$.
Therefore ${\cal H}_{\cal E}^1(E)$ is $\mu$-semi-stable with respect to $H$.

(II) ${\cal H}^1_{\cal E}(E)$ is semi-stable:
Assume that ${\cal H}^1_{\cal E}(E)$ is not semi-stable.
Then there is an exact sequence
\begin{equation}
0 \to F_1 \to {\cal H}^1_{\cal E}(E) \to F_2 \to 0
\end{equation}
such that (i) $F_2$ is stable and
(ii) $-a(F_2)/l(F_2)<-a({\cal H}^1_{\cal E}(E))/l({\cal H}^1_{\cal E}(E))
=-l/a$,
where $v(F_2)=l(F_2)w_1-a(F_2)\omega_Y$.
Since $ -a(F_2)/l(F_2)<-l/a<0$,
Lemma \ref{lem:G^2-2} and \ref{lem:G^0-2} imply that
$\widehat{{\cal H}}^0_{\cal E}(F_2)=\widehat{{\cal H}}^2_{\cal E}(F_2)=0$.
We also obtain that
$\widehat{{\cal H}}^0_{\cal E}(F_1)=0$.
Hence we have an exact sequence
\begin{equation}
 0 \to \widehat{{\cal H}}^1_{\cal E}(F_2) \to E \to 
\widehat{{\cal H}}_{\cal E}^1(F_1) \to 0.
\end{equation}
Since $\widehat{\cal H}_{\cal E}^1({\cal H}^1_{\cal E}(E))=E$,
$\widehat{\cal H}^2_{\cal E}(F_1)=0$.
Thus $\WIT_1$ also holds for $F_1$.
By (ii), we see that
\begin{equation}
\begin{split}
\frac{-a}{l}-\frac{-a(\widehat{\cal H}^1_{\cal E}(F_2))}
{l(\widehat{{\cal H}}_{\cal E}^1(F_2))} 
&=
\frac{-a}{l}+\frac{l(F_2)}{a(F_2)}\\
&=\frac{-a a(F_2)+l l(F_2)}{la(F_2)}
<0.
\end{split}
\end{equation}
This means that $E$ is not semi-stable.
Therefore ${\cal H}_{\cal E}^1(E)$ is semi-stable.
\qed

\vspace{1pc}

{\it Acknowledgement.}
I learned Dekker's thesis from G. van der Geer and T. Katsura.
I would like to thank them very much. 
I would also like to thank E. Markman, S. Mukai, K. O'Grady and
M.-H. Saito for
valuable discussions
and Hern\'{a}ndez Ruip\'{e}rez for explaining \cite{H-M:1}.
Most part of this paper was written during my stay 
at Max Planck Institut f\"{u}r Mathematik in 1998-99.
I would like to thank Max Planck Institut f\"{u}r Mathematik
for support and hospitality. 
Finally I would like to thank the referee for valuable suggestions
to improve this paper.


\begin{thebibliography}{[BBH2]}

\bibitem[A]{A:1}
Anghel, C.,
{\it  Fibr\'{e}s vectoriels stables avec $\chi=0$ sur une surface 
ab\'{e}lienne simple,}
Math. Ann. {\bf 315} (1999), 497--501
\bibitem[A-K]{A-K:1}
Altman, A., Kleiman, S.,
{\it Compactifying the Picard scheme,}
Adv. in Math. {\bf 35} (1980), 50--112
\bibitem[BBH1]{BBH:1}
Bartocci, C., Bruzzo, U., Hern\'{a}ndez Ruip\'{e}rez, D.,
{\it A Fourier-Mukai transform for stable bundles on $K3$ surfaces,}
J. Reine Angew. Math. {\bf 486} (1997), 1--16
\bibitem[BBH2]{BBH:2}
Bartocci, C., Bruzzo, U., Hern\'{a}ndez Ruip\'{e}rez, D.,
{\it Existence of $\mu$-stable vector bundles on $K3$ surfaces and 
the Fourier-Mukai transform,}
 Algebraic geometry (Catania, 1993/Barcelona, 1994), 245--257, 
Lecture Notes in Pure and Appl. Math.,
200, Dekker, New York, 1998
\bibitem[B]{B:1}
Beauville, A.,
{\it Vari\'{e}t\'{e}s K\"{a}hleriennes dont la premi\`{e}re classe de Chern
est nulle,}
J. Diff. Geom. {\bf 18} (1983), 755--782
\bibitem[Br1]{Br:1}
Bridgeland, T.,
{\it Fourier-Mukai transforms for elliptic surfaces,}
J. reine angew. Math. {\bf 498} (1998), 115--133 
\bibitem[Br2]{Br:2}
Bridgeland, T.,
{\it Equivalences of triangulated categories and Fourier-Mukai
transforms,} 
Bull. London Math. Soc. {\bf 31} (1999), 25--34,
math.AG/9809114
\bibitem[D]{D:1}
Dekker, M.,
{Moduli spaces of stable sheaves on abelian surfaces,}
Thesis, Universiteit van Amsterdam 
\bibitem[F]{Fu:1}
Fulton, W.,
{\it Intersection Theory,}
Erg Math. (3. Folge) Band 2, Springer Verlag, 1984
\bibitem[F-L]{F-L:1}
Fahlaoui, R.,
Laszlo, Y., 
{\it Transform\'{e}e de Fourier et stabilit\'{e} sur les surfaces 
ab\'{e}liennes,}
Compositio Math. {\bf 79} (1991), 271--278. 
\bibitem[G]{Gi:1}
Gieseker, D.
{\it On the moduli of vector bundles on an algebraic surface,}
Ann. of Math. {\bf 106} (1977), 45--60
\bibitem[G\"{o}]{Go:1}
G\"{o}ttsche, L.,
{\it Hilbert schemes of zero-dimensional subschemes of smooth varieties,}
 Lecture Notes in Mathematics, 1572. Springer-Verlag, Berlin, (1994)
\bibitem[G-H]{G-H:1}
G\"{o}ttsche, L., Huybrechts, D.,
{\it Hodge numbers of moduli spaces of stable bundles on K3 surfaces,}
Internat. J. Math. {\bf 7} (1996), 359--372
\bibitem[H-M]{H-M:1}
Hern\'{a}ndez Ruip\'{e}rez, D., Mu\~{n}oz Porras, J. M.,
{\it Structure of the moduli space of stable sheaves on elliptic fibrations,}
Preprint Universidad de Salamanca, (1998)
\bibitem[Hi]{Hi:1}
Hironaka, H.,
{\it An example of a non-K\"{a}hlerian complex-analytic deformation of
K\"{a}hlerian complex structures,}
 Ann. of Math. {\bf 75} (1962), 190--208 
\bibitem[Hu]{H:2}
Huybrechts, D.,
{\it Compact Hyperk\"{a}hler Manifolds: Basic Results,}
alg-geom/9705025
\bibitem[J-M]{J-M:1}
Jardim, M., Maciocia, M.,
{\it A Fourier-Mukai approach to spectral data for instantons,}
math.AG/0006054 
\bibitem[K]{K:1}
Knutson, D., 
{\it Algebraic Spaces,}
Lecture Notes in Math. {\bf 203}, Springer-Verlag
\bibitem[L-B]{L-B:1}
Lange, H., Birkenhake, Ch.,
{\it Complex Abelian Varieties,}
Springer-Verlag
\bibitem[Li]{Li:1}
Li, J.,
{\it Compactification of moduli of vector bundles over algebraic
surfaces,}
 Collection of papers on geometry, analysis and mathematical physics, 
World Sci. Publishing, River Edge, NJ, (1997), 98--113
\bibitem[Ma1]{Ma:1}
 Maruyama, M.,
 {\it Moduli of stable sheaves II,}
 J. Math. Kyoto Univ. {\bf 18} (1978),  557--614
\bibitem[Ma2]{Ma:3}
 Maruyama, M.,
{\it Moduli of algebraic vector bundles,}
in preparation
\bibitem[Mr]{Mr:1}
Markman, E.,
{\it Brill-Noether duality for moduli spaces of sheaves
on K3 surfaces,}
math.AG/9901072 
\bibitem[Mu1]{Mu:1}
Mukai, S.,
{\it Semi-homogeneous vector bundles on an Abelian variety,}
 J. Math. Kyoto Univ. {\bf 18}
(1978), 239--272
\bibitem[Mu2]{Mu:2}
Mukai, S.,
{\it Duality between $D(X)$ and $D(\hat{X})$ with its application
to Picard sheaves,}
Nagoya Math. J., {\bf 81} (1981), 153--175
\bibitem[Mu3]{Mu:3}
Mukai, S.,
{\it Symplectic structure of the moduli space of sheaves on an 
abelian or K3 surface,}
Invent. math. {\bf 77}
(1984), 101--116
\bibitem[Mu4]{Mu:4}
Mukai, S.,
{\it On the moduli space of bundles on K3 surfaces I,}
Vector bundles on Algebraic Varieties, Oxford, 1987, 341--413 
\bibitem[Mu5]{Mu:5}
Mukai, S.,
{\it Fourier functor and its application to the moduli of bundles
on an Abelian variety,}
Adv. Studies in Pure Math. {\bf 10} 
(1987), 515--550
\bibitem[Mu6]{Mu:6}
Mukai, S.,
{\it Moduli of vector bundles on K3 surfaces, and symplectic manifolds,}
Sugaku Expositions, {\bf 1} (1988), 139--174 
\bibitem[Mu7]{Mu:7}
Mukai, S.,
{\it Abelian variety and spin representation (in Japanese),}
Proceedings of symposium ``Hodge theory and algebraic geometry
(Sapporo, 1994)'', 110--135:
English translation, Univ. of Warwick preprint, 1998
\bibitem[Mu8]{Mu:8}
Mukai, S.,
{\it Duality of polarized K3 surfaces,}
Proceedings of Euroconference of Algebraic Geometry, 1996
to appear 
\bibitem[Mu9]{Mu:9}
Mukai, S.,
{\it Non-Abelian Brill-Noether theory and Fano 3-folds,}
alg-geom/9704015
\bibitem[N]{N:1}
Nikulin, V. V.,
{\it Integral symmetric bilinear forms and some of their
applications,}
English translation,
Math. USSR Izvestija, {\bf 14} (1980), 103--167
\bibitem[O]{O:1}
O'Grady, K.,
{\it The weight-two Hodge structure of moduli spaces of sheaves
on a K3 surface,}
J. Algebraic Geom., {\bf 6} (1997), no. 4, 599--644
\bibitem[Or]{Or:1}
Orlov, D.,
{\it Equivalences of derived categories and K3 surfaces,}
J. Math. Sci. (NY), {\bf 84} (1997) 1361--1381
\bibitem[S]{S:1}
Simpson, C.,
{\it Moduli of representations of the fundamental group
of a smooth projective variety I,}
Publ. Math. I.H.E.S. {\bf 79} (1994), 47--129
\bibitem[Y1]{Y:1}
Yoshioka, K.,
{\em The Betti numbers of the moduli space of stable sheaves 
of rank 2 on $P^2$,}
J. reine angew. Math. {\bf 453} (1994), 193--220
\bibitem[Y2]{Y:2}
Yoshioka, K., 
{\it Chamber structure of polarizations and the moduli of stable sheaves on a ruled surface,}
Internat. J. Math. {\bf 7} (1996), 411--431
\bibitem[Y3]{Y:3}
Yoshioka, K.,
{\it Some notes on the moduli of stable sheaves on elliptic surfaces,}
Nagoya Math. J. {\bf 154} (1999), 73--102
\bibitem[Y4]{Y:4}
Yoshioka, K.,
{\it A note on the universal family of moduli of stable sheaves,}
J. reine angew. Math. {\bf 496} (1998), 149--161
\bibitem[Y5]{Y:5}
Yoshioka, K.,
{\it Some examples of Mukai's reflections on K3 surfaces,}
J. reine angew. Math. {\bf 515} (1999), 97--123     
\bibitem[Y6]{Y:6}
Yoshioka, K.,
{\it Albanese map of moduli of stable sheaves on abelian surfaces,}
math.AG/9901013,
{\it Some examples of isomorphisms induced by Fourier-Mukai functors,}
math.AG/9902105
\bibitem[Y7]{Y:8}
Yoshioka, K.,
{\it Irreducibility of moduli spaces of vector bundles on K3 surfaces,}
math.AG/9907001
\end{thebibliography}
\end{document}